\documentclass[paper=a4,12pt]{scrreprt}
\usepackage[utf8]{inputenc}
\usepackage{dsfont}
\usepackage{graphicx,amsmath,amsfonts,amssymb,amsthm,bbm,mathrsfs}

\usepackage{fancyhdr}
\makeatletter
\renewcommand*\l@chapter[2]{%
  \ifnum \c@tocdepth >\m@ne
    \addpenalty{-\@highpenalty}%
    \vskip 1.0em \@plus\p@
    \setlength\@tempdima{3.0em}%
    \begingroup
      \parindent \z@ \rightskip \@pnumwidth
      \parfillskip -\@pnumwidth
      \leavevmode \sectfont
      \advance\leftskip\@tempdima
      \hskip -\leftskip
      #1\nobreak\hfil \nobreak\hb@xt@\@pnumwidth{\hss #2}\par
      \penalty\@highpenalty
    \endgroup
  \fi}
\renewcommand*\l@section{\@dottedtocline{1}{0em}{3.0em}}
\renewcommand*\l@subsection{\@dottedtocline{2}{3.0em}{3.0em}}
\renewcommand*\l@subsubsection{\@dottedtocline{3}{5.0em}{3.0em}}
\renewcommand*\l@paragraph{\@dottedtocline{4}{7.5em}{3.0em}}
\renewcommand*\l@subparagraph{\@dottedtocline{5}{10.0em}{3.0em}}
\makeatother

\renewcommand{\thesection}
{\Roman{section}}
\renewcommand{\theequation}%
{\thesection.\arabic{equation}}
\newcommand{\secct}[1]{\section{#1}
\setcounter{equation}{0}}

\renewcommand{\epsilon}{\varepsilon}
\newcommand{\eps}{{\varepsilon}}
\newcommand{\vphi}{{\varphi}}

\usepackage{accents}

\newcommand{\ut}[1]{\underaccent{\tilde}{#1}}
\renewcommand{\circledcirc}[1]{\ut{#1}}

\newcommand\underrel[2]{\mathrel{\mathop{#2}\limits_{#1}}}

\newcommand{\const}{{\mathrm{const}}}

\newcommand{\cirS}{\mathop{\bigcirc\kern -.73em {\scriptstyle\mathrm{S}}}}

\theoremstyle{plain}

\theoremstyle{definition}

\newtheorem*{bemx}{Bemerkungen und Beispiele}

\begin{document}
\bibliographystyle{plain}

\begin{titlepage}
        \begin{center}
        \Huge{Lectures on Minimal Surfaces} \\[4em]
        \Large{Jens Hoppe}\\
        [1em]
        \large{Institut des Hautes \'{E}tudes Scientifiques\\
        	Bures sur Yvette\\
        	91440 France}
        \\[11em]
        \large{. }
        \\[3em]
        \large{.}\\
        \large{Abstract}\\
              {\normalsize Some elementary considerations are presented concerning Catenoids
              and their stability, separable minimal hypersurfaces, minimal surfaces obtainable by rotating shapes, determinantal varieties, minimal tori in $S^3$, the minimality in $\mathbb{R}^{n k}$ of the ordered set of k orthogonal equal-length n-vectors, and U(1)-invariant minimal 3-manifolds.}
        \large{.}
        \\[3em]
        {.}
      \end{center}

\end{titlepage}



\setcounter{page}{1}

\renewcommand*\contentsname{}
\tableofcontents

\newpage



\secct{Solitonic Catenoids} \label{sec-I}

Let us start with the following simple question: what is the surface of least area connecting $ 2 $ circles of radius  $r$ lying ( one above the other ) in parallel planes a distance $ d $ apart?\\
As the answer ( if it exists) should clearly be of the form
\begin{equation}\label{1}
\stackrel{\rightharpoonup}{x}(z, \varphi)=\left(\begin{array}{c}
f(z)\cos \varphi\\f(z)\sin\varphi\\z
\end{array}\right)\quad\varphi\in [0, 2\pi]
\end{equation}
one may determine the optimal surface by minimizing, among functions $ f $ taking the value $r$ at the boundary  (say $z = \pm \frac{d}{2}$)
\begin{equation}\label{2}
\mathbb{A}[\stackrel{\rightharpoonup}{x}]=\mathbb{A}[f]=2\pi\int_{-\frac{d}{2}}^{+\frac{d}{2}}f\sqrt{1+f'^2} ~dz,
\end{equation}
the area of an axially symmetric surface of the form (\ref{1}). Varying $ f $
($f\to f+\eps$, $\eps(\pm \frac{d}{2})=0$) gives
\[
A(f,f'):=f\sqrt{1+f'^2}\to f\sqrt{1+f'^2}+\eps \sqrt{1+f'^2}+\frac{f\eps' f'}{\sqrt{1+f'^2}}+\ldots, resp
\]
\[
\sqrt{1+f'^2}=\frac{\partial A}{\partial f}\stackrel{!}{=}\frac{d}{dz}\left(\frac{\partial A}{\partial f'}\right), ~i.e.~~1+f'^2=ff'',
\]
hence
\begin{equation}\label{3}
f(z)=a\cosh \left(\frac{z}{a}+c\right).
\end{equation}
The boundary conditions at $\pm\frac{d}{2}$ imply $c=0$, and $a \cosh \frac{d}{2a}=r$, which written as
\begin{equation}\label{4}
\frac{\cosh w}{w}=\rho ; ~~w:=\frac{d}{2a}, ~~\rho:=\frac{2r}{d}
\end{equation}
is easily seen to generically, for $\rho>\bar{\rho}$, have $ 2 $ solutions $w_1<w_2$ ( corresponding to an outer resp.inner 'Catenoid'), of closest approach to the z-axis $a_1>a_2$, while no solution exists for
$\rho<\bar{\rho}$ ( large distance, resp. small radius ) and exactly one solution $w_0$, $\frac{\cosh w_0}{w_0}=\bar{\rho}$, at critical distance/ radius.\\
 Evaluating (\ref{2}) for $(\ref{3})_{c=0}$  gives
\begin{align}
\mathbb{A}[f] &=4\pi a^2 \int_0^{\frac{d}{2a}} (\cosh r)^2 dr\nonumber\\
&=\pi a^2 \left( \frac{d}{a}+\sinh \frac{d}{a}\right)\label{5}\\
&=\frac{\pi d^2}{2}\left(\frac{1}{w}+\frac{\sinh 2 w}{2w^2}\right)=: A(w)\nonumber
\end{align}
Naturally one would like to compare the two areas
\begin{align}\label{6}
A_i:=A(w_i)=\frac{\pi d^2}{2w_i}(1+\rho \sinh w_i),
\end{align}
as well as determine whether they are really local minima, or saddle points resp. maxima. Using (\ref{6}), and then again (\ref{4}), one finds that $A_1<A_2$, as
\begin{align}
\frac{2w_1w_2}{\pi d^2}(A_2-A_1)&=(w_1-w_2+\cosh w_1\sinh w_2-\cosh w_2 \sinh w_1)\nonumber\\
&=\sinh(w_2 - w_1)-(w_2-w_1)>0,
\end{align}
while calculating the second variation of (\ref{2})
\[
(f\to f+\eps, f=a\cosh \frac{z}{a}, \eps(z)=\tilde{\eps}(v:=\frac{z}{a}))
\]
gives
\begin{equation}\label{8}
\delta^2 A = \pi\int_{-w}^{+w}\frac{\tilde{\eps}(v)}{\cosh v}\left(-\partial_v^2 - \frac{2}{\cosh^2 v}\right)\frac{\tilde{\eps}}{\cosh v}dv
\end{equation}
as
\begin{align}\label{9}
\mathbb{A}[f+\eps]&=2\pi\int_{-\frac{d}{2}}^{+\frac{d}{2}}f\sqrt{1+f'^2}\left(1+\frac{\eps}{f}+\ldots\right)
\left(1+\frac{\eps'f'}{1+f'^2}+\frac{\eps'^2}{2(1+f'^2)^2}+\ldots\right) dz\nonumber\\
&=\mathbb{A}[f]+\delta A+...\nonumber\\
&=\ldots+\frac{2\pi}{2}\int_{-w}^{+w}\left(\frac{\tilde{\eps}'^2}{\cosh^2 v}+(\tilde{\eps}^2)'\frac{\sinh v}{\cosh v}\right)dv\nonumber\\
&=\ldots+\frac{2\pi}{2}\int_{-w}^{+w}\left\{-\frac{\tilde{\eps}\tilde{\eps}''}{\cosh^2}-\tilde{\eps}^2\left(\frac{\sinh}{\cosh^3}\right)'-\frac{\tilde{\eps}^2}{\cosh^2}\right\}\nonumber\\
&=\ldots+\frac{2\pi}{2}\left(\int_{-w}^w\frac{\tilde{\eps}^2}{\cosh^2 v} dv-\int_{-w}^w\frac{\tilde{\eps}\tilde{\eps}''}{\cosh^2 v} dv-3
\int_{-w}^w\frac{\tilde{\eps}^2}{\cosh^4 v} dv\right),
\end{align}
always using that the variation vanishes at the boundary, $\tilde{\epsilon}(\pm w)=0$. ( and leaving out terms of order $ 3 $ and higher)\\
The stability properties of the Catenoids therefore depend on the spectrum of the operator $J_w$ given as
\begin{equation}\label{10}
J:=-\partial^2_v-\frac{2}{\cosh^2 v}
\end{equation}
acting on functions vanishing at the boundary of the interval $I_w:=[-w,+w]$.

Recalling that $w=\frac{d}{2a}$ is determined by solving
\[
\cosh w=\frac{2r}{d}w=:\rho w,
\]
note that
(for $\rho=\bar{\rho}$) the straight line  $\rho\cdot w$ will be tangent to $ \cosh w $, hence
$\bar{\rho}=\sinh w_0$, and therefore
\begin{equation}\label{11}
w_0\tanh w_0=1.
\end{equation}
This observation is important as it allows one to conclude that
$J^{(0)}:= J_{w_0}$ is non-negative, as
\begin{equation}\label{12}
\psi_0(v):=1-v\tanh v,
\end{equation}
which can easily be seen to be annihilated by $ J^{(0)}$, being non-negative on $I_{w_0}$ and vanishing on its boundary , must be the groundstate of $J^{(0)}$, hence $J^{(0)}\geq 0$.\\
While this in particular means that surfaces corresponding to
\begin{align}\label{13}
f_{\gamma}(z)=a_0 \cosh \frac{z}{a_0}\left(1+\gamma\left(1-\frac{z}{a_0}\tanh\frac{z}{a_0}\right)\right)
\end{align}
will have, up to ( and including ) second order in $\gamma\ll 1$, the same area than the critical Catenoid, the main virtue of knowing the lowest eigenvalue of $J^{(0)}$ to be zero is that it allows one to conclude that ( for $\rho>\bar{\rho}$) the outer Catenoid is stable while the inner one is unstable!
For that one can either invoke the fact that the lowest eigenvalue of $ \mathcal{J} $ has to increase ( decrease ) when the length of the interval is decreased ($w_0\searrow w_1$) resp increased ($w_0\nearrow w_2$)-or explicitly argue as follows:\\
As is well known from integrable systems ( see e.g. chapter $ 15 $ of \cite{Hoppe92}), while apparently less common knowledge in the context of minimal surfaces,
$J$ can factorized, and be related to the free operator $\tilde{J}:= -\partial^2:$
\begin{align}\label{14}
&J=-\partial^2 - \frac{2}{\cosh^2}=(-\partial+\tanh)(\partial+\tanh) -1=:L^{\dagger}L-1\\
&LL^{\dagger}-1=(\partial+\tanh)(-\partial+\tanh)-1=-\partial^2= \tilde{{J}}\nonumber.
\end{align}
Hence ( forgetting for the moment the boundary condition, i.e. only on the level of solutions to differential equations ) solutions of $J\psi_E=E\psi_E$ can be constructed as
\begin{align}\label{15}
\psi_E=L^{\dagger}\phi_E~, ~~-\partial^2\phi_E=E\phi_E.
\end{align}
E.g. for $E=0$ ($\phi_0^+= const~,~ \phi_0^-=(const)\cdot v$)\\
one obtains
(for $\phi_0^-=-v$)
\begin{equation}\label{16}
\psi_0^{(+)}=L^{\dagger}\phi_0^-=(-\partial +\tanh v)(-v)=1-v\tanh v ,
\end{equation}
explaining a way to derive the explicit form (\ref{12}). \footnote{Many thanks to J.Choe for pointing out that, geometrically,(\ref{12}) is the projection of the position vector on to the surface normal. } \\
To explicitly construct the instability - mode of the inner Catenoid consider
\begin{equation}\label{17}
\psi^+_k(v)=L^{\dagger}\left(\frac{-\sinh kv}{k}\right)=\cosh kv-\tanh v \cdot \frac{\sinh kv}{k};
\end{equation}
 the normalization is taken to smoothly reduce to (\ref{12}) as $k\to 0$, $\psi_k(0)=1$, and the superscript ( left out from now on ) indicating the parity of the function . While for generic $k$, $\psi_k$ will not vanish on the given boundary
($\pm v=w_2=w_0+\epsilon, \epsilon>0$) but for some ('minimal') $k>0$ one will have
\begin{equation}\label{18}
\psi_k(\pm w_2)=\cosh kw_2- \tanh w_2 \frac{\sinh kw_2}{k}=0,
\end{equation}
$\psi_k(v)>0$ on ($-w_2, +w_2$).\\
To conclude this, note that
$\psi_k$, for fixed $k$ (and restricting to $v\geq 0$), is monotonically decreasing ( at least for $ k $ not too large; for simplicity let us consider small $k>0$, $\epsilon=( w_2 - w_0)\ll 1$, close to the critical case ):
\begin{equation}\label{19}
\psi'_k=k\sinh kv-\frac{1}{\cosh^2 v} \frac{\sinh kv}{k}-\tanh v \cosh kv<0.
\end{equation}
For $k=0$, $\psi_0$ vanishes at $w_0$ and then, in the (small) interval $[w_0, w_2]$ becomes negative. To conclude that ( for fixed $ v $) $\psi_k$ is monotonically increasing with $k$ ( near zero ),
\begin{equation}\label{19}
\frac{\partial\psi_k}{\partial k}(v)=v\sinh kv- v\tanh v \frac{\cosh kv}{k}+
\tanh v \frac{\sinh kv}{k^2}>0,
\end{equation}
one calculates the Taylor-expansion ( $v$ fixed )
\begin{equation}\label{20}
\frac{\partial\psi_k}{\partial k}(v)=\frac{1}{k}\cdot 0+\frac{kv^2}{3}\left(
\underbrace {3-v\tanh v}_{>0}\right)+O(k^3).
\end{equation}
Taylor expanding $\psi_k(\pm(w_0+\eps))=0$, cp. (\ref{18}),
\begin{equation}\label{21}
1+\frac{k^2(w_0+\eps)^2}{2}-\left(\tanh w_0+\frac{\eps}{\cosh^2 w_0}\right)\left(w_0(1+\frac{\eps}{w_0})+\frac{1}{6}k^2(w_0+\eps)^3\right)+O(k^4)\stackrel{!}{=}0
\end{equation}
yields
\begin{equation}\label{22}
k^2=\frac{3\eps}{w_0},
\end{equation}
while Taylor expanding (\ref{17}) gives
\begin{equation}\label{23}
\psi_k(v)=(1-v\tanh v)+\frac{k^2v^2}{6}\left(3 - v\tanh v\right)+O(k^4).
\end{equation}
As an interesting exercise one may compare/double check/ these results with $ 2 $ ( not-completely- ordinary. non-standard ) perturbation theory calculations:\\
Making in $J\psi_E=E\psi_E$ the Ansatz  $\psi_E=\psi_0+\delta_E$ one gets, via (\ref{14}),
\begin{equation}\label{24}
(L^{\dagger}L-1)\delta_E=E(\psi_0+\delta_E)\approx E\psi_0
\end{equation}
resp. ( acting with $ L $ on both sides )
\begin{equation}\label{25}
-\partial^2(L\delta_\eps)\approx EL(1-v\tanh v)=-Ev,
\end{equation}
i.e. ( integrating $2$ times and approximating $\delta_E$ by the exact solution of the approximated equation (\ref{24}))
\begin{equation}\label{26}
L\delta_\eps=E\left(\frac{v^3}{6}+\alpha+\beta v\right) \stackrel{!}{=}
\frac{1}{\cosh}(\delta_\eps\cosh)',
\end{equation}
from which one deduces
\begin{equation}\label{27}
\delta_\eps=\frac{e}{\cosh v}-\frac{E v^2}{6}(3-v\tanh v)-(\beta+1)(1-v\tanh v)E
\end{equation}
( setting $\alpha=0$ for parity-reasons ).\\
Calculating/ checking
\begin{align}\label{28}
(L^{\dagger}L-1)\delta_\eps&=-\frac{e}{\cosh v}+\left(\partial^2+\frac{2}{\cosh^2}\right)\frac{Ev^2}{6}(3-v\tanh v)\nonumber\\
&=\frac{-e}{\cosh v}+E(1-v\tanh v)\stackrel{!}{=}E_\eps \psi_0
\end{align}
one finds $e=0$, hence
\begin{equation}\label{29}
\psi_\eps=(1-E(\beta+1))(1-v\tanh v)-\frac{Ev^2}{6}(3-v\tanh v);
\end{equation}
while $\beta\neq -1$ just changes the normalization ( hence now
$\beta=-1$ for simplicity ) $\psi_\eps(\underbrace{w_0+\eps}_{w_\eps})\stackrel{!}{=}0$ gives
\begin{equation}\label{30}
E(\eps)=-\frac{(w_\eps\tanh w_\eps-1)\cdot 6}{w^2_\eps(3-w_\eps\tanh w_\eps)}\approx-\frac{3\eps}{w_0},
\end{equation}
in accordance with (\ref{23}), resp. (\ref{22}).\\
The non-standard part of this derivation is the use of the $LL^{\dagger}$ structure that allows one to calculate the perturbed wavefunction without having to use all eigenstates of the unperturbed problem.\\
The subtlety, on the other hand, when wanting to use the standard first-order formula for the perturbed eigenvalue, is that $ J_{w_\eps}$ is a perturbation of $ J^{(0)}= J_{w_0}$ only via the boundary condition \footnote{thanks to R.Hempel for pointing out to me the corresponding general treatment of boundary-perturbations given in \cite{Kato80}.}, to effectively have $ J_{w_\eps}$ in the standard form
$ J^{(0)}+\eps J'$, note that $E_\eps$ is the minimum of
\begin{equation}\label{31}
\frac{\int_0^{w_\eps}(\psi'^2-\frac{2\psi^2}{\cosh^2 v})dv}{\int_0^{w_\eps}\psi^2 dv}
\end{equation}
subject to $\psi(w_\eps)=0$.\\

Introducing
\begin{equation}\label{32}
\hat{v}:=\frac{v}{1+\frac{\eps}{w_0}}\in [0, w_0],~~ \hat{\psi}(\hat{v}(v))=\psi(v)
\end{equation}
(\ref{31}) becomes equal to

\begin{equation}\label{33}
\frac{1}{(1+\frac{\eps}{w_0})^2}
\frac{\int_0^{w_0}(\hat\psi'^2-\frac{2\hat\psi^2(1+\frac{\eps}{w_0})^2}{\cosh^2 (\hat v (1+\frac{\eps}{w_0}))})d\hat v}
{\int_0^{w_0}\hat\psi^2 d\hat v}.
\end{equation}
With $\cosh(\hat{v}+\frac{\eps \hat{v}}{w_0})\approx\cosh \hat{v}\left(
1+\frac{\eps\hat{v}}{w_0}\tanh \hat{v}\right)$ one gets
\begin{equation}\label{34}
J^{\prime}=\frac{-2}{\cosh^2}\left(\frac{2}{w_0}-\frac{2}{w_0}\hat{v}\tanh\hat{v}\right)=
-\frac{4}{w_0\cosh^2 \hat{v}}(1-\hat{v}\tanh\hat{v}),
\end{equation}
so that the standard first-order formula for $E_\eps$ gives
\begin{equation}\label{35}
E_\eps=\eps\frac{\int_0^{w_0}\psi_0^2 J'}{\int_0^{w_0}\psi_0^2}
\underrel{\underrel{\psi_0=(1-v\tanh v)}{\downarrow}}{=}-4\frac{\psi\eps}{w_0}\frac{(J_0-3J_1+3J_2-J_3)}{K_0-2K_1+K_2}
\end{equation}
where
\begin{equation}\label{(36)}
J_n:=\int_0^{w_0}\frac{v^n(\tanh v)^n}{\cosh^2 v}dv, ~~
K_n=\int_0^{w_0}v^n(\tanh v)^n dv.
\end{equation}
With the help of
\begin{equation}\label{37}
J_n=\frac{\tanh w_0}{n+1}-\frac{n}{n+1}K_{n+1}+\frac{n}{n+1}J_{n-1}
\end{equation}
one finds
\begin{equation}\label{38}
J_0=\frac{1}{w_0}, J_1=\frac{1}{w_0}-\frac{w_0}{2}, ~~J_2=\frac{1}{w_0}-\frac{w_0}{3} - \frac{2}{3}K_1, ~~J_3=\frac{1}{w_0}-\frac{w_0}{4}-\frac{1}{2}K_1-\frac{3}{4}K_2
\end{equation}
so that
\begin{equation}\label{39}
J_0-3J_1+3J_2-J_3=\frac{1}{4}w_0^3,
\end{equation}
and with $K_0=w_0$, $K_2=\frac{1}{3}w_0^3-w_0+2K_1$ ( note that the non-elementary $K_1$ cancels both in the numerator as well as the denominator ) (\ref{35}) becomes
\begin{equation}\label{40}
E_\eps=\frac{-4\eps}{w_0}\cdot \frac{3}{4},
\end{equation}
in agreement with (\ref{22}) and (\ref{30}).
\bigskip\\

Consider now fluctuations around the Catenoid $(\times\mathbb{R})$ as a minimal (hyper-)surface in $\mathbb{R}^{3,1}$, i.e. a time-independent $(\dot{z}=0)$ stationary point of the volume-functional (see e.g. \cite{Hoppe94})

\begin{equation}\label{42}
2\pi\int rdrdt \sqrt{1-\dot{z}+z'^{2}}=\int rdzdt \left| r' \right| \sqrt{1-\frac{\dot{r}^{2}}{r'^{2}}+\frac{1}{r'^{2}}}=\int rdtdz \sqrt{1+r'^{2}-\dot{r}^{2}},
\end{equation}

for (in general time-dependent) axially symmetric hypersurfaces in Minkowski-space,

\begin{equation}\label{43}
x^{\mu}=x^{\mu}(t,r,\varphi)=\begin{pmatrix}
t \\
r\cos\varphi \\
r\sin \varphi \\
z(r,t)
\end{pmatrix} = \begin{pmatrix}
t \\
r(z,t)\cos\varphi \\
r(z,t)\sin \varphi \\
z
\end{pmatrix} = x^{\mu}(t,z,\varphi),
\end{equation}

with $(r,t)\rightarrow(z(r,t),t)$ implying

\begin{equation}\label{44}
	\begin{pmatrix}
	z' & \dot{z} \\
	0 & 1
	\end{pmatrix} = \left( \frac{\partial(z,t)}{\partial(r,t)}\right) = \left( \frac{\partial (r,t)}{\partial (z,t)}\right)^{-1} = \begin{pmatrix}
	r' & \dot{r} \\
	0 & 1
	\end{pmatrix}^{-1} = \frac{1}{r'}\begin{pmatrix}
	1 & -\dot{r} \\
	0 & r'
	\end{pmatrix} = \begin{pmatrix}
	\frac{1}{r'} & \frac{-\dot{r}}{r'} \\
	0 & 1
	\end{pmatrix},
\end{equation}

and the variation of the rhs of (\ref{42}) giving

\begin{equation}\label{45}
	r''(1-\dot{r}^{2})-\ddot{r}(1+r'^{2})+2\dot{r}r'\dot{r}=\frac{1}{r}(1+r'^{2}-\dot{r}^{2}),
\end{equation}

the stationary Cateniod $ r(z,t)=r(z)=\cosh z $ indeed being a solution, as $ r''=\frac{1}{r}(1+r'^{2}) $. Linearization of (\ref{45}) around $ \cosh z=r(z,t)-\epsilon(z,t) $ gives

\begin{equation}\label{46}
	\ddot{\epsilon}+D\epsilon=0
\end{equation}

with

\begin{equation}\label{47}
	D=-\frac{1}{\cosh^{2}z}\left( \partial_{z}^{2}-2\tanh z\, \partial_{z}+1\right).
\end{equation}

While $ D $ has two zero-modes,

\begin{align}\label{48}
	\epsilon_{+}(z)&=\cosh z-z\sinh z\nonumber\\
	\epsilon_{-}(z)&=\sinh z
\end{align}	

(corresponding to the 2 parameters in the time-independent solutions of (\ref{45}), $ \frac{1}{a}\cosh(az+b) $), it is interesting to note that $D$ has (exactly) one positive parity eigenfunction with energy $ \lessapprox -\frac{8}{15}\left( =\frac{\left\langle \psi,D\psi\right\rangle }{\left\langle \psi,\psi\right\rangle }, \text{ where }\psi=\frac{1}{\cosh z}\right) $, while being positive on negative parity eigenfunctions, as

\begin{align}\label{49}
	D&=\frac{1}{\cosh z}\left( \partial_{z}+\frac{1}{\sinh z\cosh z} \right) \left( -\partial_{z}+\frac{1}{\sinh z\cosh z}\right) \frac{1}{\cosh z}\nonumber\\
	&=\left( \partial\frac{1}{\cosh}+\frac{1}{s}\right)\left( -\frac{1}{\cosh}\partial+\frac{1}{s}\right).
\end{align}

This having been noticed at least 10 years ago [J. Hoppe, unpublished note to G.Huisken], the question of non-linear stability was taken up,and answered, more recently \footnote{thanks to J.Szeftel for discussions and bringing \cite{Donninger2016} to my attention}. Let us mention a few facts/things related to the endeavour of trying to find a closed expression for the unstable mode of (\ref{47}), resp. (expressed in the coordinate $ y=\sinh z $, hence $\partial_{z}=\cosh z\, \partial_{y}, \partial_{z}^{2}=y\, \partial_{y}+(1+y^{2})\, \partial_{y}^{2}$, and compensating $dz=\frac{1}{\sqrt{1+y^{2}}}dy$ by conjugation with $(1+y^{2})^{\pm\frac{1}{4}}$),

\begin{align}\label{50}
	\tilde{D}&=-(1+y^{2})^{-\frac{5}{4}}\left( (1+y^{2})\, \partial_{y}^{2}-y\, \partial_{y}+1\right)(1+y^{2})^{\frac{1}{4}} \nonumber\\
	&=-\partial _{y}^{2}-\frac{1}{4(1+y^{2})}-\frac{5}{4}\frac{1}{(1+y^{2})^{2}}=-\partial_{y}^{2}+\tilde{V}(y),
\end{align}

which also follows as

\begin{align}\label{51}
	(1+y^{2})^{-\frac{1}{4}}\left[ \left( -\partial_{y}+\frac{y}{1+y^{2}}\right) \partial_{y}-\frac{1}{1+y^{2}}\right] (1+y^{2})^{\frac{1}{4}} \nonumber\\
	=\left( -\partial_{y}+\frac{1}{2}\frac{y}{1+y^{2}}\right) \left( \partial_{y}+\frac{1}{2}\frac{y}{1+y^{2}}\right) -\frac{1}{1+y^{2}},
\end{align}

noting (cp. (\ref{8}), \cite{Graf2019}, \cite{Hoppe92})

\begin{align}\label{52}
	D&=\frac{1}{\cosh z} \left( -\partial_{z}^{2}-\frac{2}{\cosh^{2}z} \right) \frac{1}{\cosh z}\nonumber\\
	&=\frac{1}{\cosh z}\left[ \left( -\partial_{z}+\tanh z\right) \left( \partial_{z}+\tanh z \right) -1 \right] \frac{1}{\cosh}.
\end{align}

One reformulation of trying to solve

\begin{equation}\label{53}
	\tilde{D}\tilde{\phi}=-\kappa^{2}\tilde{\phi},
\end{equation}

$ \int\tilde{\phi}^{2}dy<\infty,$ $ \tilde{\phi}(y)\neq0,$ arises from the factorization

\begin{equation}\label{54}
	\left( -\partial_{y}^{2}+\tilde{V}+\kappa^{2} \right)\stackrel{!}{=}\left( -\partial_{y}+U \right)\left( \partial_{y}+U\right) \geq 0
\end{equation}

giving the Riccati-equation

\begin{equation}\label{55}
U'=U^{2}-\tilde{V}-\kappa^{2}=U^{2}+\frac{B}{1+y^{2}}+\frac{D}{(1+y^{2})^{2}}-\kappa^{2}
\end{equation}

(with $ B=\frac{1}{4}, D=\frac{5}{4}$ in the case of interest), resp.(using that one can choose the eigenfunctions of (\ref{50}) to be either odd , or - as in the case of the ground state - even, in both cases having $ U=-\frac{\tilde{\phi}'(y)}{\tilde\phi(y)} $ to be odd, and , with $\tilde{\phi}(y)=\chi(x:=y^{2}) $,

\begin{equation}\label{56}
	W(x):=-2yU(y)=+\frac{\chi'(x)}{\chi(x)}
\end{equation}

having to satisfy

\begin{equation}\label{57}
	4x(W'+W^{2})+2W+\frac{B}{1+x}+\frac{D}{(1+x)^{2}}-\kappa^{2}=0.
\end{equation}

For $ \kappa=0 $ (and $ B=\frac{1}{4}, D=\frac{5}{4}$) a particular solution is

\begin{equation}\label{58}
W_{0}(x)=\frac{1}{4x}\left( \frac{2+x}{1+x}\right) =\frac{1}{4x}\left( 1+\frac{1}{1+x}\right),
\end{equation}

corresponding to the non-normalizable zero-mode $ \varepsilon_{-}(z)=\sinh z $, resp. $\chi_{0}\sim\frac{\sqrt{x}}{(x+1)^{\frac{1}{4}}}$, resp. $\frac{y}{(1+y^{2})^{\frac{1}{4}}}=\tilde{\phi}_{0}$, and the Ansatz

\begin{equation}\label{59}
	W=W_{0}+\frac{1}{Y}
\end{equation}

then gives $ Y'-\frac{1}{2x}Y-2W_{0}Y=1 $, resp.

\begin{align}\label{60}
Y&=\tilde{C}\frac{x^{\frac{3}{2}}}{\sqrt{x+1}}+\frac{x^{\frac{3}{2}}}{\sqrt{x+1}}\int\frac{\sqrt{x+1}}{x^{\frac{3}{2}}}\\
&=\tilde{C}\frac{x^{\frac{3}{2}}}{\sqrt{x+1}}+\frac{x^{\frac{3}{2}}}{\sqrt{x+1}}\left( 2\ln\left( \sqrt{x}+\sqrt{1+x}\right)\nonumber -2\sqrt{\frac{1+x}{x}}\right),
\end{align}

i.e.

\begin{equation}\label{61}
	W(x)=\frac{1}{4x}\left( \frac{2+x}{1+x}\right) +\dfrac{1}{2\frac{x^{\frac{3}{2}}}{\sqrt{x+1}}\left( \frac{\tilde{C}}{2}-\sqrt{\frac{1+x}{x}}+\ln \left(\sqrt{x}+\sqrt{x+1} \right)\right)},
\end{equation}

which indeed agrees with the expression one gets from

\begin{align}\label{62}
	-U(y)&=\frac{\tilde{\phi}'_{0}}{\tilde{\phi}_{0}}=\dfrac{\frac{C}{2}\left( \frac{2+y^{2}}{(1+y^{2})^{\frac{5}{4}}}\right)-\frac{A}{2}\left[\frac{y}{(1+y^{2})^{\frac{3}{4}}}+\frac{2+y^{2}}{(1+y^{2})^{\frac{5}{4}}}\ln \left(y+\sqrt{1+y^{2}}\right)\right]}{C\frac{y}{(1+y^{2})^{\frac{1}{4}}}+A\left[ \left( 1+y^{2}\right)^{\frac{1}{4}}-\frac{y}{\left(1+y^{2}\right)^{\frac{1}{4}}}\ln \left(y+\sqrt{1+y^{2}}\right)\right]} \nonumber\\
	&\stackrel{=}{{\tiny A\neq 0}}\frac{1}{2}\dfrac{\frac{C}{A}\left( \frac{2+y^{2}}{1+y^{2}}\right)-\frac{y}{\sqrt{1+y^{2}}}-\frac{2+y^{2}}{1+y^{2}}\ln \left(y+\sqrt{1+y^{2}}\right)}{\frac{C}{A}y+\sqrt{1+y^{2}}-y\ln \left(y+\sqrt{1+y^{2}}\right)}
\end{align}

$\left(\text{for } A=0, \ -U=\frac{1}{2}\frac{2+y^{2}}{y(1+y^{2})},\text{ hence }W=\frac{U}{-2y}=\frac{1}{4x}\frac{2+x}{1+x}\right)$.
To solve (\ref{57}) for $ \kappa\neq 0 $, however, seems to be just as difficult as directly trying to find the groundstate of $ D $, which, using that one knows explicitly (see e.g. \cite{Hoppe92})  the exact eigenfunctions of

\begin{equation}\label{63}
	H=-\partial_{z}^{2}-\frac{2}{\cosh^{2}z},
\end{equation}

\begin{align*}
\psi^{(0)}=\frac{1}{\sqrt{2}}\frac{1}{\cosh z}, \ \psi_{k}(z)=-\left(ik+\tanh z\right) e^{-ikz},
\end{align*}

satisfying

\begin{align}\label{64}
H\psi^{(0)}&=-\psi^{(0)}, \quad \int\psi^{(0)2}=1,\quad \int\psi^{(0)}\psi_{k}=0\nonumber\\
H\psi_{k}&=k^{2}\psi_{k}, \quad \int \psi_{k}\psi_{k'}=(k^{2}+1)\delta(k-k')
\end{align}


\begin{align}\label{65}
	\psi^{0}(z)\psi^{0}(z')+\int_{-\infty}^{+\infty}\frac{dk}{k^{2}+1}\psi_{k}(z)\psi^{*}_{k}(z')=\delta(z-z'),
\end{align}

one could formulate as trying to find constants $ C_{-1} $ and $ C(k) $ satisfying

\begin{align}\label{66}
&\dfrac{-C_{-1}}{\sqrt{2}\cosh z}-\int_{-\infty}^{+\infty}C(k)\dfrac{ik+\tanh z}{\sqrt{k^{2}+1}}e^{-ikz}k^{2}dk\nonumber\\
&=-\kappa^{2}\cosh^{2}z\left(\frac{C_{-1}}{\sqrt{2}\cosh z}-\int_{-\infty}^{+\infty}C(k)\dfrac{ik+\tanh z}{\sqrt{k^{2}+1}}e^{-ikz}dk\right),
\end{align}

with the expression in brackets (on the rhs) being , when multiplied by $ \cosh z$, square-integrable.

\newpage
\secct{Separable Minimal Hypersurfaces and Rotating Shapes} \label{sec-III}

For surfaces representable as graphs over ( parts of ) $\mathbb{R}^2$ the area is expressed as
\begin{equation}\label{a1}
A[z]=\int\int \sqrt{1+z_x^2+z_y^2} ~dx dy,
\end{equation}
 whose stationary points correspond to solutions $z(x,y)$ of
\begin{equation}\label{a2}
z_{xx}(1+z_y^2)+z_{yy}(1+z_x^2)=2z_xz_yz_{xy}.
\end{equation}
Inserting the Ansatz $z(x,y)=\zeta\left(f(x)+g(y)\right)$, and denoting the inverse of $ -\zeta $ by $ h $ one finds an equation for separable surfaces,
\begin{equation}\label{a3}
\Sigma:=\{\stackrel{\rightharpoonup}{x}\in \mathbb{R}^3 | f(x)+g(y)+h(z)=0\},
\end{equation}
to be 'minimal' :
\begin{equation}\label{a4}
f''(x)\left(g'^2(y)+h'^2(z)\right)+g''(y)\left(h'^2(z)+f'^2(x)\right)+h''(z)\left(f'^2(x)+g'^2(y)\right)=0,
\end{equation}
to hold for all $(x,y,z)\in \Sigma$ (i.e. on $\Sigma$).\\
While it is easy to verify the Catenoid
\begin{equation}\label{a5}
x^2+y^2-\left(\cosh z\right)^2=0,
\end{equation}
as a solution ( on $\Sigma_{5}$ ) of (\ref{a4}),

other elementary minimal surfaces ( and in fact, after (\ref{a5}) historically the first known ones ) are the helicoid,
\begin{equation}\label{a6}
y\cos z = x\sin z,
\end{equation}
Scherk's first,
\begin{equation}\label{a7}
e^z \cos x = \cos y,
\end{equation}
and second,
\begin{equation}\label{a8}
\sin z = \sinh x \sinh y,
\end{equation}
surface.\\
Exercise:\\
Show that (\ref{a6})/(\ref{a7})/(\ref{a8}) are of the form (\ref{a3}), $\sum_{i=1}^3 f_i(x_i)=0$, with each of the functions satisfying
\begin{equation}\label{a9}
f_i'^2=a_i+b_i e^{\kappa f_i}+c_i e^{-\kappa f_i},
\end{equation}
and derive the general conditions on the coefficients appearing in (\ref{a9}) to guarantee (\ref{a4}), i.e. $\sum_{i\neq j} f_i'' f_j'^2=0$, on (\ref{a3}).

Note that varying the 'area' ( volume )
\begin{equation}\label{a10}
A[u]:=\int \delta (u(\stackrel{\rightharpoonup}{x}))\vert \stackrel{\rightharpoonup}{\nabla} u\vert ~ d^{N}x
\end{equation}
of a hypersurface described as a level set,
\begin{equation}\label{a11}
\Sigma:=\{ \stackrel{\rightharpoonup}{x}\in \mathbb{R}^N | u( \stackrel{\rightharpoonup}{x})=0\},
\end{equation}
yields the equation
\begin{equation}\label{a12}
(\nabla u)^2 \Delta u - \sum_{i,j=1}^N u_i u_j u_{ij}=0,
\end{equation}
to hold on (\ref{a11}). The separation Ansatz
\[u( \stackrel{\rightharpoonup}{x})=\sum_{k} f_k (x_k),\]
then yields that
\begin{equation}\label{a13}
\sum_{i\neq j}f_i'' f_j'^2=0
\end{equation}
should hold on $\sum f_i=0$.\\
Existence and form of the solutions heavily depend on the dimension. While for $N=3$, (\ref{a13}) was completely solved already $ 130 $ years ago \cite{Weingarten87}, and the earliest attempt for $N=4$ seems to be in the Lorentzian context \cite{Hoppe95}, a complete classification for $N\geq 4$ has been attacked (and more or less completed) recently together with J.Choe and V.Tkatjev.\\
For $N=3$, if none of the $ 3 $ functions is linear, solutions are of the form
\begin{equation}\label{a14}
f_i'^2 =a_i+b_i e^{\sqrt{\mu}f_i}+c_i e^{-\sqrt{\mu}f_i} ,
\end{equation}
$\mu\neq 0$, with the $ 9 $ constants linked by non-linear equations allowing for solutions in terms of $ 5 $ free constants; resp.(`$\mu=0$')
\begin{equation}\label{a15}
f_i'^2=a_i+b_i f_i+c_i f_i^2,
\end{equation}
with the coefficients satisfying another set of non-linear equations.\\
Apart form the fully linear case
\begin{equation}\label{a16}
f_i=\alpha_i+\beta_i x_i
\end{equation}
the, up to permutation, only other case ( cp. [$ N $] ) is
\begin{align}\label{a17}
&f_1'^2=a_1 + b_1 e^{\lambda f_1}\nonumber\\
&f_2'^2=a_2 + c_2 e^{-\lambda f_2} ~~~\lambda>0\\
&f_3(x_3)=\alpha x_3+\beta,\nonumber
\end{align}
corresponding e.g  (if $ a_1, a_2 < 0 $, $b_1, c_2 > 0 $ )
to
\begin{align}\label{a18}
\begin{array}{l}
f_1(x_1)=\frac{-2}{\lambda}ln\left(\sqrt{\frac{-b_1}{a_1}}\cos\left(\frac{x_1\lambda}{2}\sqrt{-a_1}\right)\right)\\\\
f_2(x_2)=\frac{2}{\lambda}ln\left(\sqrt{\frac{-c_2}{a_2}}\cos\left(\frac{x_2\lambda}{2}\sqrt{-a_2}\right)\right)
\end{array},
\end{align}
including (\ref{a7}).

Other choices of sign combinations give $u_{\pm}=\pm ln \cosh x$ and $\pm ln \sinh x:=v_{\pm}$ as constituents ( satisfying $u_{\pm}'^2=1-e^{\mp2u_{\pm}}$,
$u''_{\pm}=\mp(u'^2_{\pm}-1)$, resp. $v_{\pm}''=\mp(v_{\pm}'^2-1)$,
$v_{\pm}'^2=1+e^{\mp 2v_{\pm}}$, instead of
$w_{\pm}=\pm ln\cos x$, which satisfies
$w_{\pm}'^2=-1+e^{\mp w_{\pm}2}$,
$w_{\pm}''=\mp(w_{\pm}'^2+1)$).
While the solutions of
(\ref{a14}) are in general elliptic functions, special cases will yield trigonometric / hyperbolic expressions, such as
($\kappa=4$, $b_1=b_2=c_3=0$, $a_1=a_2=-a_3=b_3=c_1=c_2=1$)
\begin{equation}\label{a19}
\sin z=\sinh x \cdot \sinh y.
\end{equation}
What about $N\geq 4$?

Nonlinear solutions of the form

\begin{equation}\label{a20}
z(x_1, x_2, \ldots, x_n)=\sum_{i=1}^{n=N-1}z_i(x_i),
\end{equation}
apparently do not exist:\\
while the resulting equation
\begin{equation}\label{a21}
\sum_{i=1}^n z_i''\left(1+\sum_{j\neq i}z_j'^2\right)=0
\end{equation}
is ' trivially ' solved for $n=2$, letting
\begin{equation}\label{a22}
z_1''=c(1+z_1'^2), ~~z_2''=-c(1+z_2'^2),
\end{equation}
\begin{equation}\label{a23}
z_i''(z_j'^2+c_{ij})+z_j''\underset{i\neq j}{(z_i'^2+c_{ji})}=s_{ij} ,
\end{equation}
does not (seem to) have any non-trivial solutions once  $n>2$ ( i.e. only the linear solution
$z_i(x_i)=a_i x_i+b_i$ ).\\

For $N=4$, the Ansatz
\begin{equation}\label{a24}
f'^2_L=\alpha_L e^{+\kappa f_L}+\beta_L e^{-\kappa f_L},
\end{equation}
yields solutions of (\ref{a13}), provided
\begin{equation}\label{a25}
\alpha_L \alpha_K=\beta_{L'}\beta_{K'},
\end{equation}
($L, K, L',K'\in \{1, 2, 3, 4\} ~all~ different!$)
e.g.
\begin{align}\label{a26}
\begin{array}{l}
\alpha_1=\alpha_2=1=\beta_3=\beta_4\\\\
\alpha_3=\alpha_4=-1=\beta_1=\beta_2,
\end{array}
\end{align}
yielding elliptic functions, resp.
\begin{equation}\label{a27}
\Sigma_3:=\{\stackrel{\rightharpoonup}{x}\in \mathbb{R}^4|~\wp(x)\wp(y)=\wp(z)\wp(v)\}
\end{equation}
where $\wp$ is an elliptic Weierstrass-function, satisfying
\begin{equation}\label{a28}
\wp'^2=4\wp(\wp^2-1),
\end{equation}
with $\wp$ real, taking its minimum, $ 1 $, at half-period ( while diverging at $ 0 $, $2w$, ...) .\\
Note that $\alpha_3=\alpha_4=0=\beta_1=\beta_2$ (and $\kappa=2$, $\alpha_1=\alpha_2=1=\beta_3=\beta_4$), i.e.
$f'^2_{i=1,2}=e^{2f_i}$, $f'^2_{i=3, 4}=-e^{-2f_i}$ gives the known solution
\begin{equation}\label{a29}
x_1\cdot x_2=x_3\cdot x_4.
\end{equation}
 A slightly more elegant route to solving (\ref{a13}) is to note that
\begin{equation}\label{a30}
f'^2_L(x_L)=\epsilon_LQ_L(f_L(x_L))=F_L(v_L)=f_L(x_L),~~~L=1, 2, 3, 4 ~~, ~ \epsilon_L\in \mathbb{R},
\end{equation}
implying $2f'_Lf_L''=\epsilon_L f'_L Q'_L$ will (when $f'_L\neq 0$) solve (\ref{a13}), resp.
\begin{equation}\label{a31}
\sum_{L\neq K}\epsilon_K \epsilon_L Q_L(f'_L)Q'_K(f_K)=0,
\end{equation}
provided
\begin{equation}\label{a32}
Q_L(f_L)Q'_K(f_K)+Q_K(f_K)Q'_L(f_L)=R_{LK}(f_L+f_K)
\end{equation}
with $R_{LK}=R_{L'K'}$ (cp. (\ref{a25})) any (!) odd function of its argument, resp.
$R_{LK}=(-)^{\sigma_{LK}}R_{L'K'}$ having parity  $(-)^{\sigma_{LK}+1}$.

Examples are
\begin{align}\label{a33}
&Q(f)= af(\to x_1^2+x_2^2=x_3^2+x_4^2 ~~ eg.)\nonumber\\
&Q(f)=\cosh f~~or~~\sinh f\\
&Q(f)=\cos f~~or~~\sin f,
\end{align}
the trigonometric $Q'$s giving elliptic solutions, like (\ref{a27}).
\bigskip

\textbf{Separable Minimal Hypersurfaces in $\mathbb{R}^{N\geq 4}$}:
\[ \sum =\left\{ \stackrel{\rightharpoonup}{x}\in \mathbb{R}^N |  u(\stackrel{\rightharpoonup}{x}):=\sum_{i=1}^N f_i (x_i)=0\right\} \]
\[
(\nabla u)^2 \Delta u - u_i u_j u_{ij}=0 \Rightarrow \sum_{i\neq j} f'^2_i f''_j=0 ~~on ~~u=0.
\]
Defining  $J_i(v_i=f_i (x_i)):=f'^2_i(x_i)$ (i.e. necessarily nonnegative) the basic equation to solve is
\begin{equation}J:=\sum_{i\neq k} J_i J'_k\approx 0 \;  \text{(i.e. $=0$ on $\sum v_j=0$)}. \tag{J}  \end{equation}

Differentiating the basic equation $J\approx 0$, using that  $F\approx 0$ implies  $\partial_{v_i} F -\partial_{v_k}F\approx 0$, and  $F\approx 0$ together with $F$ not depending on one of the $v_j$ implying  $F\equiv 0$, gives ( applying $\partial_k - \partial_i$,and denoting $\Sigma_k:=\sum_{i\neq k} J_i$ , $\Sigma_k':=\sum_{i\neq k} J_i'$ )
\begin{equation}\label{a34}
J_{ki}:=J''_k \Sigma_k+J'_k \Sigma'_k-J''_i\Sigma_i-J'_i\Sigma'_i\approx 0
\end{equation}
as well as (applying $\partial_l - \partial_n$, .... all different)
\begin{equation}\label{a35}
J_{ki, ln}:=(J''_k - J''_i)(J'_l - J'_n)+(J''_l - J''_n)(J'_k - J'_i)\approx 0,
\end{equation}
and( applying again $\partial_k - \partial_i $),
\begin{equation}\label{a36}
J_{(ki)^2, ln}=(J'''_k +J'''_i)(J'_l - J'_n)+(J''_k+J_i'')(J''_l-J''_n)\approx 0,
\end{equation}
hence, multiplying by ($J'_k - J'_i$), and using (\ref{a35}),
\begin{equation}
(J'_l-J'_n)\left\{ J'''_k+J'''_i)(J'_k-J'_i) - (J''_k+J''_i)(J''_k-J''_i)\right\}\approx 0
\end{equation}
For every $ l $, $ n $ with $(J'_l - J'_n)\not\approx 0 $ one therefore has
\begin{equation}\label{a36}
(J'''_k+J'''_i)(J'_k-J'_i)-(J''_k+J''_i)(J''_k-J''_i)=0
\end{equation}
 ( for all distinct ( $ki$ ) different from ($ ln $)).\\
 Differentiating with respect to $ v_{k} $ gives
\begin{equation}\label{a37}
J''''_k J'_k - J'''_k J''_k=J''''_k J'_i -J''_k J'''_i
\end{equation}
and then w.r.t. $v_i$, finally,
\begin{equation}\label{a38}
J''''_k J''_i-J''_k J''''_i=0
\end{equation}
i.e. ( here derived for all $N>3$ and all i different from l,n for which $J_l'\neq J_n'$; for N=3 see \cite{Weingarten87}), if $J_k$ is nonlinear,
\begin{equation}\label{a39}
J''''_i=cJ''_i   ;
\end{equation}
implying ( and when inserted into (40), $e_i=e$ )
\begin{equation}\label{a40}
\begin{array}{l}
J'''_i=cJ'_i+e_i\\\\
J''_i=cJ_i+e_i v_i+d_i
\end{array}
\end{equation}
i.e. ( if there are at least two nonlinear $J_j$'s )
\begin{equation}\label{a41}
J_i=\alpha_i e^{\sqrt{c}v_i}+\beta_i e^{-\sqrt{c}v_i}-\frac{d_i}{c}-\frac{e_i}{c}v_i
\end{equation}
for  $c\neq 0$,
and
 \begin{equation}\label{a42}
 J_i=\dfrac{e_i}{6} v_i^3+\frac{d_i}{2}v_i^2+b_iv_i+a_i
 \end{equation}
for $c=0$. Inserting (\ref{a40}) into (\ref{a36}) yields $e_i=e$ and (with some separation constant d)
 \begin{align}\label{a43}
 (J''_i)^2&=c(J'_i)^2+2eJ'_i+d\nonumber\\
 &\stackrel{(\ref{a40})!}{=}(cJ_i+ev_i+d_i)^2
 \end{align}
As the form (44/45) includes linear functions one can ( if at least 2 of the $J_j$'s are non-linear ) simply insert it , for all! i and k, into (J) - finding that non-linearities actually are impossible if $N>4$ , and for N=4 ( the linear parts having to vanish because of the single positive and negative exponentials necessitating opposite signs for the linear parts ) the only possibility being

\begin{equation}\label{a44}
J_i(v_i)=\alpha_i e^{\sqrt{c}v_i}+\beta_i e^{-\sqrt{c}v_i}
\end{equation}
with ( $ iki'k' $ all different )
\begin{equation}\label{a45}
  \alpha_i \alpha_k =\beta_{i'}\beta_{k'}
\end{equation}
  ( which in particular implies that the products of the coefficients of the $v_i$-exponentials are independent of i, a condition that one also finds as a consequence of (39) ).

Let us now discuss the special case that all $J_{\alpha=i>1}$ are of the form $b_iv_i+a_i$, except $J_1(v_1)$ . It is easy to see that all the $b_i$ must be the same ( $J_1$ being non-linear ) in which case
\begin{align}\label{a49}
\sum_{\alpha\neq \beta}J'_{\alpha}J_{\beta}&=\sum_{\alpha=2}^N b \sum_{\beta\neq \alpha}(bv_{\beta}+a_{\beta})=(N-2)b(b\sum v_\beta+\sum a_\beta)\nonumber\\
&\approx (N-2)b(-bv_1+A)
\end{align}
\[
J'_1\sum_{\alpha}J_{\alpha}+J_1\sum J'_{\alpha}\approx J'_1(-bv_1+a)+(N-1)b J_1
\]
 $J(v:=bv_1-a):=J_1 (v_1)$ will have to satisfy
\begin{equation}\label{a50}
-J'v+(N-1)J=(N-2)v,
\end{equation}
hence $J(v)=v+Cv^{N-1}$, i.e.
\begin{equation}\label{a51}
f'^2_1(x_1)\stackrel{!}{=}(b f_1-a)+C(bf_1-a)^{N-1}.
\end{equation}
 As $bf_\alpha+a_\alpha=J_\alpha=f'^2_\alpha (x_\alpha)$ implies
\begin{equation}\label{a52}
f_\alpha=\frac{b}{4}(x_\alpha-t_\alpha)^2 - \frac{a_\alpha}{b}
\end{equation}
the corresponding hypersurface
$\Sigma:=\{\stackrel{\rightharpoonup}{x}\in \mathbb{R}^N, \sum f_\alpha+f_1=0\}$, given by
\begin{equation}\label{a53}
\sum(x_\alpha - t_\alpha)^2 =\frac{4}{b}(-f_1+\frac{a}{b})=-\frac{4}{b^2}(bf_1-a)=:g(x_1)=r^2(x_1),
\end{equation}
 is also given by solving the ODE ( cp.(51) )
\begin{equation}\label{a54}
g'^2=4g\left(-C\left(\frac{-b^2}{4}g\right)^{N-2}-1\right),
\end{equation}
resp.
\begin{equation}\label{a55}
r'^2+1+C\left(\frac{-b^2 r^2}{4}\right)^{N-2}=0.
\end{equation}
It is easy to check that one indeed gets the Weierstrass-function(s) $\wp(x_1)$ as solutions of $(54)_{N=4}$ and the catenoid(s)
$r(z)=\frac{1}{e}\cosh (ez+d)$ as solutions of $(55)_{N=3}$, $(e:=\frac{b}{2}\sqrt{+C})$.\\
Finally, if all $J_i$ are linear ($=b_i v_i +a_i$),

one finds the condition $\sum_{i\neq j} b_j (b_i v_i+a_i)\approx 0$, i.e. , with $B:=\sum b_i$,
\begin{equation}\label{a56}
\begin{array}{l}
Bb_i-b_i^2=const\\\\
\left(\sum a_i\right)B=\sum a_i b_i.
\end{array}
\end{equation}
Solving the quadratic equation,
\begin{equation}\label{a57}
b_i=\frac{B}{2}\pm \sqrt{\frac{B^2}{4}-const.}=\frac{B}{2}\left(1\pm \sqrt{1-c}\right)=b_{\pm},
\end{equation}
it follows that the $b_i$ can take only $ 2 $ different values, say $ r < \dfrac{N}{2} $ times $\frac{B}{2}(1+\sqrt{1-c})$ and $N-r$ times $\frac{B}{2}(1-\sqrt{1-c})$, hence ( summing \ref{a57})
\begin{align}\label{a58}
2=\frac{2}{B}\sum b_i&=r(1+\sqrt{1-c})+(N-r)(1-\sqrt{1-c})\\ \nonumber
&=N-(N-2r)\sqrt{1-c},
\end{align}
which implies $\sqrt{1-c}=\frac{N-2}{N-2r}$, i.e. after scaling ( multiplying with $N-2r$ and dividing by $ \vert B\vert $): ~ $ r $ times $\pm(N-r-1)$ and $(N-r)$ times $ \mp(r-1) $.

\bigskip\textbf{Minimal Surfaces from Rotating Shapes}
\bigskip

What kind of $M-$dimensional objects can be moved such that a higher dimensional minimal surface (in some constant-curvature embedding space; $ \mathbb{R}^{N}, \mathbb{R}^{1, N}, S^N, \ldots $) results?
This question is more or less fully understood for the lowest dimension ($M=1, \mathbb{R}^{N=3}$ being classic; this case can be reduced to $M=0$, i.e. simple point motion generating the 1-dimensional object, e.g. as being the trace of a point on a circle rolling around another circle).
On the other hand, as found more than 2 decades ago \cite{Hoppe95}, the Ansatz

\begin{align}\label{2.1}
\left( x^{\mu}(t,\varphi)\right)_{\mu=0, 1, 2} = \begin{pmatrix}
t &  \\
R(wt) & \vec{u}(\varphi)
\end{pmatrix}
\end{align}

for a minimal surface in $ \mathbb{R}^{1, 2} $, with

\begin{align}\label{2.2}
R(wt)\cdot\vec{u}(\varphi)=\begin{pmatrix}
\cos (wt) & -\sin (wt) \\
\sin (wt) & \cos (wt)
\end{pmatrix}\begin{pmatrix}
u_{1}(\varphi) \\
u_{2}(\varphi)
\end{pmatrix}
\end{align}

describing the (constant angular velocity) rotation of a parametrized planar curve

\begin{align}\label{2.3}
\vec{u}(\varphi)=r(\varphi)\begin{pmatrix}
\cos \theta (\varphi) \\
\sin \theta (\varphi)
\end{pmatrix} ,
\end{align}

and leading to the equations

\begin{align}\label{2.4}
\partial_{\alpha}\left(\dfrac{1}{\sqrt{\vec{u}'^2(1-w^2r^2)+w^2(\vec{u}\times\vec{u}')^2}}\begin{pmatrix}
-\vec{u}'^2 & w(\vec{u}\times\vec{u}') \\
w(\vec{u}\times\vec{u}') & 1-w^2r^2
\end{pmatrix}^{\alpha\beta}\partial_{\beta}x^{\mu} \right)_{\mu=0, 1, 2} =0,
\end{align}

due to (\ref{2.4})$_{\mu=0} $ implying

\begin{align}\label{2.5}
\partial_{\varphi}\left( \underbrace{\dfrac{wr\sin\phi}{\sqrt{1-w^2r^2\cos^2\phi}}}_{=:\gamma_{0}}\right) =0
\end{align}

where $\phi=\sphericalangle(\vec{u}, \vec{u}')$ is the angle between $ \vec{u} $ and $ \vec{u}' $, allows one to conclude that the shape of the curve $ \vec{u} $ is given by the simple equation

\begin{align}\label{2.6}
w^2r^2(1+\gamma \sin^2 \phi)=1,\quad \gamma+1=\frac{1}{\gamma_{0}^2}\ (<1, \text{if }w^2r^2<1)
\end{align}

where $ \gamma $ is a constant of integration. This derivation of the shape is somewhat simpler than the standard one (calculating the mean-curvature from (\ref{2.1}), giving

\begin{align}\label{2.7}
\vec{u}'^2w^2(\vec{u}, A\vec{u}')=(w^2r^2-1)(\vec{u}', A\vec{u}''),
\end{align}

with $
A=\begin{pmatrix}
0 & -1 \\
1 & 0
\end{pmatrix},$ resp. deriving a second-order differential equation for $ \theta $ as a function of $ r $, which (being of first order Bernoulli-type in $ g'=-\frac{\tilde{\theta}'}{w} $) can be linearized, and twice integrated,  the very last step being equivalent to solving (\ref{2.6}) after using that

\begin{align}\label{2.8}
\sin^2\phi=\frac{r^2\theta'^2}{r'^2+r^2\theta'^2}(\varphi)=\frac{w^2r^2g'^2}{1+w^2r^2g'^2}(r),
\end{align}

i.e.

\begin{align}\label{2.9}
g'^2(r)=\frac{1-\text{w}}{\text{w}(\text{w}(1+\gamma)-1)}, (\text{w}:=w^2r^2).
\end{align}

In order to compare with the corresponding Euclidean calculation (see e.g. \cite{Nitsche89}) one could rewrite (\ref{2.1}) by substituting

\begin{align}\label{2.10}
t=\frac{v}{w}+g(u(\varphi))
\end{align}

and then notice that

\begin{align}\label{2.11}
R(wt)\vec{u}(\varphi)=R(v)R(wg(u(\varphi)))\left| \vec{u}\right| \begin{pmatrix}
\cos\theta(\varphi) \\
\sin\theta(\varphi)
\end{pmatrix}=R(v)u\begin{pmatrix}
1 \\
0
\end{pmatrix}
\end{align}

if one defines $g$ to undo the rotation $R(\theta)$ by which $\begin{pmatrix}
	\cos\theta \\
	\sin\theta
\end{pmatrix} $ results from $\begin{pmatrix}
1 \\
0
\end{pmatrix} $ -- and then choosing $\varphi=r=\left| \vec{u}\right| =u$.
The class of solutions (of (\ref{2.6}), resp. (\ref{2.4}), resp. (\ref{2.9})) considered in \cite{Hoppe95} were $ (n\neq k, nk>0) $

\begin{align}\label{2.12}
\vec{u}(\varphi)=\frac{1}{2n}\begin{pmatrix}
\cos n\varphi \\
\sin n\varphi
\end{pmatrix} +\frac{1}{2k}\begin{pmatrix}
\cos k\varphi \\
\sin k\varphi
\end{pmatrix}
\end{align}

for which, with

\begin{align*}
c&:=\cos\left( \frac{n-k}{2}\varphi\right), s=\sin\left( \frac{n-k}{2}\varphi\right)\\
w^2(\vec{u}^2)&=1+\frac{4nk}{(n-k)^2}c^2
\end{align*}

\begin{align}\label{2.13}
w^2=\frac{4n^2k^2}{(n-k)^2}, \ \vec{u}'^2=c^2, \sin^2\phi=\frac{c^2(1+\frac{4nk}{(n-k)^2})}{1+\frac{4nk}{(n-k)^2}+c^2}.
\end{align}

Note that they can be written ( for later convenience ) in the form

\begin{align}\label{2.14}
2\vec{u}_{\pm}(\varphi)=\left( \frac{1}{n}R^n(\varphi)\pm\frac{1}{k}R^k(\varphi)\right) \begin{pmatrix}
1 \\
0
\end{pmatrix},
\end{align}

making the calculations leading to (\ref{2.13}) very simple, just using

\begin{align*}
R'(\varphi)&=AR(\varphi)\\
R^{T}&=R^{-1}\\
A&=\begin{pmatrix}
0&-1 \\
1&0
\end{pmatrix}=-A^T=-A^{-1};
\end{align*}

e.g. $2\vec{u}'_{\pm}=A(R^{n}\pm R^{k})\begin{pmatrix}
1 \\
0
\end{pmatrix}$

\begin{align}\label{2.15}
4\vec{u}'^2&=(1,0)(R^{-n}\pm R^{-k})(-A^2)(R^n\pm R^k)\begin{pmatrix}
1 \\
0
\end{pmatrix}\nonumber\\
&=(1,0)(1+1\pm R^{n-k}\pm R^{k-n})\begin{pmatrix}
1 \\
0
\end{pmatrix}\nonumber\\
&=2(1\pm \cos((n-k)\varphi))=4c^2,
\end{align}

resp. $ 4s^2 $. The above given frequency $w=\frac{2nk}{k-n}$  is special to (\ref{2.12}), alone for the following reason:


\begin{align}\label{2.16}
\vec{x}(t,\varphi)&:=R(wt)\vec{u}(\varphi)\nonumber\\
&=\left( \frac{1}{2n}R^n\left( \varphi+\frac{w}{n}t\right) +\frac{1}{2k}R^k\left( \varphi+\frac{w}{k}t\right) \right)\begin{pmatrix}
1 \\
0
\end{pmatrix} \nonumber\\
&\stackrel{!}{=}\left( \frac{1}{2n}R^n\left( \tilde{\varphi}+\tilde{t}\right) +\frac{1}{2k}R^k\left( \tilde{\varphi}-\tilde{t}\right) \right)\begin{pmatrix}
1 \\
0
\end{pmatrix}\nonumber\\
&=:\tilde{\vec{x}}\left(\tilde{t}=t,\tilde{\varphi}:=\varphi-\frac{n+k}{n-k}t\right)=\frac{1}{2}\vec{v}_{+}(\tilde{\varphi}+\tilde{t})-\frac{1}{2}\vec{v}_{-}(\tilde{\varphi}-\tilde{t})
\end{align}

is of the same type as $\vec{u}$ (only $ \varphi\rightarrow\tilde{\varphi}+\tilde{t} $, $ \varphi\rightarrow\tilde{\varphi}-\tilde{t} $ by a reparametrization of $ \varphi $ alone (!), in the 2 terms) and $ \begin{pmatrix}
t \\
\tilde{\vec{x}}
\end{pmatrix}  $ a difference  of 2 Null-curves in $\mathbb{R}^{1,2}$.

\begin{align}\label{2.17}
\begin{pmatrix}
t \\
\tilde{\vec{x}}
\end{pmatrix}=\frac{1}{2}\begin{pmatrix}
\varphi+t=\varphi_{+}\\
\vec{v}_{+}(\varphi_{+})
\end{pmatrix} -\frac{1}{2}\begin{pmatrix}
\varphi-t=\varphi_{-}\\
\vec{v}_{-}
\end{pmatrix}.
\end{align}

While a lot is known about Null curves in relation with minimal surfaces in $ \mathbb{R}^{1,2} $, and their Weierstraß representations, the crucial question is whether (and if yes, how) any of these structures can be also used for $ M>1 $. (The $\varphi_{\pm}$ decomposition certainly is special to $ M+1=2 $).
The most appealing seems to be that the technical simplifications following from (\ref{2.14}) are matched by the important geometric property of the Epicycloids (having $\left|n-k\right|$ cusps) as being obtained by rolling circles (with one point marked) around circles.
Consider (see e.g.[``Epicycloids'', Wikipedia]) rolling a circle of radius $a$ around a circle of radius $b$ (centered at $\vec{0}$); then

\begin{align}\label{2.18}
\vec{x}(\psi)&=R(\psi)\underbrace{R_{\vec{c}:=\vec{a}+\vec{b}}\left(\frac{b}{a}\psi \right)\vec{b}}_{\vec{c}+R\left( \frac{b}{a}\psi\right) (\vec{b}-\vec{c})}\nonumber\\
&=R(\psi)\left(\vec{a}+\vec{b}-\underbrace{R(\psi)\left(\frac{b}{a}\psi \right)}_{R\left( \left( 1+\frac{b}{a}\right) \psi\right) }\vec{a}\right)\nonumber\\
&=(a+b)\begin{pmatrix}
\cos \psi \\
\sin \psi
\end{pmatrix} -a\begin{pmatrix}
\cos(1+\frac{b}{a}) \psi \\
\sin(1+\frac{b}{a}) \psi
\end{pmatrix}\nonumber\\
&\stackrel{\psi=\chi a}{=}(a+b)\begin{pmatrix}
\cos \chi a \\
\sin \chi a
\end{pmatrix} -a\begin{pmatrix}
\cos(a+b) \chi \\
\sin(a+b) \chi
\end{pmatrix}\\
&=a(a+b)\left\lbrace \frac{1}{a}\begin{pmatrix}
\cos \chi a \\
\sin \chi a
\end{pmatrix}-\frac{1}{a+b}\begin{pmatrix}
\cos \chi (a+b)\\
\sin \chi(a+b)
\end{pmatrix} \right\rbrace. \nonumber
\end{align}

For $a=n$ and $ a+b=k $ this is proportional to

\begin{align}\label{2.19}
\vec{u}_{-}=\frac{1}{2n}\begin{pmatrix}
c_n \\
s_n
\end{pmatrix} -\frac{1}{2k}\begin{pmatrix}
c_k \\
s_k
\end{pmatrix}
\end{align}

(i.e. not (\ref{2.12}), but with a relative sign); choosing $ \chi=\varphi+\frac{\pi}{n-k} $ (however) one obtains

\begin{align}\label{2.20}
\vec{u}_{-}\left( \varphi+\frac{\pi}{n-k}\right) &=\left( \frac{1}{2n}R^n(\varphi)R\left( \frac{n\pi}{n-k}\right) -\frac{1}{2k}R^k(\varphi)R\left( \frac{k\pi}{n-k}\right)\right)  \begin{pmatrix}
1 \\
0
\end{pmatrix}\nonumber \\
&=R\left( \frac{n\pi}{n-k}\right) \left\lbrace\left( \frac{1}{2n}R^n(\varphi)+ \frac{1}{2k}R^k(\varphi)\right)  \right\rbrace \begin{pmatrix}
1 \\
0
\end{pmatrix}\nonumber\\
&=R\left( \frac{n\pi}{n-k}\right) \vec{u}_{+}(\varphi).
\end{align}

The shape equation (\ref{2.6}), $ w^2r^2(1+\gamma\sin^2\phi)=1 $, for the curve $ \vec{u}(\varphi)=\begin{pmatrix}
u_1(\varphi) \\
u_2(\varphi)
\end{pmatrix}=r(\varphi)\begin{pmatrix}
\cos\theta(\varphi) \\
\sin\theta(\varphi)
\end{pmatrix}  $, $ \phi:=\sphericalangle(\vec{u},\vec{u}') $, can for $ r'\neq 0\neq r $ be written as

\begin{align}\label{2.21}
\frac{1}{w^2r^2}-1=\gamma \sin^2\phi=\gamma\frac{r^2\theta'^2(\varphi)}{r'^2+r^2\theta'^2}=\gamma\frac{r^2\tilde{\theta}'^2(r)}{1+r^2\tilde{\theta}'^2(r)},
\end{align}

when taking $r$ to parametrized the curve, i.e. using $ \theta(\varphi)=\tilde{\theta}(r(\varphi))\Rightarrow\theta'=\frac{d\tilde{\theta}}{dr}\cdot r $, implying

\begin{align*}
r^2\tilde{\theta}'^2=\frac{w-1}{1-\delta w},
\end{align*}

with $ \delta:=1+\gamma=\frac{1}{\gamma_{0}^2} $ and $ \text{w}:=w^2r^2, $ i.e.

\begin{align}\label{2.22}
\pm \int d\tilde{\theta}=\int \sqrt{\frac{\text{w}-1}{1-\delta \text{w}}}\frac{dr}{r}=\frac{1}{2\sqrt{\delta}}\int\frac{d\text{w}}{\text{w}}\sqrt{\frac{\text{w}-1}{\gamma_{0}^2-\text{w}}}=:\gamma_{0}J.
\end{align}

Calculating $J$ with the substitution

\begin{align}\label{2.24}
v:=\sqrt{\frac{\text{w}-1}{\gamma_{0}^2-\text{w}}}=\left| \frac{\sin\psi}{\cos\psi}\right|=\tan\psi,\ \psi \in \left( 0,\frac{\pi}{2}\right),
\end{align}

with

\begin{align*}
\text{w}=\gamma_{0}\sin^2\psi+\cos^2\psi,
\end{align*}

one obtains

\begin{align}\label{2.23}
\pm(\tilde{\theta}(r)-\theta_0)=\gamma_{0}\arctan v-\arctan(\gamma_{0}v),
\end{align}

resp.

\begin{align*}
\pm \tan (\tilde{\theta}-\theta_0)&=\frac{\tan(\gamma_{0}\psi)-\gamma_{0}\tan\psi}{1+\gamma_{0}\tan\psi+\tan(\gamma_{0}\psi)}\\
&=\frac{(a-b)\sin_0\cos-(a+b)\sin\cos_0}{(a-b)\cos_0\cos+(a+b)\sin\sin_0},\ \gamma_{0}=\frac{a+b}{a-b},\  \frac{b}{a}>0,
\end{align*}
with $\sin:=\sin\left(\frac{a-b}{2}\varphi\right), \cos:=\cos\left(\frac{a-b}{2}\varphi\right)$, resp. $\sin_0:=\sin\left(\frac{a+b}{2}\varphi\right), \cos_0:=\cos\left(\frac{a+b}{2}\varphi\right)$,
\begin{align}\label{2.25}
=\frac{b\sin a\varphi-a\sin b\varphi}{b\cos a\varphi-a\cos b\varphi},
\end{align}

which is! the tangent of $ \theta $ of

\begin{align}\label{2.26}
\lambda\vec{u}(\varphi)=\frac{1}{2a}\begin{pmatrix}
\cos a\varphi \\
\sin a\varphi
\end{pmatrix} - \frac{1}{2b}\begin{pmatrix}
\cos b\varphi \\
\sin b\varphi
\end{pmatrix},
\end{align}

cp (\ref{2.18}).




\newpage
\secct{Minimal Tori in $S^3$, Stiefel Manifolds, and Determinantal Varieties} \label{sec-III}
 One way of formulating the problem of finding minimal (hyper) surfaces in spheres (rather than in $\mathbb{R}^N$) is to subject the usual parametric area functional
 
\begin{equation}\label{b1}
\begin{array}{ll}
A[\stackrel{\rightharpoonup}{x}]=\int \sqrt{g} ~
\mathrm{d}^{M}\vphi\\\\
g=det(\partial_a \stackrel{\rightharpoonup}{x}\partial_b \stackrel{\rightharpoonup}{x})_{a,b =1, \ldots, M}
\end{array}
\end{equation}

to the constraint  $\stackrel{\rightharpoonup}{x}^2(\vphi)=1$, i.e. to consider

\begin{equation}\label{b2}
S[\stackrel{\rightharpoonup}{x}]:=A[\stackrel{\rightharpoonup}{x}]+\frac{1}{2}\int \mathrm{d}^{M}\varphi ~~\lambda(\varphi)(\stackrel{\rightharpoonup}{x}^2-1),
\end{equation}
whose stationary points $\stackrel{\rightharpoonup}{x}(\varphi)$ are then easily seen to satisfy
\begin{equation}\label{b3}
\Delta \stackrel{\rightharpoonup}{x}:=\frac{1}{\sqrt{g}}\partial_a\sqrt{g} g^{ab}\partial_b \stackrel{\rightharpoonup}{x}=\lambda \stackrel{\rightharpoonup}{x}=-M \stackrel{\rightharpoonup}{x},
\end{equation}
where the last equality in (\ref{b3}) easily follows by multiplying the inner part with $\stackrel{\rightharpoonup}{x}$ ~ ( and noting that $\stackrel{\rightharpoonup}{x}^2(\varphi)=1$ implies
$\stackrel{\rightharpoonup}{x} \partial_a \stackrel{\rightharpoonup}{x}=0$ and, differentiating again,
$\partial_b \stackrel{\rightharpoonup}{x} \partial_a \stackrel{\rightharpoonup}{x}+\stackrel{\rightharpoonup}{x} \partial^2_{ab} \stackrel{\rightharpoonup}{x}=0$,
hence $\stackrel{\rightharpoonup}{x}\Delta\stackrel{\rightharpoonup}{x}=g^{ab}\stackrel{\rightharpoonup}{x} \partial^2_{ab}\stackrel{\rightharpoonup}{x}=-M$).\\

The celebrated clifford torus CT being
\begin{equation}\label{b4}
\stackrel{\rightharpoonup}{x}^T(\varphi^1 \varphi^{M=2})=\frac{1}{\sqrt{2}}( \cos \varphi^1, \sin\varphi^1, \cos\varphi^2, \sin\varphi^2)
\end{equation}
let us try (as in \cite{Arnlind}) to find solutions as graphs over CT, i.e. of the form
\begin{equation}\label{b5}
\stackrel{\rightharpoonup}{x}(\varphi^1\varphi^2)=\left(\begin{array}{l}
\cos \theta \cos \varphi^1\\
\cos \theta \sin \varphi^1\\
\sin \theta \cos \varphi^2\\
\sin \theta \sin \varphi^2\\
\end{array}\right)=:\left(\begin{array}{l}
cc_1\\cs_1\\sc_2\\ss_2
\end{array}\right).
\end{equation}

\begin{equation}\label{b6}
\begin{array}{l}
\partial_1 \stackrel{\rightharpoonup}{x}=c
\left(\begin{array}{c}
-s_1\\c_1\\0\\0
\end{array}\right)+\theta_1
\left(\begin{array}{ll}
{}& c_1\\
-s &{}\\
{}& s_1\\\\
{} &c_2\\
c & {}\\
{}&s_2
\end{array}\right)=c \stackrel{\rightharpoonup}{e_1}+\theta_1\stackrel{\rightharpoonup}{e}
\\\\
\partial_2 \stackrel{\rightharpoonup}{x}=s
\left(\begin{array}{c}
0\\0\\-s_2\\c_2
\end{array}\right)+\theta_2 \stackrel{\rightharpoonup}{e}, ~~~\theta_a:=\frac{\partial \theta}{\partial \varphi^a}
\end{array}
\end{equation}
imply
\begin{equation}\label{b7}
g_{ab}=\left(\begin{array}{cc}
c^2+\theta_1^2 & \theta_1\theta_2\\\\
\theta_1\theta_2 & s^2+\theta_2^2
\end{array}\right),
\end{equation}
and the Ansatz $\stackrel{\rightharpoonup}{m}=\alpha_1 \stackrel{\rightharpoonup}{e_1}+
\alpha_2\stackrel{\rightharpoonup}{e_2}+\alpha\stackrel{\rightharpoonup}{e}$,
automatically orthogonal to $\stackrel{\rightharpoonup}{x}$, and easily giving
\begin{equation}\label{b8}
\stackrel{\rightharpoonup}{m} || sc\stackrel{\rightharpoonup}{e}-s\theta_1 \stackrel{\rightharpoonup}{e_1}
-c\theta_2 \stackrel{\rightharpoonup}{e_2},
\end{equation}
when requiring orthogonality with (\ref{b6}), implies that $h_{ab}:=\frac{\stackrel{\rightharpoonup}{m}\partial^2_{ab}x}{\vert
\stackrel{\rightharpoonup}{m}\vert}$ will be proportional to

\begin{equation}\label{b9}
\left(\begin{array}{cc}
s^2 c^2 +sc\theta_{11}+2s^2 \theta_1^2 ~~& sc\theta_{12}+(s^2-c^2)\theta_1\theta_2\\\\
sc\theta_{12}+(s^2-c^2)\theta_1\theta_2 & -s^2 c^2 +sc\theta_{22}-2c^2 \theta_2^2
\end{array}\right)
\end{equation}
so that  $g^{ab} h_{ab}=0$ ( which is the non-trivial content of (\ref{b3}),as itscomponents in the direction(s) of $\partial_1 \stackrel{\rightharpoonup}{x}$, $ \partial_2 \stackrel{\rightharpoonup}{x} $ and $ \stackrel{\rightharpoonup}{x} $ are trivially satisfied) corresponds to
\begin{align}\label{b10}
(s^2+\theta_2^2)(s^2c^2+sc\theta_{11}+2s^2\theta_1^2)&+(c^2+\theta_1^2)(-s^2c^2+sc\theta_{22}
-2c^2\theta_2^2)\nonumber\\&-2\theta_1\theta_2(sc\theta_{12}+(s^2-c^2)\theta_1\theta_2)=0,
\end{align}
also following directly by varying (cp. (\ref{b7}), (\ref{b2}), (\ref{b5}))
\begin{equation}\label{b11}
S[\theta]:=\int\sqrt{c^2s^2+s^2\theta_1^2+c^2\theta_2^2}d\varphi^1 d\varphi^2.
\end{equation}
For $\theta=\theta (t:=k\varphi^1+l\varphi^2)$, (\ref{b10}) reduces to the ( highly non-linear ) ODE
\begin{align}\label{b12}
sc \{ sc  \stackrel{..}{\theta}(k^2 s^2+l^2 c^2)&+ \stackrel{.}{\theta}^2 \left[
(l^2 - k^2)s^2c^2+2s^4k^2 - 2c^4 l^2\right] \\ \nonumber
&+ s^2c^2(s^2-c^2) \} =0
\end{align}
corresponding to stationary points of (cp. (\ref{b11}), the overall sign put in for later convenience )
\begin{align}\label{b13}
\int \left(L:=-\sqrt{c^2s^2+(k^2s^2+l^2c^2)\stackrel{.}{\theta}^2}\right) dt.
\end{align}

Switching now to physical terminology, interpreting $\theta$ as the t(ime)-dependent position of a `particle' moving in some `potential' one identifies a conserved quantity, resp. integration of (\ref{b12}), via the `Hamiltonian'
\begin{equation}\label{b14}
K:=\stackrel{.}{\theta} \frac{\partial L}{\partial \stackrel{.}{\theta}}-L=\ldots=\frac{s^2c^2}{\sqrt{c^2s^2+
(k^2s^2+l^2c^2)\stackrel{.}{\theta}^2}}= const =: E,
\end{equation}
i.e.
\begin{equation}\label{b15}
\stackrel{.}{\theta}^2+\frac{c^2s^2}{k^2s^2+l^2c^2}\left(1-\frac{c^2s^2}{E^2}\right)=0,
\end{equation}
with the second term an `effective potential' (of the mass $\frac{1}{2}$ particle ) for a `zero-energy' solution ( with respect to which the positive integration constant $E$ should perhaps better be denoted by  $\kappa^2$, to reflect the somewhat dangerous double interpretation of `Energy') as the rhs (\ref{b15}) and, expressing (\ref{b14}) in terms of $\pi:=\frac{\partial L}{\partial \stackrel{.}{\theta}}$,
\begin{equation}\label{b16}
\vert \sin \theta\cos \theta\vert\sqrt{1-\frac{\pi^2}{k^2s^2+l^2c^2}}.
\end{equation}

That (\ref{b15})is a consequence of (\ref{b12}) ( the reverse is trivially verified ) can of course also be checked directly ( without referring to physics terminology ): (\ref{b12}) is, after dividing by $s^2c^2 =: r(\theta)$,  and multiplying by  $2 \stackrel{.}{\theta}$, of the form
($L=-\sqrt{r + f \stackrel{.}{\theta}^2}$,
$f(\theta):=k^2 s^2 + l^2 c^2$) ~~~$ 2f \stackrel{..}{\theta} \stackrel{.}{\theta}+f' \stackrel{.}{\theta}^3=\stackrel{.}{\theta}r'+2f \frac{r'}{r}\stackrel{.}{\theta}^3$,
\begin{align}\label{b17}
\stackrel{.}{\theta}\left(r'+2(f\stackrel{.}{\theta}^2)\frac{r'}{r}\right)=\left(f\stackrel{.}{\theta}^2\right)^.
=\frac{d}{dt}\left(G(\theta(t))\right),
\end{align}
giving the first order ODE
\[
G'=r'+2G\frac{r'}{r}, ~ resp
\]
\[
r^2\left(\frac{G}{r^2}\right)'=r'
\]
whose solution is
$\frac{G}{r^2}=-\frac{1}{r}+const.$, i.e.
\begin{equation}\label{b18}
Cr^2 -r=G(\theta(t))=f\stackrel{.}{\theta}^2
\end{equation}
with $C=\frac{1}{E^2}>0$, as the right hand side is manifestly non-negative.\\
While the case $k=0$ (resp. $l=0$),
solvable with the help of elliptic integrals, is well discussed in the differential geometry community (cp. \cite{Brendle2013}) the case
$k=l\neq 0$ can be solved in terms of elementary functions as fallows: with
\[
\alpha(\varphi:=\varphi^1+\varphi^2):=2\theta\left(k(\varphi^1+\varphi^2)\right),
\]
$a:=\frac{1}{2E}$,
(\ref{b15}) reads
\begin{equation}\label{b19}
\stackrel{.}{\alpha}^2+\sin^2\alpha\left(1-a^2\sin^2\alpha\right)=0.
\end{equation}
Assuming $a=\cosh \gamma\geq 1$, ($E\leq \frac{1}{2}$) in order to have turning-points, i.e.  $\stackrel{.}{\alpha}=0$ for some $t$ ) the particle ($\alpha$) oscrillates between
$\alpha_{-}=\arcsin 2E\underset{(=)}{<} \frac{\pi}{2}$, and
$\alpha_{+}=\pi-\arcsin 2E$, while direct integration of  (\ref{b19}), yields
\begin{align}\label{b20}
\varphi-\varphi_0=\pm\int\frac{d\alpha}{\sin \alpha\sqrt{a^2 \sin^2\alpha-1}}=\mp \arctan\frac{\cos \alpha}{\sqrt{a^2\sin^2\alpha-1}}
\end{align}

\begin{equation}\label{b21}
\begin{array}{l}
\sin \alpha(\varphi)=\frac{1}{\sqrt{\cos^2(\varphi-\varphi_0)+\cosh^2\gamma\sin^2(\varphi-\varphi_0)}}\\\\
\cos\alpha(\varphi)=\frac{-\sinh \gamma\sin (\varphi-\varphi_0)}{ \sqrt{\cos^2(\varphi-\varphi_0)+\cosh^2\gamma\sin^2(\varphi-\varphi_0)}}
\end{array}
\end{equation}

i.e. a one parameter class  ($e:=\sinh \gamma$) of `minimal' (extremal) tori in $S^3$,
\begin{equation}\label{b22}
\stackrel{\rightharpoonup}{x}_e(\varphi^1\varphi^2)=\frac{1}{\sqrt{2}}\left(
\begin{array}{ll}
{}&\cos \varphi^1\\
\sqrt{1-\frac{e\sin (\varphi-\varphi_0)}{\sqrt{1+e^2\sin^2}}}&{}\\
{}& \sin \varphi^1\\\\
{}&\cos\varphi^2\\
\sqrt{1+\frac{e\sin (\varphi-\varphi_0)}{\sqrt{1+e^2\sin^2}}}&{}\\
{}&\sin\varphi^2
\end{array}\right).
\end{equation}
While it is easy to see that they are without intersections ( i.e. embedded) the consequence, namely that they must be congment to (\ref{b4}), is stuming.\\
As (\ref{b20}) appears identically in the equation for geodesics on $S^2$ ( i.e. great circles ), and the Hopf map ( s.b.) applied to (\ref{b22}) gives
\begin{equation}\label{b23}
\left(\begin{array}{c}
\cos 2\theta\\\\
\sin 2\theta \cos(\varphi^2 - \varphi^1)\\\\
\sin 2\theta \sin (\varphi^2 -\varphi^1)
\end{array}\right)
\end{equation}
( and the signs of $k=\pm l$ were never used above, so that one could as well have defined $ t $ to be  $\varphi^2-\varphi^1$)(rather than $\varphi^2+\varphi^1$ for which one would need a `conjugate' Hopf map) one may view the solutions (\ref{b22}) as a one parameter clars of inverse images of great circles on $S^2$.\\
The problem however rests in the vast freedom in the construction, which becomes apparent when trying to fix the details, using quaternions:
\begin{align}\label{b24}
q=(q_0, q_1, q_2, q_3)\hat{=}\hat{q}=
\left(\begin{array}{cc}
q_0-iq_3~~& -iq_1 -q_2\\
-iq_1+q_2 & q_0+iq_3
\end{array}\right)
\cong q_0+iq_1+j q_2+k q_3,
\end{align}
using either matrix multiplications for $\hat{q}$, or
\begin{equation}\label{b25}
(t,\stackrel{\rightharpoonup}{x})\cdot(s, \stackrel{\rightharpoonup}{y})=(ts-\stackrel{\rightharpoonup}{x}\stackrel{\rightharpoonup}{y}, t\stackrel{\rightharpoonup}{y}+s\stackrel{\rightharpoonup}{x}+\stackrel{\rightharpoonup}{x}\times\stackrel{\rightharpoonup}{y}),
\end{equation}
resp. $i^2=j^2=k^2=-1$, $ij=k$, ... .\\

Defining an anti automorphism via
\begin{equation}\label{b26}
\tilde{q}:=q_0-iq_1+jq_2+hq_3
\end{equation}
the Hopf-map from $S^3$ to $S^2$  can be given as
\begin{equation}\label{b27}
\pi(q):=\tilde{q}q,
\end{equation}
nicely fitting with the action of $S^3$ onto itself by right-multiplication ($q\to qr$), and
\begin{equation}\label{b28}
\pi(qr)=\tilde{r}\pi(q)r
\end{equation}
defining an action of  $S^3$ on $S^2$ ( $q_1\equiv 0$).
In coordinates one finds
\begin{equation}\label{b29}
\left(\begin{array}{c}
t'\\y'\\z'
\end{array}\right)
=
\left(\begin{array}{c}
t^2+x^2-y^2-z^2\\
2(tz-xy)\\
2(ty+xz)
\end{array}\right)
\end{equation}
for $\pi(q)$ and for (\ref{b28})
\begin{equation}\label{b30}
\left(\begin{array}{c}
t''\\y''\\z''
\end{array}\right)
=
\left(\begin{array}{ccc}
1-2(r_2^2+r_3^2) ~& 2(r_1 r_3-r_0 r_2)~&-2(r_0r_3+r_1r_2\\
2(r_0r_2+r_1r_3)& 1-2(r_1^2+r_2^2)&
2(r_0r_1-r_2r_3)\\
2(r_0 r_3-r_1r_2)&-2(r_0r_1+r_2r_3)&
1-2(r_1^2+r_3^2)
\end{array}\right)
\left(\begin{array}{c}
t'\\y'\\z'
\end{array}\right)
\end{equation}
i.e. written as an ordinary SO(\ref{b3}) transformation on $S^2$.
As $\pi(e^{i\rho}q)=\pi(q)$, and for the great circle
\begin{equation}\label{b31}
\left(\begin{array}{c}
\cos 2\theta\\\sin 2\theta\cos \\ \sin 2\theta\sin
\end{array}\right)
=
\left(\begin{array}{ccc}
\cos\beta~& 0~&\sin\beta\\
0 & 1 & 0\\
-\sin\beta & 0& \cos\beta
\end{array}\right)
\left(\begin{array}{c}
0\\ \cos\tilde{\varphi} \\ \sin\tilde{\varphi}
\end{array}\right),
\end{equation}
with the constant angle $ \beta $ related to the parameter $e$ in (\ref{b22}), while a rotation
$\left(\begin{array}{ccc}
1~& 0~&0\\
0 & \cos\tilde{\gamma} & -\sin\tilde{\gamma}\\
0 & \sin\tilde{\gamma}& \cos\tilde{\gamma}
\end{array}\right)$,
with possibly non-trivial $\tilde{\gamma}(\varphi^1\varphi^2)$, geometrically leaves the equator invariant, one finds $ 3 $ different transformations to possibly act on the ordinary Clifford-torus,
\begin{equation}\label{b32}
\underset{\sim}{\stackrel{\rightharpoonup}{x}^T}(\tilde{\varphi^1},\tilde{\varphi^2})=\frac{1}{\sqrt{2}}(\cos\tilde{\varphi^1}, \sin \tilde{\varphi^1}, \cos \tilde{\varphi^2}, \sin\tilde{\varphi^2}),
\end{equation}
namely:\\
$\cos \rho +i\sin\rho $ from the left ( as that does not change $\pi(q)$), $(-\sin\tilde{\gamma}+i\cos\tilde{\gamma})$ from the right-leaving the equation invariant, while transforming
$(c_1  s_1  c_2  s_2)^T$ to
\begin{equation}\label{b33}
\left(\begin{array}{cccc}
-\sin\tilde{\gamma}~& 0~ & 0~& \cos\tilde{\gamma}\\
0 & -\sin\tilde{\gamma} & \cos\tilde{\gamma} & 0\\
0 & -\cos\tilde{\gamma}& -\sin\tilde{\gamma} & 0\\
\cos\tilde{\gamma}& 0& 0& -\sin\tilde{\gamma}
\end{array}\right)
\left(\begin{array}{c}
c_1\\s_1\\c_2\\s_2
\end{array}\right),
\end{equation}
as well as $\cos\frac{\beta}{2}-\sin\frac{\beta}{2}\cdot k$ from the right (corresponding to (\ref{b31}). Leaving out (\ref{b33}) ( i.e. choosing
$-\sin\tilde{\gamma}=1$) ) for the moment, one would find
\begin{equation}\label{b34}
\cos\rho \left(\begin{array}{c}
c_1c+s_2s\\cs_1-sc_2\\cc_2+ss_1\\cs_2-sc_1
\end{array}\right)+
\sin\rho \left(\begin{array}{c}
sc_2-cs_1\\cc_1+ss_2\\sc_1-cs_2\\cc_2+ss_2
\end{array}\right)
\stackrel{!}{=}
\sqrt{2}
\left(\begin{array}{ll}
{}&c_1\\
\cos\theta&{}\\
{}& s_1\\\\
{}&c_2\\
\sin\theta&{}\\
{}&s_2
\end{array}\right),
\end{equation}
which is ` almost ' ( but not quite) solvable, (note that even if it was, the question why $\rho(\varphi^1\varphi^2)$ does not destroy minimality would still have to be answered).\\
Note that the right-action $q\to qr$ may also be written as
\begin{equation}\label{b35}
q\to
\left(\begin{array}{cc}
r_0~ & -\stackrel{\rightharpoonup}{r}^T
\\\\\stackrel{\rightharpoonup}{r}&(r_0  \mathbf{1}- \stackrel{\rightharpoonup}{r}\times)
\end{array}\right)
\left(\begin{array}{c}
q_0\\\\\stackrel{\rightharpoonup}{q}.
\end{array}\right)
\end{equation}
 Instead of trying to match the freedom in the quaternion description of the Hopf-map relation between geoderics on$S^2$ and minimal surfaces in $S^3$ let us note that in [ ACH13 ] an explicit reparametrisation,
\begin{equation}\label{b36}
\tilde{\varphi^1}=\varphi^1+\int^{\varphi}u=\varphi^1+f_1(\varphi),~~
\tilde{\varphi^2}=\varphi^2+\int^{\varphi}v=\varphi^2+f_2(\varphi),
\end{equation}
 was constructed that transforms the metric induced from (\ref{b32}),  $\tilde{g}_{ab}=\frac{1}{2}\delta_{ab}$, into the metric induced from (\ref{b22}), (cp. (\ref{b7}), (\ref{b22}), (\ref{b15}).\\
$c=\cos\theta(\varphi)$, $s=\sin\theta(\varphi)$; $\frac{1}{2E}=\sqrt{1+e^2}$)
\begin{align}\label{b37}
g_{ab}&\stackrel{(\ref{b7})}{=}\left(\begin{array}{cc}
c^2+\stackrel{.}{\theta}^2 ~~& \stackrel{.}{\theta}^2\\\\
\stackrel{.}{\theta}^2 & s^2+\stackrel{.}{\theta}^2
\end{array}\right)
\stackrel{(\ref{b15})}{=}\left(
\begin{array}{cc}
c^4+\frac{s^4c^4}{E^2} ~~& s^2c^2\left(\frac{s^2c^2}{E^2}-1\right)\\\\
s^2c^2\left(\frac{s^2c^2}{E^2}-1\right) & s^4+\frac{s^4c^4}{E^2}
\end{array}\right)\nonumber\\
&=
\left(\begin{array}{cc}
g+c^4 ~~&g- s^2c^2\\\\
g-s^2c^2& g+s^4
\end{array}\right)=J^T(\tilde{g}_{ab})J=\frac{1}{2}J^T J\nonumber\\
&=\frac{1}{2} \left(\begin{array}{cc}
1+2u+u^2+v^2~~&u+v+u^2+v^2\\\\
u+v+u^2+v^2& 1+2v+u^2+v^2
\end{array}\right)\nonumber\\
&=\frac{1}{2} \left(\begin{array}{cc}
(1+u)^2+v^2~~&u(1+u)+v(1+v)\\\\
u(1+u)+v(1+v)& (1+v)^2+u^2
\end{array}\right),
\end{align}
\\
\\
\begin{equation}\label{b38}
J:=\left(\dfrac{\partial \tilde{\varphi}^a}{\partial {\varphi}^b}
\right)\stackrel{(\ref{b36})}{=}
\left(\begin{array}{cc}
1+u~~&u\\\\
v& 1+v
\end{array}\right),
\end{equation}
from which

\begin{equation}\label{b39}
u=-s^2+\frac{s^2c^2}{E}=\sqrt{g} -s^2, ~~v=-c^2+\frac{s^2c^2}{E}=\sqrt{g} -c^2,
\end{equation}
follows. Note that
$u-v=c^2-s^2$, and
\begin{equation}\label{b40}
\sqrt{g}=\sqrt{s^2c^2+\stackrel{.}{\theta}^2}=\frac{s^2c^2}{E}.
\end{equation}

Let us now show that also the second fundamental forms are transformed into each other under the transformation (\ref{b36}), with $u$ and $v$ given by (\ref{b39}), resp.
\begin{equation}\label{b41}
\begin{array}{l}
u=\frac{1}{2}\left(1+\frac{e\sin}{\sqrt{1+e^2\sin^2}} \right) \left[ \sqrt{1+e^2}(
1-\frac{e\sin}{ \sqrt{1+e^2\sin^2}})-1 \right] \\\\
~~=\frac{1}{2}\left[-1-\frac{e\sin}{\sqrt{1+e^2\sin^2}}+\frac{\sqrt{1+e^2}}{1+e^2\sin^2}\right] \\\\
v=\frac{1}{2}\left(1-\frac{e\sin}{\sqrt{1+e^2\sin^2}}\right)\left[\sqrt{1+e^2}\left(
1+\frac{e\sin}{\sqrt{1+e^2\sin^2}}\right)-1\right] \\\\
~~=\frac{1}{2}\left[-1+\frac{e\sin}{\sqrt{1+e^2\sin^2}}+\frac{\sqrt{1+e^2}}{1+e^2\sin^2}\right] \\\\
\sin=\sin( \varphi-\varphi^\circ ).
\end{array}
\end{equation}
 Namely, using (\ref{b9}), (\ref{b12}) and  (\ref{b15}) ($k=l=1$ for simplicity) one finds
\begin{equation}\label{b42}
h_{ab}=\frac{s^2c^2}{E}
\left(\begin{array}{cc}
2c^2~~&c^2-s^2\\\\
c^2-s^2& -2s^2
\end{array}\right)=\sqrt{g}
\left(\begin{array}{cc}
2c^2~~&c^2-s^2\\\\
c^2-s^2& -2s^2
\end{array}\right)
\end{equation}
for the minimal surfaces (\ref{b22}). Note that due to
$-4s^2c^2-(c^2-s^2)^2=-1$, $h=-g$ ( hence $K=-1, R=0$) is manifest, just as $g^{ab}h_{ba}=0$.\\
Transforming on the other hand the second fundamental form of (\ref{b32})
\begin{equation}\label{b43}
\tilde{h}_{ab}=-\frac{1}{2}
\left(\begin{array}{cc}
-1~~&0\\\\
0& +1
\end{array}\right),
\end{equation}
choosing the orientation via  $\stackrel{\rightharpoonup}{n}:=\frac{1}{\sqrt{2}}
\left(\begin{array}{c}
-\tilde{c}_1\\
-\tilde{s}_1\\
+\tilde{c}_2\\
+\tilde{s}_2
\end{array}\right)
$
one finds
\begin{equation}\label{b44}
\begin{array}{ll}
h_{ab}:&= -\frac{1}{2}J^T \left(\begin{array}{cc}
-1& 0\\
0 & +1
\end{array}\right) J\\\\
{}& =-\frac{1}{2}\left(\begin{array}{cc}
v^2-(u+1)^2& v(v+1)-u(u+1)\\
v(v+1)-u(u+1) & (v+1)^2-u^2
\end{array}\right)\\\\
{}& =+\frac{1}{2}\left(\begin{array}{cc}
(u+1)^2-v^2& u(u+1)-v(v+1)\\
u(u+1)-v(v+1) & u^2-(v+1)^2
\end{array}\right)\\\\
{}&={ \frac{1}{2}\left(\begin{array}{cc}
(\sqrt{g}+c^2)^2-(\sqrt{g}-c^2)^2&(\sqrt{g}-s^2)(\sqrt{g}+c^2)-(\sqrt{g}-c^2)(\sqrt{g}+s^2)\\\\
(\sqrt{g}-s^2)(\sqrt{g}+c^2)-(\sqrt{g}-c^2)(\sqrt{g}+s^2) & (\sqrt{g}-s^2)^2-(\sqrt{g}+s^2)^2
\end{array}\right)}\\\\
{}&=\sqrt{g}\left(\begin{array}{cc}
2c^2& (c^2-s^2)\\
(c^2-s^2) & -2s^2
\end{array}\right).
\end{array}
\end{equation}
With both first and second fundamental forms coinciding, insertion of (\ref{b36}) into (\ref{b32}), yielding
\begin{align}\label{b45}
\stackrel {\rightharpoonup}{ \underset{\sim}{x}}&=\cos \varphi_1 \frac{1}{\sqrt{2}}\left(
\begin{array}{c}
\cos f_1 \\\sin f_1\\ \cos(f_2+\varphi)\\\sin(f_2+\varphi)
\end{array}\right)
+
\sin \varphi_1 \frac{1}{\sqrt{2}}\left(
\begin{array}{c}
-\sin f_1 \\\cos f_1\\ \sin(f_2+\varphi)\\-\cos(f_2+\varphi)
\end{array}\right)\\ \nonumber
\\ \nonumber
&=\cos\varphi_1 \stackrel {\rightharpoonup}{ \underset{\sim}{e_{1}}}(\varphi)
+\sin\varphi_1 \stackrel {\rightharpoonup}{ \underset{\sim}{e_{2}}}(\varphi),
\end{align}
should differ from (\ref{b22}), resp. ($\varphi^2=\varphi-\varphi^1$)
\begin{align}\label{b45}
\stackrel {\rightharpoonup}{x}&=\cos \varphi_1 \left(
\begin{array}{c}
\cos\theta(\varphi) \\ 0\\ \sin \theta\cos \varphi\\\sin \theta\sin \varphi
\end{array}\right)
+
\sin \varphi_1\left(
\begin{array}{c}
0\\\cos \theta(\varphi)\\ \sin \theta\sin \varphi\\-\sin \theta\cos \varphi
\end{array}\right)\\ \nonumber
\\ \nonumber
&=\cos\varphi_1 \stackrel {\rightharpoonup}{e_{1}}(\varphi)
+\sin\varphi_1 \stackrel {\rightharpoonup}{e_{2}}(\varphi),
\end{align}
only by a fixed ($\varphi$-independent) orthogonal transformation $S$, i.e. should hold that
\begin{equation}\label{b47}
 \stackrel {\rightharpoonup}{ \underset{\sim}{e_i}}(\varphi)=S \stackrel {\rightharpoonup}{e}_i(\varphi).
\end{equation}
 While at first glance hard to believe (as e.g. implying that the $ 4 $ components of  $\stackrel {\rightharpoonup}{ \underset{\sim}{e_{1}}}(\varphi)$, resp.
$\stackrel {\rightharpoonup}{ \underset{\sim}{e_{2}}}(\varphi)$, must be linearly dependent \footnote{I am grateful to J.Arnlind, M.Bordemann and B.Durhuus for helpful discussions when trying to resolve this puzzle} i.e. define a
$\varphi$-independent hyperplane -one of the components of $\stackrel {\rightharpoonup}{e_{1}}(\varphi)$, resp.
$\stackrel {\rightharpoonup}{e}_2(\varphi)$, being zero) help comes from the (proof of) the fundamental theorem for surfaces, which ( given equality of first and second fundamental forms ) constructs $\stackrel {\rightharpoonup}{x}$, resp. $\stackrel {\rightharpoonup}{ \underset{\sim}{x}}$, from given initial conditions. Choosing for simplicity ( and without lose of generality ) $\varphi^\circ=0$ and $\varphi_1^\circ=0=\varphi_2^\circ$ as well as the integration constants defining $f_1$ and $f_2$ such that $f_1(0,0)=0=f_2(0,0)$, and noting that (cp. (\ref{b22})) $\theta_0=\theta(0,0)=\frac{\pi}{4}$ as well as $\stackrel{.}{\theta}_0:=\stackrel{.}{\theta}(0,0)=\frac{e}{2}$,  whereas \\
 $ f_1'(0,0) = u_0 = u(0,0) = \frac{1}{2}\left(\sqrt{1+e^2}-1 \right) = v(0,0) = v_0 = f_2'(0,0) $, one finds

 \begin{align}
 \begin{array}{l}
 \stackrel {\rightharpoonup}{x}_0:=\stackrel {\rightharpoonup}{x}(0,0)=
 \stackrel {\rightharpoonup}{ \underset{\sim}{x}}(0,0)= \stackrel {\rightharpoonup}{ \underset{\sim}{x_0}}
 =\frac{1}{\sqrt{2}}\left(\begin{array}{c}
 1\\0\\1\\0
 \end{array}\right)=: \stackrel {\rightharpoonup}{ \underset{\sim}{n_0}}= \stackrel {\rightharpoonup}{n}_0\\\\
  \stackrel {\rightharpoonup}{x}_+:=\left((\partial_1+\partial_2) \stackrel {\rightharpoonup}{x}\right)(0,0)
  =\frac{1}{\sqrt{2}}\left(\begin{array}{c}
 -e\\1\\e\\1
 \end{array}\right)=:\sqrt{1+e^2}\frac{1}{\sqrt{2}}\left(\begin{array}{c}
-\sin \epsilon\\\cos \epsilon\\\sin \epsilon\\\cos \epsilon
 \end{array}\right)=:\sqrt{1+e^2}\stackrel {\rightharpoonup}{n}_+\\\\
 \end{array}. \nonumber
 \end{align}
 \begin{align}
 \begin{array}{l}
 \stackrel {\rightharpoonup}{ \underset{\sim}{x}}_+:=\left((\partial_1+\partial_2)  \stackrel {\rightharpoonup}{ \underset{\sim}{x}}\right)(0,0)=\sqrt{1+e^2}\frac{1}{\sqrt{2}}\left(\begin{array}{c}
0\\1\\0\\1
 \end{array}\right)=:\sqrt{1+e^2}\stackrel {\rightharpoonup}{ \underset{\sim}{n}}_+\\\\
 \stackrel {\rightharpoonup}{x}_-:=\left((\partial_1-\partial_2)  \stackrel {\rightharpoonup}{x}\right)(0,0)=\frac{1}{\sqrt{2}}\left(\begin{array}{c}
0\\1\\0\\-1
 \end{array}\right)=:\stackrel {\rightharpoonup}{n}_-\\\\
  \stackrel {\rightharpoonup}{ \underset{\sim}{x}}_-:=\left((\partial_1-\partial_2)  \stackrel {\rightharpoonup}{\underset{\sim}{x}}\right)(0,0)=:\stackrel {\rightharpoonup}{\underset{\sim}{n}}_-=\stackrel {\rightharpoonup}{n}_-
 \end{array}.
 \end{align}
 Hence one is looking for an orthogonal transformation ($R$) leaving $\stackrel {\rightharpoonup}{n}_0$ and $\stackrel {\rightharpoonup}{n}_-$ fixed, transforming
$\stackrel {\rightharpoonup}{ \underset{\sim}{n}}_+$ into $\stackrel {\rightharpoonup}{n}_+$, i.e.
\begin{equation}\label{b49}
R \stackrel {\rightharpoonup}{ \underset{\sim}{n}}_+=\stackrel {\rightharpoonup}{n}_+=
\cos \epsilon \stackrel {\rightharpoonup}{ \underset{\sim}{n}}_++\sin\epsilon
\left(\begin{array}{c}
-1\\0\\1\\0
 \end{array}\right) \frac{1}{\sqrt{2}}=: \stackrel {\rightharpoonup}{ \underset{\sim}{n}},
\end{equation}
hence ( choosing $R$ to be a $\epsilon-$ rotation with ($ \stackrel {\rightharpoonup}{ \underset{\sim}{n}}_+,  \stackrel {\rightharpoonup}{ \underset{\sim}{n}})$ -plane )
\begin{equation}\label{b50}
R \stackrel {\rightharpoonup}{ \underset{\sim}{n}}=-\sin \epsilon  \stackrel {\rightharpoonup}{ \underset{\sim}{n}}_++\cos\epsilon  \stackrel {\rightharpoonup}{ \underset{\sim}{n}}.
\end{equation}

Instead of directly verifying $(\ref{b47})_{(S=R)}$ it is instructive to consider the action of $R$ on (\ref{b32}), resp.

\begin{equation}\label{b51}
\begin{array}{ll}
\stackrel {\rightharpoonup}{ \underset{\sim}{x}}(\tilde{\varphi^1}, \tilde{\varphi^2})&=\frac{1}{2}
(\cos\tilde{\varphi^1}+\cos\tilde{\varphi^2})
\left(\begin{array}{c}
1\\0\\1\\0
 \end{array}\right) \frac{1}{\sqrt{2}}-\frac{1}{2}(\cos\tilde{\varphi^1}-\cos\tilde{\varphi^2})
\left(\begin{array}{c}
-1\\0\\1\\0
 \end{array}\right)\frac{1}{\sqrt{2}}\\\\
 {}&+\frac{1}{2}(\sin\tilde{\varphi^1}+\sin\tilde{\varphi^2})
\left(\begin{array}{c}
0\\1\\0\\1
 \end{array}\right)\frac{1}{\sqrt{2}}+
 \frac{1}{2}(\sin\tilde{\varphi^1}-\sin\tilde{\varphi^2})
\left(\begin{array}{c}
0\\1\\0\\-1
 \end{array}\right)\frac{1}{\sqrt{2}}
 \\\\
{}&= \frac{1}{2}\left\{ (\tilde{c}_1+\tilde{c}_2)
 \stackrel {\rightharpoonup}{ \underset{\sim}{n_0}}-
 (\tilde{c}_1-\tilde{c}_2)\stackrel {\rightharpoonup}{ \underset{\sim}{n}}+
 (\tilde{s}_1+\tilde{s}_2)  \stackrel {\rightharpoonup}{ \underset{\sim}{n_+}}+ (\tilde{s}_1-\tilde{s}_2)  \stackrel {\rightharpoonup}{ \underset{\sim}{n_-}}\right\}
\end{array}
\end{equation}
i.e.
\begin{equation}\label{b52}
\begin{array}{ll}
R \stackrel {\rightharpoonup}{ \underset{\sim}{x}}&=\frac{1}{2}\left\{ (\tilde{c}_1+\tilde{c}_2)
 \stackrel {\rightharpoonup}{ \underset{\sim}{n_0}}-
  (\tilde{c}_1-\tilde{c}_2)\left( \cos\epsilon \stackrel {\rightharpoonup}{ \underset{\sim}{n}}
 - \sin\epsilon \stackrel {\rightharpoonup}{ \underset{\sim}{n_+}} \right) \right.
\\\\
{}& \left. ~ + ( \tilde{s}_1+\tilde{s}_2) \left( \cos\epsilon \stackrel {\rightharpoonup}{ \underset{\sim}{n_+}}
 + \sin\epsilon \stackrel {\rightharpoonup}{ \underset{\sim}{n}} \right)+
( \tilde{s}_1-\tilde{s}_2 )  \stackrel {\rightharpoonup}{ \underset{\sim}{n_-}} \right\},
\end{array}
\end{equation}

which is supposed to equal (\ref{b22}),
\begin{equation}\label{b53}
\frac{1}{2}\left(\begin{array}{cccc}
1+c~~ & -s~~& 1-c~~& -s\\\\
s & 1+c & -s & -1+c \\\\
1-c & s & 1+c & s\\\\
s & c-1 & -s & 1+ c
\end{array}\right)
\left(\begin{array}{c}
\tilde{c}_1\\ \tilde{s}_1 \\ \tilde{c}_2 \\ \tilde{s}_2
\end{array}\right)\frac{1}{\sqrt{2}}
\stackrel{!}=
\left(\begin{array}{c}
\cos\theta c_1\\ \cos\theta s_1 \\ \sin\theta c_2\\ \sin\theta s_2
\end{array}\right),
\end{equation}
the orthogonal $ 4 \times 4 $ matrix on the left ($c=\cos\epsilon=\frac{1}{\sqrt{1+e^2}}, s=\sin\epsilon=\frac{e}{\sqrt{1+e^2}}$)
being the looked for ( cp. (\ref{b47})) matrix $S(=R)$. Simple inversion
\begin{equation}\label{b54}
\frac{1}{\sqrt{2}}\left(\begin{array}{c}
\tilde{c}_1\\ \tilde{s}_1 \\ \tilde{c}_2 \\ \tilde{s}_2
\end{array}\right)
=
\frac{1}{2}\left(\begin{array}{cccc}
1+c~~ & s~~& 1-c~~& s\\\\
-s & 1+c & s & c-1 \\\\
1-c & -s & 1+c & -s\\\\
-s &-1+c & s & 1+ c
\end{array}\right)
\left(\begin{array}{c}
\cos\theta c_1\\ \cos\theta s_1 \\ \sin\theta
\left(\begin{array}{c}
c_1\cos\varphi+s_1\sin\varphi\\
c_1\sin\varphi-s_1\cos\varphi
\end{array}\right)
\end{array}\right)
\end{equation}
then, writing the l.h.s. as
\begin{equation}\label{b55}
\frac{1}{\sqrt{2}}\left(\begin{array}{cc}
\begin{array}{cc}
\cos f_1~~&-\sin f_1\\
\sin f_1 & \cos f_1
\end{array} & 0\\\\
0 &\left(\begin{array}{cc}
\cos(f_2+\varphi) ~~& +\sin(f_2+\varphi)\\
\sin(f_2+\varphi) & -\cos (f_2+\varphi)
\end{array}\right)
\end{array}\right)
\left(\begin{array}{c}
c_1\\ s_1\\ c_1\\ s_1
\end{array}\right)
\end{equation}
turns out to be consistently equivalent to
\begin{equation}\label{b56}
\begin{array}{l}
\sqrt{2}\cos f_1=(1+c)\cos\theta +(1-c)\sin \theta\cos\varphi+s\sin\theta\sin\varphi\\\\
\sqrt{2}\sin f_1=-s\cos\theta +s\sin \theta\cos\varphi-(1-c)\sin\theta\sin\varphi\\\\
\sqrt{2}\cos(f_2+\varphi)=(1-c)\cos\theta +\sin \theta[(1+c)\cos\varphi-s\sin\varphi]\\\\
\sqrt{2}\sin(f_2+\varphi)=-s\cos\theta +\sin \theta[(1+c)\sin\varphi+s\cos\varphi].
\end{array}
\end{equation}
 As the sum of the first $ 2 $ squares ( crucially using $e\tan 2\theta\sin\varphi=-1$ ) gives indeed $ 2 $,
 one could in principle forget eqs (\ref{b37})-(\ref{b44}) and simply define $f_1$ and $f_2$ (cp. (\ref{b36})) by (\ref{b56}).\\
 The hyperplane property for
$\stackrel {\rightharpoonup}{ \underset{\sim}{e_i}}(\varphi)$ finally follows, as (\ref{b56}) implies
\begin{align*}
s\cos(f_2+\varphi)+(1-c)\sin(f_2+\varphi)&=\sqrt{2}s\sin\theta\cos\varphi\\
&=s\cos f_1+(1+c)\sin f_1
\end{align*}

\bigskip
Pedestrian\footnote{an elegant, earlier, proof \cite{Chunpong} has been communicated to me by J.Choe.} level-set proof that the \textbf{Stiefel-manifolds} $ \sum_{n,k} $ are minimal in the corresponding sphere $ S^{nk-1}(\sqrt{k}), $ resp. the cone $ \hat{\sum}_{n,k} $minimal in $ \mathbb{R}^{nk} $: consider $ k $ orthogonal $ n- $ vectors $ \vec{x}^{(1)}, \vec{x}^{(2)}, \ldots, \vec{x}^{(n)} \in \mathbb{R}^{n},$ of equal length, i.e. the constraints $ (a,b=1\cdots k) $

\begin{eqnarray}\label{3.1}
	u^{(a\neq b)}=u^{\left[ ab\right] }:=\vec{x}^{(a)}\vec{x}^{(b)}\stackrel{!}{=}0 \left( \#=\frac{k(k-1)}{2}\right) \\
	\label{3.2}
	v^{(a)}:=\frac{1}{2}\left( \left( \vec{x}^{(a)}\right) ^{2}-\left( \vec{x}^{(a+1)}\right) ^{2}\right)\stackrel{!}{=}0 \begin{pmatrix}
	\#=k-1 \\
	a=1\cdots k-1
	\end{pmatrix}
\end{eqnarray}

defining $ \hat{\sum}_{n,k}\subset\mathbb{R}^{nk} (dim=nk-(k-1)(1+\frac{k}{2})=nk-\frac{k(k+1)}{2}+1) $.

\begin{align}
Area (\hat{\sum}_{n,k})=\int d^{nk}x\prod_{a=1}^{k-1}\delta\left( v^{(a)}(x)\right) \prod_{a<b}\delta\left(u^{\left[ab \right]}(x)\right)\sqrt{\det M},
\end{align}

where $ M $ is the matrix formed by the scalar products of the gradients of the $ K=\frac{1}{2}(k-1)(k+2)=\frac{k(k+1)}{2}-1 $ constraints $ W^{(A)}_{A=1\cdots K} $ $\left( \text{i.e. } M_{AB}:=\vec{\nabla}W^{(A)}\cdot\vec{\nabla} W^{(B)}\right) $.
The minimality condition then reads (see e.g. \cite{Gee90})

\begin{align}\label{3.3}
\delta\left( \prod_{B}W^{(B)}\right) \cdot Tr(P\cdot\partial^{2}W^{(A)})\stackrel{!}{=}0, \quad A=1\cdots K
\end{align}

where $ P $ is the projector onto the subspace orthogonal to all the $ \nabla W^{(A)} $ (and projects onto the tangent-space of $ \hat{\sum}_{n,k} $ when $ \chi:=\prod_{A}W^{(A)}=0 $),

\begin{align}\label{3.4}
P_{ij}=\delta_{ij}-\partial_{i}W^{(A)}(M^{-1})^{AB}\partial_{j}W^{(B)}.
\end{align}

For $ \hat{\sum}_{n,k}, $ resp. (\ref{3.1})/(\ref{3.2}) one can calculate (\ref{3.4}), resp. $\hat{P}:=P|_{\chi=0}$ explicitly, and prove that (\ref{3.3}) is satisfied.

Simple(st) non-trivial example: $ n=3, k=2 $

\begin{align*}
u^{\left[12\right]}\left( \vec{\mathfrak{x}}\right) &=u:=x_{1}x_{4}+x_{2}x_{5}+x_{3}x_{6}\left( =\vec{x}\cdot\vec{y}\right) \stackrel{!}{=}0\\
v^{(1)}\left( \vec{\mathfrak{x}}\right)&=v:=\frac{1}{2}\left( x_{1}^{2}+x_{2}^{2}+x_{3}^{2}-x_{4}^{2}-x_{5}^{2}-x_{6}^{2}\right) =\frac{1}{2}\left(\vec{x}^{2}-\vec{y}^{2}\right) \stackrel{!}{=}0
\end{align*}

(2 orthogonal vectors of equal length)

\begin{align*}
\vec{\nabla}u&=\begin{pmatrix}
\vec{y} \\
\vec{x}
\end{pmatrix} , \quad \vec{\nabla}v=\begin{pmatrix}
\vec{x} \\
-\vec{y}
\end{pmatrix} \quad \left( \vec{\nabla}\vec{u}\right) ^{2}=\left( \vec{\nabla}\vec{v}\right) ^{2}=\underbrace{r^{2}}_{=\vec{x}^{2}+\vec{y}^{2}}=2s^{2}, \quad \vec{\nabla}u\cdot \vec{\nabla}v=0\\
M&=r^{2}\mathds{1},\quad M^{-1}=\frac{1}{r^{2}}\mathds{1}\\
P&=\begin{pmatrix}
\mathds{1} & 0 \\
0 & \mathds{1}
\end{pmatrix} - \frac{1}{r^{2}} \left(
\begin{array}{ccc|ccc}
y_{1}^{2} & & y_{\alpha}y_{\beta} &x_{1}y_{1}&&y_{\alpha}x_{\beta}\\
&y_{2}^{2} & &&x_{2}y_{2} \\
&& y_{3}^{2} & &&x_{3}y_{3} \\ \hline
&&& x_{1}^{2}&&x_{\alpha}x_{\beta}\\
&&&&x_{2}^{2}\\
&&&&&x_{3}^{2}
\end{array}\right)\\
&-\frac{1}{r^{2}}\left(
\begin{array}{ccc|ccc}
x_{1}^{2} & & x_{\alpha}x_{\beta} &-x_{1}y_{1}&&-x_{\alpha}y_{\beta}\\
&x_{2}^{2} & &&-x_{2}y_{2} \\
&& x_{3}^{2} & &&-x_{3}y_{3} \\ \hline
&&& y_{1}^{2}&&y_{\alpha}y_{\beta}\\
&&&&y_{2}^{2}\\
&&&&&y_{3}^{2}
\end{array}\right)\\
\partial^{2}u&=\begin{pmatrix}
0 & \mathds{1} \\
\mathds{1} & 0
\end{pmatrix}\quad \partial^{2}v=\begin{pmatrix}
\mathds{1} &  \\
 & -\mathds{1}
\end{pmatrix} \quad\delta(v)Tr(P\partial^{2}u)=0, \ \delta(u)Tr(P\partial^{2}v)=0
\end{align*}

even without $\delta(u)\cdot\delta(v)$

\begin{align*}
K=\frac{2\cdot 3}{2}-1=2, \quad dim\hat{\Sigma}_{3,2}=6-2=4, \quad dim\Sigma_{3,2}=3.
\end{align*}

Another example \footnote{special thanks to G.Linardopoulos and T.Turgut for discussions going back to 2016/17}: $ n=3, k=3 $

\begin{align*}
u_{3}&=u^{[12]}=\vec{x}\cdot\vec{y}, \ u_{1}=u^{[23]}=\vec{y}\cdot\vec{z}, \ u_{2}=u^{[13]}=\vec{x}\cdot\vec{z}\\
v^{(1)}&=\frac{1}{2}\left( \vec{x}^2-\vec{y}^2\right), \ v^{(2)}=\frac{1}{2}\left( \vec{y}^2-\vec{z}^2\right)\\
\end{align*}
\begin{align*}
\vec{\nabla}u^{[12]}&=\begin{pmatrix}
\vec{y} \\
\vec{x} \\
0
\end{pmatrix} ,\ \vec{\nabla}u^{[23]}=\begin{pmatrix}
0 \\
\vec{z} \\
\vec{y}
\end{pmatrix} , \ \vec{\nabla}u^{[13]}=\begin{pmatrix}
\vec{z} \\
0 \\
\vec{x}
\end{pmatrix}\\
\vec{\nabla}v^{(1)}&=\begin{pmatrix}
\vec{x} \\
-\vec{y} \\
0
\end{pmatrix}, \ \vec{\nabla}v^{(2)}=\begin{pmatrix}
0 \\
\vec{y} \\
-\vec{z}
\end{pmatrix}\\
K&=5, \ dim\hat{\sum_{\ 3,3}}=4\subset\mathbb{R}^{9},\ dim\Sigma_{\ 3,3}=3 \subset S^{8}
\end{align*}

on $ \hat{\sum}_{3,3}, $ with $ s^{2}=\frac{r^{2}}{3}=(\vec{x}^{2}=\vec{y}^{2}=\vec{z}^{2}) $, and $\hat{M}:=M|_{\chi=0}$

\begin{align*}
\hat{M}&=s^{2}\left(
\begin{array}{ccc|cc}
2&&&\\
&2&&0&\\
&&2\\ \hline
0&&&2&-1\\
&&&-1&2
\end{array}\right), \ \hat{M}^{-1}=\frac{1}{s^{2}}\left(
\begin{array}{ccc|cc}
\frac{1}{2}&&&\\
&\frac{1}{2}&&0&\\
&&\frac{1}{2}\\ \hline
0&&&\frac{2}{3}&\frac{1}{3}\\
&&&\frac{1}{3}&\frac{2}{3}
\end{array}\right)=\widehat{M^{-1}}\\
\partial^{2}u_{3}&=\begin{pmatrix}
0 & \mathds{1} & 0 \\
\mathds{1} & 0 & 0 \\
0 & 0 & 0
\end{pmatrix}, \ \partial^{2}u_{1}=\left(
\begin{array}{c|cc}
0&0&0\\ \hline
0&0&\mathds{1}\\
0&\mathds{1}&0
\end{array}\right)\\
\partial^{2}v_{1}&=\begin{pmatrix}
\mathds{1} &  &  \\
& -\mathds{1} &  \\
&  & 0
\end{pmatrix}, \ \partial^{2}v_{2}=\begin{pmatrix}
0 & 0 & 0 \\
0 & \mathds{1} & 0 \\
0 & 0 & -\mathds{1}
\end{pmatrix} \\
\hat{P}&=\begin{pmatrix}
\mathds{1} &  &  \\
& \mathds{1} &  \\
&  & \mathds{1}
\end{pmatrix}-\frac{1}{2s^{2}}\left[\left(
\begin{array}{c|c|c}
y_{\alpha} y_{\beta} & y_{\alpha} x_{\beta} & 0 \\ \hline
x_{\alpha} y_{\beta} & x_{\alpha} x_{\beta} & 0 \\ \hline
0 & 0 & 0
\end{array}\right) + \left(\begin{array}{c|c|c}
0 & 0 & 0 \\ \hline
0 & z_{\alpha} z_{\beta} & z_{\alpha} y_{\beta} \\ \hline
0 & y_{\alpha} z_{\beta} & y_{\alpha} y_{\beta}
\end{array}\right) \right. \\
& \left. + \left(\begin{array}{c|c|c}
z_{\alpha} z_{\beta} & 0 & z_{\alpha} x_{\beta} \\ \hline
0 & 0 & 0 \\ \hline
x_{\alpha} z_{\beta} & 0 & x_{\alpha} x_{\beta}
\end{array}\right) \right] - \frac{1}{3s^2} \left(
\begin{array}{c|c|c}
0 & x_{\alpha} y_{\beta} & -x_{\alpha} z_{\beta} \\ \hline
0 & -y_{\alpha} y_{\beta} & y_{\alpha} z_{\beta} \\ \hline
0 & 0 & 0
\end{array}\right) - \frac{1}{3s^2} \left(
\begin{array}{c|c|c}
0 & 0 & 0 \\ \hline
y_{\alpha} x_{\beta} & -y_{\alpha} y_{\beta} & 0 \\ \hline
-z_{\alpha} x_{\beta} & z_{\alpha} y_{\beta} & 0
\end{array}\right) \\
& - \frac{2}{3s^2} \left(
\begin{array}{c|c|c}
x_{\alpha} x_{\beta} & -x_{\alpha} y_{\beta} & 0 \\ \hline
-y_{\alpha} x_{\beta} & y_{\alpha} y_{\beta} & 0 \\ \hline
0 & 0 & 0
\end{array}\right) - \frac{2}{3s^2} \left(
\begin{array}{c|c|c}
0 & 0 & 0 \\ \hline
0 & y_{\alpha} y_{\beta} & -y_{\alpha} z_{\beta} \\ \hline
0 & -z_{\alpha} y_{\beta} & z_{\alpha} z_{\beta}
\end{array}\right) \\
&\left\langle k=4:\begin{pmatrix}
2 & -1 & 0 \\
-1 & 2 & -1 \\
0 & -1 & 2
\end{pmatrix}^{-1}=\frac{1}{4}\begin{pmatrix}
3 & 2 & 1 \\
2 & 4 & 2 \\
1 & 2 & 4
\end{pmatrix}   \right\rangle
\end{align*}

To verify (\ref{3.3}) for the general case (any $k\leq n\in\mathbb{N}$) is completely straightforward, except for the fact that $ \vec{\nabla}v^{(a)}\cdot \vec{\nabla}v^{(a+1)}\neq 0. $

\begin{align*}
\hat{M}&=s^{2}\left(
\begin{array}{c|ccccc}
2\cdot \mathds{1}\\ \hline
&2&-1&&&0\\
&-1&2&-1\\
&&-1&2&-1\\
&&&&\ddots&\\
&0&&&-1&2
\end{array}\right)
\end{align*}

however has the explicit inverse

\begin{align*}
\hat{M}^{-1}&=\frac{1}{ks^{2}}\left(
\begin{array}{c|c}
\frac{k}{2}\cdot \mathds{1}&0\\ \hline
0&Q
\end{array}\right)=\widehat{M^{-1}}
\end{align*}

with the $ (k-1)\times(k-1) $ matrix $Q$ having the matrix elements

\begin{align}\label{3.5}
Q_{a'b'}=min(a',b')\cdot k-a'b'.
\end{align}

With

\begin{align*}
\partial_{ci}u^{(a<b)}&=\delta_{ca}x_{bi}+\delta_{cb}x_{ai}\\
\partial_{ci}v^{(a')}&=\delta_{ca'}x_{a'i}-\delta_{ca'+1}x_{a'+1,i}
\end{align*}

the projector $P$ becomes:

\begin{align}\label{3.6}
\hat{P}_{ci,dj}&=\delta_{cd}\delta_{ij}-\frac{1}{2s^2}\sum(\delta_{ca}x_{bi}+\delta_{cb}x_{ai})(\delta_{da}x_{bj}+\delta_{db}x_{aj})\nonumber\\
&-\frac{1}{ks^2}\sum^{k-1}_{a,b=1}(\delta_{ca}x_{ai}-\delta_{c,a+1}x_{a+1,i})(k\cdot min(a,b)-a\cdot b)(\delta_{db}x_{bj}-\delta_{d,b+1}x_{b+1,j}).
\end{align}

When acting on

\begin{align*}
\left( \partial^{2}u^{(a<b)}\right) _{ci,dj}=\delta_{ij}(\delta_{ca}\delta_{db}+\delta_{cb}\delta_{da})
\end{align*}

\begin{align}\label{3.7}
\left( \partial^{2}v^{(e)}\right) _{ci,dj}=\delta_{ij}(\delta_{ce}\delta_{de}+\delta_{c,e+1}\delta_{d,e+1})
\end{align}

and then taking the trace,
the only non-trivial part is the action of the last term in (\ref{3.6}) on (\ref{3.7}), which is proportional to
something that on the constraint manifold vanishes ( here calculated/displayed for non-boundary e ),
\begin{align*}
&\frac{1}{ks^2}\vec{x}^{2}_{e}\left[(ek-e^2)+(k(e-1)-(e-1)^2)-2(k(e-1)-(e-1)e)\right]\\
&-\frac{1}{ks^2}\vec{x}^{2}_{e+1}\left[(k(e+1)-(e+1)^2)+(ek-e^2)-2(ek-e(e+1))\right] \approx 0. \\
\end{align*}

\bigskip\textbf{Minimality of Determinantal Varieties}
\bigskip

Consider\footnote{Thanks to J. Choe for several discussions, in particular for having raised the question, whether non-quadratic matrices can give rise to minimal submanifolds just as quadratic ones do.}

\begin{align}
\label{new1}
\Sigma_{p>q}:=\left\lbrace A=\left( \vec{a}_1\vec{a}_2\cdots\vec{a}_q\right) \mid\vec{a}_i\in\mathbb{R}^p, rank A=q-1\right\rbrace ;
\end{align}

as the linear dependence of the vectors $\vec{a}_i$ may, after permutation, be taken to be of the form

\begin{align}
\label{new2}
\vec{a}_q=\lambda_1\vec{a}_1+\cdots+\lambda_{q-1}\vec{a}_{q-1},
\end{align}

clearly

\begin{align}
\label{new3}
\dim\Sigma_{p>q}=p(q-1)+(q-1)=(p+1)(q-1)<pq-1,
\end{align}

i.e. $\Sigma_{p>q}$ is of codimension $>1$. 

Let us prove minimality for the simplest non trivial example,

\begin{align}
\label{new4}
\Sigma_{3,2}=\left\lbrace A=\left( \vec{a}_1,\vec{a}_2\right) \mid\vec{a}_1,\vec{a}_2\in\mathbb{R}^3, rank A=1\right\rbrace 
\end{align}

in 2 different ways:

1. Taking

\begin{align}
\label{new5}
x(a^1,a^2,a^3,\lambda):=\begin{pmatrix}
\vec{a}=a^{\alpha}\vec{v}_{\alpha} \\ 
\lambda\vec{a}
\end{pmatrix}, 
\end{align}

with $\left\lbrace \vec{v}_{\alpha}\right\rbrace _{\alpha=1,2,3}$ being an orthonormal basis of $\mathbb{R}^{3}$, one has

\begin{equation}
\label{new6}
\begin{split}
\partial_{\alpha}x&=\frac{\partial x}{\partial a^{\alpha}}=\begin{pmatrix}
\vec{v}_{\alpha} \\ 
\lambda\vec{v}_{\alpha}
\end{pmatrix} , \quad\partial_{\lambda}x=\partial_{4}x=\begin{pmatrix}
0\\
\vec{a}
\end{pmatrix}\\
G_{AB}&=\begin{pmatrix}
(1+\lambda^2)\mathds{1}_{3\times3} & \lambda\vec{a} \\ 
\lambda\vec{a}^T & \vec{a}^T\vec{a}	\end{pmatrix}, \quad  G^{AB}=\begin{pmatrix}
\frac{1}{1+\lambda^2}\left( \mathds{1}+\lambda^2\frac{\vec{a}\vec{a}^T}{\vec{a}^2}\right) & -\frac{\lambda\vec{a}}{\vec{a}^2} \\ 
-\frac{\lambda\vec{a}^T}{\vec{a}^2} & \frac{1+\lambda^2}{\vec{a}^2}
\end{pmatrix}, 
\end{split}
\end{equation}

and, defining 2 orthogonal normal vectors

\begin{align}
\label{new7}
n=\frac{1}{\sqrt{1+\lambda^2}}\begin{pmatrix}
\lambda\vec{e}\\
-\vec{e}
\end{pmatrix},\quad n'=\frac{1}{\sqrt{1+\lambda^2}}\begin{pmatrix}
\lambda\vec{e}\, '\\
-\vec{e}\, '
\end{pmatrix},
\end{align}

where $\vec{e}$ and $\vec{e}\, '$ are 2 orthonormal vectors orthogonal to $\vec{a}$,

\begin{align}
	\label{new8}
	H_{AB}=\frac{1}{\sqrt{1+\lambda^2}}\begin{pmatrix}
	0_{3\times 3} & -\vec{e} \\ 
	-\vec{e}^T & 0
	\end{pmatrix} ,\quad H_{AB}'=\frac{1}{\sqrt{1+\lambda^2}}\begin{pmatrix}
	0_{3\times 3} & -\vec{e}\, ' \\ 
	-\vec{e}\, '^T & 0
	\end{pmatrix}
	\end{align}
	
	as
	
	\begin{align}
	\label{new9}
	\partial_{\alpha\beta}^{2}x=0,\quad\partial_{\lambda}^2x=0,\quad\partial_{\alpha\lambda}^2x=\begin{pmatrix}
	0\\
	\vec{v}_{\alpha}
	\end{pmatrix};
	\end{align}
	
	hence
	
	\begin{align}
	\label{new10}
	H\sim G^{AB}H_{AB}=0, \quad H'\sim G^{AB}H'_{AB}=0.
	\end{align}
	
	\newpage
	2. Using the Singular Value Decomposition
	
	\begin{align}
	\label{new11}
	A=\sigma\, \vec{u}\, \overrightarrow{v}^T,
	\end{align}
	
	where $\sigma^2$, $\vec{u}(\theta, \varphi)\in S^2\subset\mathbb{R}^3$, and $\overrightarrow{v}(\psi)\in S^1\subset \mathbb{R}^2$ are the non-zero eigenvalue and unit eigenvector of $AA^T$ resp $A^TA$ one has 
	
	\begin{equation}
	\label{new12}
	\begin{split}
	x&=\begin{pmatrix}
	\sigma\cos\psi\, \vec{u}\\
	\sigma\sin\psi\, \vec{u}
	\end{pmatrix}, \quad \vec{u}=\begin{pmatrix}
	\sin\theta\cos\varphi\\
	\sin\theta\sin\varphi\\
	\cos\theta
	\end{pmatrix}\\
	\partial_{\sigma}x&=\begin{pmatrix}
	c\vec{u}\\
	s\vec{u}
	\end{pmatrix},\quad
	\partial_{\theta}x=\sigma\begin{pmatrix}
	c\,\partial_{\theta}\vec{u}\\
	s\,\partial_{\theta}\vec{u}
	\end{pmatrix}\\
	\partial_{\varphi}x&=\sigma\begin{pmatrix}
	c\, \partial_{\varphi}\vec{u}\\
	s\, \partial_{\varphi}\vec{u}
	\end{pmatrix},\quad
	\partial_{\psi}x=\sigma\begin{pmatrix}
	-s\,\vec{u}\\
	c\,\vec{u}
	\end{pmatrix},
	\end{split}
	\end{equation}
	
	i.e.
	
	\begin{align}
	\label{new13}
	G_{AB}=\begin{pmatrix}
	1 &  & 0 \\ 
	& \sigma^2\begin{pmatrix}
	1 & 0 \\ 
	0 & \sin^2\theta
	\end{pmatrix}  &  \\ 
	0 &  & 1
	\end{pmatrix} \text{ and }G^{AB}=\begin{pmatrix}
	1 &  & 0 \\ 
	& \frac{1}{\sigma^2}\begin{pmatrix}
	1 & 0 \\ 
	0 & \frac{1}{\sin^2\theta}
	\end{pmatrix}  &  \\ 
	0 &  & 1
	\end{pmatrix} 
	\end{align}
	
	are diagonal. Hence only the diagonal entries of the second fundamental forms are needed, - which are all zero, due to 
	
	\begin{align}
	\label{new14}
	n=\begin{pmatrix}
	-s\partial_{\theta}\vec{u}\\
	c\partial_{\theta}\vec{u}
	\end{pmatrix},\quad n'=\begin{pmatrix}
	-s\partial_{\varphi}\vec{u}\\
	c\partial_{\varphi}\vec{u}
	\end{pmatrix}
	\end{align}
	
	both being orthogonal to 
	
	\begin{align}
	\label{new15}
	\partial_{\sigma}^2x=0,\quad \partial_{\theta}^2x=-x,\quad \partial_{\varphi}^2x=-\begin{pmatrix}
	c\,\sin\theta\begin{pmatrix}
	\cos\varphi \\ 
	\sin\varphi \\ 
	0
	\end{pmatrix}  \\ 
	s\,\sin\theta\begin{pmatrix}
	\cos\varphi \\ 
	\sin\varphi \\ 
	0
	\end{pmatrix}
	\end{pmatrix}, \quad\partial_{\psi}^2x=-x. 
	\end{align}

Higher dimensional cases are of course less trivial: already for the 10-dimensional space $\sum_{4,3}\subset\mathbb{R}^{12}$ one has to either deal with the dependencies in the SVD approach,

\begin{align}
\label{new16}
A=\sigma\,\vec{u}\,\overrightarrow{v}^T\,+\,\sigma'\,\vec{u}\,'\,\overrightarrow{v}\,'^T,
\end{align}

namely $\vec{u}\perp\vec{u}\,'$ and $\overrightarrow{v}\perp\overrightarrow{v}'$; or invert the metric

\begin{equation}
\label{new17}
\begin{split}
\hat{G}_{\hat{B}\hat{C}}=\begin{pmatrix}
\left(1+\lambda^2\right)\mathds{1} & \lambda\lambda'\mathds{1} & \lambda\vec{a} & \lambda\vec{a}' \\ 
\lambda\lambda'\mathds{1} & \left(1+\lambda'^2\right)\mathds{1} & \lambda'\vec{a} & \lambda'\vec{a}' \\ 
\lambda\vec{a}^T & \lambda'\vec{a}^T & \vec{a}^2 & \vec{a}\vec{a}' \\ 
\lambda\vec{a}'^T & \lambda'\vec{a}'^T & \vec{a}'\vec{a} & \vec{a}'^2
\end{pmatrix}=\begin{pmatrix}
g_{\beta\gamma} & u_{\beta} & u_{\beta}' \\ 
u^T_{\gamma} & \vec{a}^2 & \vec{a}\vec{a}' \\ 
u_{\gamma}^T & \vec{a}'\vec{a} & \vec{a}'^2
\end{pmatrix} =\begin{pmatrix}
G_{BC} & u'_B \\ 
u'^T_C & \vec{a}'^2
\end{pmatrix}
\end{split}
\end{equation}

which however, for $\vec{a},\vec{a}'\in\mathbb{R}^{\text{any }p>3}$, can be done, noting that

\begin{equation}
\label{new18}
\begin{split}
\mid\hat{G}\mid&=\mid G\mid\left( \vec{a}'^2-u'_BG^{BC}u'_{C}\right)=\frac{\mid G\mid}{\hat{\rho}} =\frac{g}{\rho\hat{\rho}}\\
&=\left( 1+\lambda^2+\lambda'^2\right) ^{p-2}\left( \vec{a}^2\vec{a}'^2-\left( \vec{a}\vec{a}'\right)^2\right), 
\end{split}
\end{equation}

and

\begin{align}
\label{new19}
G^{BC}=\begin{pmatrix}
g^{\beta\gamma}+\rho u^{\beta}u^{\gamma} & -\rho u^{\beta}\\
-\rho u^{T\gamma} & \rho
\end{pmatrix}, \quad \rho=\vec{a}^2-u_{\alpha}g^{\alpha\beta}u_{\beta}.
\end{align}

Letting $\vec{e},$ $\vec{e}\,', \dots \vec{e}\,^{(p-3)}$ be orthonormal vectors orthogonal to $\vec{a}$ and $\vec{a}'$, it is easy to see that an orthonormal basis of vectors normal to $\sum_{p,3}$ can be chosen to be 

\begin{align}
\label{new20}
n=\frac{1}{\mu^2}\begin{pmatrix}
\lambda\vec{e}\\
\lambda'\vec{e}\\
-\vec{e}
\end{pmatrix},\quad n'=\frac{1}{\mu^2}\begin{pmatrix}
\lambda\vec{e}'\\
\lambda'\vec{e}'\\
-\vec{e}'
\end{pmatrix}, \cdots,
\end{align}

(where $\mu^2:=(1+\lambda^2+\lambda'^2)$), and due to the $\partial^{2}_{1\lambda}x, \dots, \partial^{2}_{p\lambda'}x$ being the only non-vanishing second derivatives of 

\begin{align}
\label{new21}
x=\begin{pmatrix}
\vec{a} \\ 
\vec{a}' \\ 
\lambda\vec{a}+\lambda'\vec{a}'
\end{pmatrix} ,
\end{align}

one gets

\begin{align}
\label{new22}
-\mu^2\hat{H}_{\hat{A}\hat{B}}=\begin{pmatrix}
0_{2p\times 2p} & \vec{e} & \vec{0} \\ 
 & \vec{0} & \vec{e} \\ 
\begin{array}{cc}
\vec{e}^T & \vec{0}^T \\ 
\vec{0}^T & \vec{e}^T
\end{array}  && 0_{2\times 2}
\end{pmatrix} , \quad -\mu^2\hat{H}'_{\hat{A}\hat{B}}=\begin{pmatrix}
0 & \begin{array}{cc}
\vec{e}\,' & \vec{0} \\ 
\vec{0} & \vec{e}\,'
\end{array}  \\ 
& 0
\end{pmatrix} ,\cdots
\end{align}

for the second fundamental forms.
While $O(p)-$covariance suggests that the relevant parts of $\hat{G}^{\hat{A}\hat{B}}$ are linear combinations of $\vec{a}$ and $\vec{a}'$ (so that$\hat{G}^{\hat{A}\hat{B}}\hat{H}_{\hat{A}\hat{B}}=0=\hat{G}^{\hat{A}\hat{B}}\hat{H}'_{\hat{A}\hat{B}}\cdots$), one can also explicitly verify that, using (\ref{new17}-\ref{new19}). In the general case one presumably has that the determinant of the induced metric on $\sum_{p>q}$ is\footnote{I am grateful to M. Hynek for discussions and checking (\ref{new23}) in a few cases by symbolic computer calculations.}

\begin{align}
\label{new23}
\mid \hat{G}_{..}\mid=\left( 1+\lambda^2+\lambda'^2+\dots+\lambda^2_{q-1}\right) ^{p-q+1}\cdot\det\left( \vec{a}^{T}_{i}\vec{a}_{j}\right) _{i,j=1\dots q-1}.
\end{align}

\newpage
\secct{U(1)-Invariant Minimal 3-Manifolds}

While it is quite possible that even the general problem of extremal hypersurfaces in Euclidean- and Minkowski-space(s) is of integrable nature, the place to start is clearly the case of only 2 independent variables involved. Let me present several aspects of this interesting problem.

\bigskip\textbf{Orthonormal Parametrization}
\bigskip

Start with (cp. \cite{Hoppe95})

\begin{align}
\label{4.1}\ddot{\vec{x}}=\frac{1}{\rho}\partial_{A}\left( \frac{gg^{AB}}{\rho}\partial_{B}\vec{x}\right),
\end{align}

\begin{align}
\label{4.2}
\frac{d}{dt}\left( \rho:=\frac{\sqrt{g:=(\det\left(  g_{AB}=\partial_{A}\vec{x}\cdot\partial_{B}\vec{x})\right) }}{\sqrt{1-\dot{\vec{x}}^2}} \right) =0,
\end{align}

which are the minimal surface equations for time-like manifolds in Minkowski-space,

\begin{align*}
\frac{1}{\sqrt{G}}\partial_{\alpha}\left( \sqrt{G}G^{\alpha\beta}\partial_{\beta}x^{\mu}\right) =0, \quad \mu =0, \cdots, D-1, \  \alpha,\beta=0,\cdots, M,
\end{align*}

when choosing $x^0\left( \varphi^0 \varphi^1\cdots\varphi^{M}\right) = \varphi^0 =: t, $ i.e.

\begin{align}
\label{4.3}
x^{\mu}=\begin{pmatrix}
t \\ 
\vec{x}(t,\varphi^{1}\cdots\varphi^{M})
\end{pmatrix} ,
\end{align}

as well as 

\begin{align}
\label{4.4}
G_{0A}=\dot{x}^{\mu}\partial_{A}x^{\nu}\eta_{\mu\nu}=-\dot{\vec{x}}\partial_{A}\vec{x}\stackrel{!}{=}0, \quad A=1,\dots, M.
\end{align}

The Ansatz

\begin{align}
\label{4.5}
\vec{x}=\begin{pmatrix}
r_{1}(t,u) \vec{n}_{1}(\phi_{1}) \\ 
\vdots \\ 
r_{l}(t,u) \vec{n}_{l}(\phi_{l}) \\ 
\vec{z}(t,u)
\end{pmatrix} 
\end{align}

with the $\vec{n}_i(\phi_{i})$ parameterizing minimal embeddings into unit spheres, and $\vec{z}\in \mathbb{R}^{n}$, gives

\begin{align}
\label{4.6}
\left( G_{\alpha\beta}\right) = \begin{pmatrix}
1-\dot{\vec{r}}^{\, 2}-\dot{\vec{z}}^{\,2} &  &  \\ 
& -\vec{r}^{\, '2}-\vec{z}^{\, '2} &  \\ 
&  & - \begin{pmatrix}
r_{1}^{2}\left(\stackrel{1}{g}_{a_{1}b_{1}}\right) &  &  \\ 
& \ddots &  \\ 
&  & r_{l}^{2}\left(\stackrel{l}{g}_{a_{l}b_{l}}\right)
\end{pmatrix} 
\end{pmatrix} 
\end{align}

i.e.

\begin{align*}
g=\left( \vec{r}\ '^{2}+ \vec{z}\ '^{2} \right) \prod^{l}_{i=1}r_{i}^{2d_i}\stackrel{i}{g}=\left( \vec{r}\ '^2+\vec{z}\ '^2\right) \tilde{g}(t, u)\overline{g}(\phi_{1}\cdots \phi_{l})
\end{align*}

\begin{align}
\label{4.7}
g^{AB}=\begin{pmatrix}
\frac{1}{\vec{r}\ '^{2}+\vec{z}\ '^{2}} &  &  &  \\ 
& \frac{1}{r_{1}^{2}}\left( \stackrel{1}{g}\! ^{a_{1}b_{1}}\right)  &  &  \\ 
&  & \ddots &  \\ 
&  &  & \frac{1}{r_{l}^{2}}\left( \stackrel{l}{g}\!^{a_{l}b_{l}}\right) 
\end{pmatrix} 
\end{align}

\begin{align}
\label{4.8}
\rho=\frac{\sqrt{\vec{r}\ '^{2}+\vec{z}\ '^{2}}}{\sqrt{1-\dot{\vec{r}}^2-\dot{\vec{z}}^2}}\cdot\sqrt{\tilde{g}(t, u)}\cdot\sqrt{\overline{g}}=:\tilde{\rho}(u)\cdot\sqrt{\overline{g}(\phi)}
\end{align}

and (\ref{4.1}) reads

\begin{equation}
\begin{split}
\label{4.9}
\ddot{\vec{x}}&=\frac{1}{\rho}\partial_{u}\tilde{g}\frac{\sqrt{\overline{g}}}{\tilde{\rho}}\partial_{u}\vec{x}+\sum^{l}_{i=1}\frac{1}{\tilde{\rho}\sqrt{\overline{g}}}\partial_{a_{i}}\frac{\left( \vec{r}\ '^{2}+\vec{z}\ '^{2}\right) }{\tilde{\rho}\sqrt{\overline{g}}}\frac{\tilde{g}\overline{g}}{r_{i}^{2}}\stackrel{i}{g}\!^{a_{i}b_{i}}\partial_{b}\vec{x}\\
&=\frac{1}{\tilde{\rho}}\partial_{u}\frac{\tilde{g}}{\tilde{\rho}}\partial_{u}\vec{x}+\frac{\left( \vec{r}\ '^{2}+\vec{z}\ '^{2}\right) }{\tilde{\rho}^{2}}\tilde{g}\sum_{i}\frac{1}{r^2_i}\Delta_i \vec{x}.
\end{split}
\end{equation}

Reparameterizing the u-dependence (later dropping the $\sim$) via (cp. e.g. \cite{H2009})

\begin{align}
\label{4.10}
v(u):=\int\limits^{u}\tilde{\rho}(w)\ dw, \ r_{i}(u)=\tilde{r}_{i}\left(v(u)\right), \ \vec{z}(u)=\vec{\circledcirc{z}}(v(u))
\end{align}

(resp., if wanting $v$ to have a fixed range, absorbing the arising constant in scaling $t$) one gets (denoting $\frac{d}{dt}$ by $'$ and using $\Delta_i \vec{n}_i=-d_i\vec{n}_i$)

\begin{equation}
\label{4.11}
\begin{split}
\ddot{\vec{z}}&=\left( \vec{z}\, ' \prod r_{i}^{\, 2d_i}\right)',\\ \ddot{r}_{j}&= \left( r_j ' \prod r_{i}^{\, 2d_i}\right)'-\left(\vec{r}^{\,'2}+\vec{z}^{\,'2}\right) \frac{d_j}{ r_{j}}\prod r_{i}^{\, 2d_i},
\end{split}
\end{equation}

to be solved subject to

\begin{equation}
\label{4.12}
\begin{split}
\dot{\vec{r}}\vec{r}'+\dot{\vec{z}}\vec{z}'=0,\\
\dot{\vec{r}}^2+\dot{\vec{z}}^2+\tilde{g}\left( \vec{r}\,'^2+\vec{z}\,'^2\right) =\epsilon^2.
\end{split}
\end{equation}

Let us look at the simplest non-trivial case ($l=1=n, D=4,$ since \cite{H82} discussed several times; see
\cite{EHHS}, as well as various comments in \cite{Hoppe94}, \cite{9405001}, \cite{9407103}; \cite{WW} observed that (\ref{4.5}) must develop singularities in finite (positive, or negative) time due to (\ref{4.11}) trivially implying that the acceleration of the $v$-integral of the $r_i$'s is negative; singularity-formation similarity-equations were given already in \cite{EH2009}, and further analyzed in \cite{EHHS} - which also contains, for both spherical and toroidal topologies, pictures of numerical simulations, including singularity formations
at $r=0$ as well as $r\neq 0$ ),

\begin{align}
\label{4.13}
\vec{x}(t,\varphi,\theta)=\begin{pmatrix}
r(t,\varphi)\cos\theta \\ 
r(t,\varphi)\sin\theta \\ 
z(t,\varphi)
\end{pmatrix} ,
\end{align}

\begin{align}
\label{4.14}
\ddot{z}=\left( z'r^2\right) ', \ddot{r}=\left( r^2r'\right) '-(r'^2+z'^2)r,
\end{align}

subject to (and, as long as $\dot{\vec{x}}$ and $\vec{x}'$ are non-zero, actually implied by)

\begin{align}
\label{4.15}
\dot{r}r'+\dot{z}z'=0, \quad \dot{r}^2+\dot{z}^2+r^2\left(r'^2+z'^2\right)=\epsilon^2,
\end{align}

resp.

\begin{align}
\label{4.16}
\left( \dot{r}\pm rr'\right)^2+\left( \dot{z}\pm rz'\right)^2=\epsilon^2;
\end{align}

note that one could directly define

\begin{align}
\label{4.17}
\frac{1}{\epsilon}\left( rr'\pm \dot{r}\right) =:\pm\sin u_\pm, \quad\frac{1}{\epsilon}\left( rz'\pm \dot{z}\right) =:\cos u_\pm.
\end{align}

The hodograph-transformation $t\varphi\leftrightarrow r,z$ 

\begin{align}
\label{4.18}
\begin{pmatrix}
\dot{r}r' \\ 
\dot{z}z'
\end{pmatrix} = \frac{1}{t_{r}\varphi_{z}-t_{z}\varphi_{r}}\begin{pmatrix}
\varphi_{z} & -t_{z} \\ 
-\varphi_{r} & t_{r}
\end{pmatrix} 
\end{align}

transforms (\ref{4.15}) into

\begin{align}
\label{4.19}
\frac{1}{t_{r}^{2}+t_{z}^{2}}+\frac{r^{2}}{\varphi_{r}^{2}+\varphi_{z}^{2}}=\epsilon^{2} \quad \varphi_{z}t_{z}+\varphi_{r}t_{r}=0,
\end{align}

which can be solved by

\begin{equation}
\label{4.20}
\begin{split}
\begin{pmatrix}
\varphi_{z} \\ 
\varphi_{r}
\end{pmatrix} &=\frac{r}{\epsilon\cos w}\begin{pmatrix}
\cos v\\
\sin v
\end{pmatrix}=\frac{r \cosh u}{\epsilon}\begin{pmatrix}
\cos v\\
\sin v
\end{pmatrix}\\
\begin{pmatrix}
t_{z} \\ 
t_{r}
\end{pmatrix}&=\frac{\coth u}{\epsilon}\begin{pmatrix}
-\sin v\\
\cos v
\end{pmatrix}=\frac{1}{\epsilon \sin w}\begin{pmatrix}
-\sin v\\
\cos v
\end{pmatrix}.
\end{split}
\end{equation}

The compatibility conditions $\left( \varphi_{zr}=\varphi_{rz}, t_{zr}=t_{rz}\right) $ are

\begin{equation}
\label{4.21}
\begin{split}
\sin u_{+}\partial_{r}u_{-}+\cos u_{+}\partial_{z}u_{-}&=-\frac{1}{2r}(\cos u_{+}+\cos u_{-})\\
\sin u_{-}\partial_{r}u_{+}-\cos u_{-}\partial_{z}u_{+}&=-\frac{1}{2r}(\cos u_{+}+\cos u_{-})
\end{split}
\end{equation}

where $u_{\pm}:=\arctan\sinh\pm v=w\pm v$ (the quadrically non-linear equations (\ref{4.15}) have thus been transformed into quasi-linear equations). $u_{\pm}\leftrightarrow r,z$ on the other hand,

\begin{align}
\label{4.22}
\begin{pmatrix}
\partial_{r}u_{+} & \partial_{z}u_{+} \\ 
\partial_{r}u_{-} & \partial_{z}u_{-}
\end{pmatrix} =\left( \frac{\partial u_{+}u_{-}}{\partial rz}\right) = \left( \frac{\partial rz}{\partial u_{+}u_{-}}\right) ^{-1}=\frac{1}{r_{+}z_{-}-r_{-}z_{+}}\begin{pmatrix}
z_{-} &-r_{-}\\
-z_{+} & r_{+}
\end{pmatrix}
\end{align}

gives

\begin{equation}
\label{4.23}
\begin{split}
	-\sin u_{-}z_{+}+\cos u_{+}r_{+}&=-\frac{1}{2r}(\cos u_{+}+\cos u_{-})(r_{+}z_{-}-r_{-}z_{+})\\
	&=\sin u_{-}z_{-}+\cos u_{-}r_{-}\\
	&=\cdots,
\end{split}
\end{equation}

implying in particular

\begin{align}
\label{4.24}
\left( \sin u_{+}\partial_{+}+\sin u_{-}\partial_{-}\right) z=\left( \cos u_{+}\partial_{+}-\cos u_{-}\partial_{-}\right) r,
\end{align}

which when multiplied by $\tan\left( \frac{u_{+}+u_{-}}{2}\right) $ allows to be interpreted as 

\begin{align}
\label{4.25}
\frac{\partial}{\partial x_{2}}z=\frac{\partial}{\partial x_{1}}r,
\end{align}

where

\begin{align}
\label{4.26}
x_{1}:=\frac{\sin u_{+}-\sin u_{-}}{\sin u_{+}\cos u_{-}+\sin u_{-}\cos u_{+}}, \quad x_{2}:=\frac{\cos u_{+}-\cos u_{-}}{\sin u_{+}\cos u_{-}+\cos u_{+}\sin u_{-}},
\end{align}

satisfying 

\begin{align}
\label{4.27}
x_{1}\cos u_{\pm}\pm x_{2}\sin u_{\pm}=\mp 1 ,
\end{align}

\begin{align}
\label{4.28}
x^2\cos u_{\pm}=\mp x_{1}-x_{2}\sqrt{x^2-1}, \quad x^2\sin u_{\pm}=-x_2\pm x_1\sqrt{\sin^2-1}
\end{align}

\begin{equation}
\label{4.29}
\begin{split}
x^2\sin u_{+}\sin u_{-}&=1-x_1^2, \quad x^2\cos u_{+}\cos u_{-}=x_2^2-1\\
x^2\left( \sin \cos u_{-}-\cos u_{+}\sin u_{-}\right) &=-2x_1x_2,\quad x^2\left( \sin\cos u_{-}+\cos u_{+}\sin u_{-}\right)=2\sqrt{x^2-1}
\end{split}
\end{equation}

\begin{align}
\label{4.30}
\frac{\sin u_{+}+\sin u_{-}}{\cos u_{+}+\cos u_{-}}=\tan\left( \underbrace{\frac{u_{+}+u_{-}}{2}}_{=w}\right) =\frac{1-\cos 2w}{\sin 2w}=\frac{\sin 2w}{1+\cos 2w}=\frac{1}{\sqrt{x^2-1}}.
\end{align}

(\ref{4.23}) then implies the Monge-Ampère equation

\begin{align}
\label{4.31}
q_2\left[ \left( 1-x_1^2\right) q_{11}+\left( 1-x_2^2\right) q_{22}-2x_1x_2q_{12}\right] =+x_2\left( x_1^2+x_2^2-1\right) \left( q_{11}q_{22}-q_{12}^2\right) 
\end{align}

if putting $z=q_{1}, \ r=q_{2}$ (cp. (\ref{4.25})). The hodograph transformation $x_{1}x_{2}\leftrightarrow r,z$

\begin{align}
\label{4.32}
\begin{pmatrix}
r_1 & r_2 \\ 
z_1 & z_2
\end{pmatrix} =(r_{1}z_{2}-r_{2}z_{1})\begin{pmatrix}
\partial_{z}x_{2} & -\partial_{z}x_{1} \\ 
-\partial_{r}x_{2} & \partial_{r}x_{1}
\end{pmatrix} 
\end{align}

and using (\ref{4.25}) implying

\begin{align}
\label{4.33}
\partial_{z}x_2=\partial_r x_1
\end{align}

\begin{equation}
\label{4.34}
\begin{split}
\partial_{\pm}&=\frac{1}{fg}(\cos_{\mp}\partial_{2}\pm \sin_{\mp}\partial_{1})\\
r_{+}z_{-}-r_{-}z_{+}&=\frac{1}{f^2g}(r_{1}z_{2}-r_{2}z_{1})
\end{split}
\end{equation}

\begin{align*}
f:=\tan\frac{u_{+}+u_{-}}{2}=\frac{1}{\sqrt{x^2-1}}, \quad g:=\sin_{+}\cos_{-}+\cos_{+}\sin_{-}=\frac{2\sqrt{x^2-1}}{x^2}
\end{align*}

gives

\begin{align*}
\left( \frac{t^2_r+t^2_z-1}{r}\right) \left( t_{r}+r\Delta t\right) =t_r^2 t_{rr}+t^2_z t_{zz}+2t_{r}t_{z}t_{rz}
\end{align*}

resp.

\begin{align}
\label{4.35}
\frac{1}{r}\left( t^{2}_{r}+t^{2}_{z}-1\right) t_{r}-\Delta t=-t_{z}^{2}t_{rr}-t_{r}^{2}t_{zz}+2t_{r}t_{z}t_{rz}.
\end{align}

(when putting $x_{2}=-t_{r}, \ x_1=-t_z$, cp. (\ref{4.33})), which is the equation of motion following from the time-graph-description $t(r,z)$ of the minimal surface, with volume

\begin{align}
\label{4.36}
\int\sqrt{t_r^2+t_z^2-1}r drdz.
\end{align}

Note that (\ref{4.26}) also directly implies

\begin{align}
\label{4.37}
x_1=\frac{r'}{\dot{r}z'-r'\dot{z}}=-t_z, \quad x_2=\frac{-z'}{\dot{r}z'-r'\dot{z}}=-t_r
\end{align}

(using (\ref{4.17}), i.e. $\sin_{\pm}=\dot{r}\pm rr', \cos_{\pm}=rz'\pm \dot{z}$, and (\ref{4.18}).

\begin{align}
\label{4.38}
\begin{pmatrix}
\dot{r} \\ 
\dot{z}
\end{pmatrix} = -\sin\psi\begin{pmatrix}
\cos\phi\\\sin\phi
\end{pmatrix}, \ \begin{pmatrix}
r' \\ 
z'
\end{pmatrix}=\frac{\cos \psi}{r}\sin\psi\begin{pmatrix}
-\sin\phi \\
\cos\phi
\end{pmatrix}
\end{align}

gives

\begin{align}
\label{4.39}
\dot{\phi}=r\psi', \ \dot{\psi}=r\phi'+\frac{\cos\psi\cos\phi}{r},
\end{align}

resp.

\begin{equation}
\label{4.40}
\begin{split}
\dot{u}_{+}-ru_{+}'=-\frac{1}{2r}\left( c_{+}+c_{-}\right) \\
\dot{u}_{-}-ru_{-}'=-\frac{1}{2r}\left( c_{+}+c_{-}\right) ,
\end{split}
\end{equation}

with $u_{\pm}(t,\varphi):=-(\psi\pm\phi)$ (note that for strings one would get $\left( \partial_{t}\mp \partial_{\varphi}\right)u_{\pm}=0 $ instead, i.e. $u_{\pm}=f_{\pm}(\varphi\pm t)$). $\varphi^{0}=t,$ $\varphi\rightarrow t$, $r(t,\varphi)$ resp.

\begin{align}
\label{4.41}
\partial_{0}=\partial_{t}-\sin\psi\cos\phi\partial_{r}, \ r\partial_{\varphi}=-\cos\psi\sin\phi\partial_{r}
\end{align}

results in

\begin{equation}
\label{4.42}
\begin{split}
\dot{\tilde{u}}_{\pm}+s_{\mp}\underbrace{\tilde{u}_{\pm}'}_{=\partial_{r}\tilde{u}_{\pm}}=-\frac{1}{2r}(c_{+}+c_{-})
\end{split}
\end{equation}

for 

\begin{align*}
u_{\pm}(t,\varphi)=\tilde{u}_{\pm}(t,r(t,\varphi)).
\end{align*}

One can (try to explicitly) introduce characteristic coordinates either in each of the quasilinear systems, or directly in (\ref{4.16}), introducing $t_{+}$ and $t_{-}$ (instead of $t$ and $\varphi$) according to \footnote{Many thanks to J. Eggers for related discussions.}

\begin{align}
\label{4.43}
\frac{\partial t}{\partial t_{\pm}}\left( \partial_{t} \pm r\partial_{\varphi}\right) :=\frac{\partial}{\partial t_{\pm}},
\end{align}

obtaining

\begin{align}
\label{4.44}
\left( \partial_{\pm}r\right)^{2}+\left( \partial_{\pm}z\right) ^{2}=\left( \partial_{\pm}t\right) ^{2}\epsilon^{2}.
\end{align}

\begin{align}
\label{4.45}
\frac{\partial\varphi}{\partial t_{\pm}}=\pm r\frac{\partial t}{\partial t_{\pm}}
\end{align}

gives

\begin{align}
\label{4.46}
2\partial_{+-}^2t+\partial_{+}R\partial_{-}t+\partial_{-}R\partial_{+}t=0
\end{align}

which may be viewed as an inhomogenuous first order linear PDE for $R:=\ln r$ for given $t(t_{+}, t_{-})$, or as a homogenuous second order linear PDE for $t$ (given in terms of $R(t_{+},t_{-})$). The second compatibility condition is hidden in (\ref{4.44}), as if taken to determine

\begin{align}
\label{4.47}
\partial_{\pm}z=\zeta_{\pm}\sqrt{\epsilon^{2}(\partial_{\pm}t)^{2}-(\partial_{\pm}r)^{2}},
\end{align}

\begin{align}
\label{4.48}
t_{+-}\left( \frac{\partial_{+}}{\sqrt{+}}t-\tilde{\zeta}\frac{\partial_{-}}{\sqrt{-}}t\right)=r_{+-}\left( \frac{\partial_{+}}{\sqrt{+}}r-\tilde{\zeta}\frac{\partial_{-}}{\sqrt{-}}r\right)  
\end{align}

has to hold (with $\tilde{\zeta}=\frac{\zeta_{+}}{\zeta_{-}}$ being either plus or minus 1). Another possibility (linking to (\ref{4.21})) to display the second (hidden) integrability condition is to use (cp. (\ref{4.17}))

\begin{equation}
\begin{split}
\label{4.49}
\partial_{\pm}r&=\epsilon\sin u_{\pm}\partial_{\pm}t,\\
&=\epsilon s_{\pm}\partial_{\pm}t\\
\partial_{\pm}z&=\pm\epsilon\cos u_{\pm}\partial_{\pm}t,\\
&=\epsilon c_{\pm}\partial_{\pm}t
\end{split}
\end{equation}

giving

\begin{align*}
\partial_{-}\left( s_{+}\partial_{+}t\right) =\partial_{+}\left( s_{-}\partial_{-}t\right),\ \partial_{-}(c_{+}\partial_{+}t)=-\partial_{+}(c_{-}\partial_{-}t) 
\end{align*}

i.e.

\begin{equation}
\label{4.50}
\begin{split}
t_{+-}(s_{+}-s_{-})+(\partial_{+}t\partial_{-}s_{+}-\partial_{-}t\partial_{+}s_{-})&=0\\
t_{+-}(c_{+}+c_{-})+(\partial_{+}t\partial_{-}c_{+}+\partial_{-}t\partial_{+}c_{-})&=0\\
\end{split}
\end{equation}

resp.

\begin{equation}
\label{4.51}
\begin{split}
\partial_{+}u_{-}&=-\frac{1}{2r}(c_{+}+c_{-})\partial_{+}t\\
\partial_{-}u_{+}&=-\frac{1}{2r}(c_{+}+c_{-})\partial_{-}t
\end{split}
\end{equation}

with

\begin{align}
\label{4.52}
\partial_{+}t\partial_{-}u_{+}=\partial_{-}t\partial_{+}u_{-}=\rho=\frac{s_{+}-s_{-}}{c_{-}-c_{+}}t_{+-}=\frac{c_{+}+c_{-}}{s_{+}+s_{-}}t_{+-}
\end{align}

and 

\begin{align}
\label{4.53}
\frac{t_{+-}}{\partial_{+}t\partial_{-}t}=-\frac{1}{2r}(s_{+}+s_{-})\\
\partial_{\pm}u_{\mp}=-\frac{1}{2r}(c_{+}+c_{-})\frac{\partial_{\pm}r}{\epsilon s_{\pm}}.
\end{align}

Note that (\ref{4.45}) also implies

\begin{align}
\label{4.55}
2\varphi_{+-}=\partial_{+}\ln r\partial_{-}\varphi+\partial_{-}\ln r\partial_{+}\varphi.
\end{align}

\bigskip\textbf{Standard Characteristic Coordinates Approach}
\bigskip

Start with either of the graph-descriptions for $U(1)$-symmetric extremal hypersurfaces in $\mathbb{R}^{1,3}$:

\begin{align}
\label{5.56}
S\left[ t(r,z)\right] =\int r\sqrt{t_{r}^{2}+t_{z}^{2}-1}\ drdz,
\end{align}

stationary points of which are solutions of (cp. (\ref{4.35})) 

\begin{align}
\label{5.57}
t_{rr}(t_{z}^{2}-1)+t_{zz}(t_{r}^{2}-1)-2t_{r}t_{z}t_{rz}=-\frac{t_{r}}{r}(t_{r}^{2}+t^{2}_{z}-1),
\end{align}

or

\begin{align}
\label{5.58}
S\left[ r(t,z)\right] =\int r\sqrt{1-\dot{r}^{2}+r'^2}\ dt dz,
\end{align}

the space-time volume expressed in terms of the radius (of the axially symmetric 3 manifold) as a function of time $t$ and height $z$, or

\begin{align}
\label{5.59}
S\left[ z(t,r)\right] =\int r\sqrt{1-\dot{z}^{2}+z'^2}\ dt dr,
\end{align}

whose variation gives (cp. e.g. \cite{Hoppe94})

\begin{align}
\label{5.60}
\ddot{z}(1+z'^2)-z''(1-\dot{z}^2)-2\dot{z}z'\dot{z}'=\frac{z'}{r}(1-\dot{z}^2+z'^2).
\end{align}

The standard procedure to introduce characteristic coordinates, e.g. for (cp. (\ref{5.58}))

\begin{align}
\label{5.61}
L\left[ r\right] =A\ddot{r}+Cr''+2B\dot{r}'=\ddot{r}(1+r'^2)-r''(1-r^2)-2\dot{r}r'\dot{r}'=-\frac{1}{r}(1-\dot{r}^2+r'^2),
\end{align}

is to rewrite the equation in the new coordinates $\xi_{\pm}(t,z)$,

\begin{equation}
\label{5.62}
\begin{split}
\dot{r}&=\dot{\xi}_{+}r_{+}+\dot{\xi}_{-}r_{-}(=\dot{\xi}^{a}r_{a}), \quad
 r'=\xi'_+r_+ +\xi'_- r_-(=\xi'^a r_a), \\ 
\ddot{r}&=\dot{\xi}^a r_{ab}\dot{\xi}^b+\cdots, \quad  r''=\xi'^a r_{ab}\xi'^b+\cdots, \quad \dot{r}'=\dot{\xi}^a r_{ab}\xi'^b+\cdots,
\end{split}
\end{equation}

demanding that the coefficients of $r_{++}$ and $r_{--}$,

\begin{align}
\label{5.63}
\alpha_{\pm}=A\ddot{\xi}^2_{\pm}+2B\dot{\xi}_{\pm}\xi'_{\pm}+C\xi'^2_{\pm},
\end{align}

are zero. For $\lambda=\lambda_{\pm}=\frac{\dot{\xi}_{\pm}}{\xi'_{\pm}}$, i.e.

\begin{align}
\label{5.64}
\dot{\xi}_{\pm}-\lambda_{\pm}\xi'_{\pm}=0
\end{align}

this gives

\begin{align}
\label{5.65}
\lambda_{\pm}=-\frac{B\pm\sqrt{B^2-AC}}{A}=\frac{\dot{r}r'\pm\sqrt{1-\dot{r}^2+r'^2}}{1+r'^2}
\end{align}

and with

\begin{align}
\label{5.66}
\begin{pmatrix}
\dot{\xi}_{+} & \xi'_+ \\ 
\dot{\xi}_{-} & \xi'_-
\end{pmatrix} =\begin{pmatrix}
t_+ & t_- \\ 
z_+ & z_-
\end{pmatrix} ^{-1}=\frac{1}{t_+z_--t_-z_+}\begin{pmatrix}
z_- & -t_- \\ 
-z_+ & t_+
\end{pmatrix} 
\end{align}

(\ref{5.64}) implies

\begin{align}
\label{5.67}
z_-+\lambda_+t_-=0=z_++\lambda_-t_+,
\end{align}

from which it is easy to see that in the case of (\ref{5.61}) (equivalently (\ref{5.57}) and (\ref{5.60})) the characteristic coordinates are such that

\begin{align}
\label{5.68}
\left(x^{\mu}\right)_+=\begin{pmatrix}
t_+ \\ 
r_+ \\ 
z_+
\end{pmatrix} \text{ and } \left(x^{\mu}\right)_-=\begin{pmatrix}
t_- \\ 
r_- \\ 
z_-
\end{pmatrix}
\end{align}

are null-vectors (in $\mathbb{R}^{1,2}$), i.e.

\begin{align}
\label{5.69}
t_{\pm}^2=r_{\pm}^2+z_{\pm}^2
\end{align}

as seen already in the first approach (cp. (\ref{4.44})) while here:

\begin{equation}
\begin{split}
\label{5.70}
r_{\pm}^2&=\left(\dot{r}t_{\pm}+r'z_{\pm}\right)^2\\
&=\dot{r}^2t_{\pm}^2+r'^2z_{\pm}^2+2\dot{r}r't_{\pm}z_{\pm}\\
&=(\dot{r}^2-1)t^2_{\pm}+(r'^2+1)z_{\pm}^2+\cdots+t_{\pm}^2-z_{\pm}^2\\
&=\underbrace{t_{\pm}^2(r'^2+1)\left\lbrace \left( \frac{z_{\pm}}{t_{\pm}}\right)^2 +\frac{2\dot{r}r'}{r'^2+1}\frac{z_{\pm}}{t_{\pm}}+\frac{\dot{r}^2-1}{r'^2+1}\right\rbrace}_{=0,\text{ due to (\ref{5.67})/(\ref{5.63})=0}} +\left( t_{\pm}^2-z_{\pm}^2\right)\\
&=0+\left( t_{\pm}^2-z_{\pm}^2\right) .
\end{split}
\end{equation}

Defining

\begin{align}
\label{5.71}
m_{\mu}:=\begin{pmatrix}
r_+z_--r_-z_+ \\ 
z_+t_--z_-t_+ \\ 
t_+r_--t_-r_+
\end{pmatrix} =\begin{pmatrix}
\left\lbrace r,z\right\rbrace  \\ 
\left\lbrace z,t\right\rbrace  \\ 
\left\lbrace t,r\right\rbrace 
\end{pmatrix} 
\end{align}

\begin{align}
\label{5.72}
\left\lbrace f,g\right\rbrace :=f_+g_--f_-g_+=\epsilon^{ab}\partial_a f\partial_b g
\end{align}

and using (\ref{5.66}), one can easily verify that

\begin{align}
\label{5.73}
\begin{split}
\begin{array}{ll}
\partial_t=\frac{1}{m_1}\epsilon^{ab}z_{a}\partial_b=\frac{\left\lbrace z,\cdot\right\rbrace }{\left\lbrace z,t\right\rbrace },  & \partial_z=\frac{-1}{m_1}\epsilon^{ab}t_{a}\partial_b=\frac{\left\lbrace \cdot,t\right\rbrace }{\left\lbrace z,t\right\rbrace }, \\ 
\dot{r}=\frac{\left\lbrace z,r\right\rbrace }{\left\lbrace z,t\right\rbrace }=-\frac{m_0}{m_1}, & r'=\frac{\left\lbrace r,t\right\rbrace }{\left\lbrace z,t\right\rbrace }=-\frac{m_2}{m_1},
\end{array} 
\end{split}
\end{align}

as well as 

\begin{equation}
\begin{split}
\label{5.74}
m_1^2+m_2^2&=2t_+t_-(x_+\cdot x_-)\\
m_0^2-m_1^2&=2z_+z_-(x_+\cdot x_-)\\
m_0^2-m_2^2&=2r_+r_-(x_+\cdot x_-),
\end{split}
\end{equation}

trivially implying (taking the sum of the 3 equalities)

\begin{align}
\label{5.75}
-\left( m_0^2-m_1^2-m_2^2\right) =-m^2=(x_+\cdot x_-)^2;
\end{align}

also,

\begin{align}
\label{5.76}
m_1\cdot m_2=(x_+\cdot x_-)S_0, \quad m_2\cdot m_0=S_1(x_+\cdot x_-), \quad m_0\cdot m_1=S_2(x_+\cdot x_-);
\end{align}

with

\begin{align*}
S^{\mu}=\begin{pmatrix}
r_+z_-+r_-z_+ \\ 
z_+t_-+z_-t_+ \\ 
t_+r_-+t_-r_+
\end{pmatrix}; 
\end{align*}

indices are raised and lowered with $\eta_{\mu \nu}=\begin{pmatrix}
1 &  &  \\ 
& -1 &  \\ 
&  & -1
\end{pmatrix} $. Whereas the normal form of hyperbolic quasi-linear second-order differential equations in two variables is known  in full generality (cp.\cite{CH}\footnote{Many thanks to T. Damour.}; note however the absolute-sign ambiguity that would arise for (\ref{5.61}), and an existing explicit result for the 'free' Born Infeld equation (\cite{BS2}, eq. (8))), it is useful to explicitely derive from scratch the form (\ref{5.61}) will take in characteristic coordinates (\ref{5.69}), in 2 different ways: first,

\begin{equation}
\begin{split}
\label{5.77}
\ddot{r}&=\frac{1}{m_1^3}\left\lbrace\begin{array}{l}
m_1\left[ \vec{z}^{2}\Delta r+r_az_{ab}z_b-z_a r_{ab}z_b-\Delta z\vec{z}\vec{r}\right] \\ 
+m_0\left[ \vec{z}^{2}\Delta t+z_az_{ab}t_b-\vec{z}\vec{t}\Delta z-z_at_{ab}z_b\right]
\end{array} \right\rbrace \\
(1+r'^2)&=\frac{1}{m^2_1}(m_1^2+m_2^2)=\frac{2}{m_1^2}(x_+\cdot x_-)\\
r''&=\frac{1}{m_1^3}\left\lbrace\begin{array}{l}
m_1\left[ t_at_{ab}r_b+\vec{t}^2\Delta r-\Delta t\vec{t}\vec{r}-t_ar_{ab}t_b\right] \\ 
+m_2\left[\Delta z\vec{t}^2+z_at_{ab}t_b-t_az_{ab}t_b-\vec{t}\vec{z}\Delta t \right]
\end{array} \right\rbrace \\
\dot{r}^2-1&=\frac{1}{m_1^2}(m_0^2-m_1^2)=\frac{2}{m_1^2}z_+z_-(x_+\cdot x_-)\\
\dot{r}'&=\frac{1}{m_1^3}\left\lbrace\begin{array}{l}
m_1\left[\Delta t\vec{z}\vec{r}+t_ar_{ab}z_b-z_at_{ab}r_b-\vec{z}\vec{t}\Delta r \right] \\ 
+m_2\left[z_az_{ab}t_b+\vec{z}^{2}\Delta t-\Delta z\vec{z}\vec{t}-z_at_{ab}z_b\right]
\end{array} \right\rbrace\\
-2\dot{r}r'&=\frac{-2}{m_1^2}m_0m_2=\frac{2}{m_1^2}(x_+\cdot x_-)(z_+t_-+z_-t_+)
\end{split}
\end{equation}

While by construction (cp. (\ref{5.63})$_{=0}$) the terms proportional to $r_{++}$ and $r_{--}$ have to cancel when inserting (\ref{5.77}) into (\ref{5.61}), this is also true (clear by Lorentz invariance, and [CH]/[BC]) for the terms proportional to $t_{++},\  t_{--},\  z_{++} \text{ and } z_{--}$ (though tedious to verify when using (\ref{5.77})). A simple way to derive what (\ref{5.61}) becomes in the characteristics is as follows: (\ref{5.58}) gives

\begin{align}
\label{5.78}
\partial_{\alpha}\left( \frac{rr^{\alpha}}{\sqrt{1-r^{\gamma}r_{\gamma}}}\right) =\left( \frac{\dot{r}r}{\sqrt{1-\dot{r}^{2}+r'^{2}}}\right) -\left( \frac{rr'}{\sqrt{\ }}\right) '\stackrel{!}{=}-\sqrt{1-\dot{r}^2+r'^2},
\end{align}

hence

\begin{equation}
\label{5.79}
\begin{split}
-\frac{1}{r}&=\sqrt{\ }\partial_{\alpha}\left( \frac{r^{\alpha}}{\sqrt{\ }}\right) =\frac{x_+x_-}{m_1}\left( \left( \frac{m_1\dot{r}}{x_+x_-}\right)^{\cdot} -\left( \frac{m_1r'}{x_+x_-}\right)' \right) \\
&=\frac{-x_+x_-}{m_1^2}\left( \left\lbrace z,\frac{m_0}{x_+x_-} \right\rbrace +\left\lbrace t,\frac{m_2}{x_+x_-} \right\rbrace \right) \\
&=\frac{-1}{m_1^2}\left[ \begin{array}{l}
z_+(r_+z_--r_-z_+)_--z_-(r_+z_--r_-z_+)_+ \\ 
+t_+(t_+r_--t_-r_+)_--t_-(t_+r_--t_-r_+)_+
\end{array} \right] \\
&+\frac{1}{x_+\cdot x_-m_1^2}\left( m_0\left\lbrace z,x_+\cdot x_-\right\rbrace +m_2\left\lbrace t,x_+\cdot x_-\right\rbrace \right).
\end{split}
\end{equation}

The last bracket is easily seen to give

\begin{align}
\label{5.80}
-r_+(x_+\cdot x_-)_--r_-(x_+\cdot x_-)_+,
\end{align}

which due to (\ref{5.69}) implying

\begin{align}
\label{5.81}
x_{+-}\cdot x_{+}=0=x_{+-}\cdot x_{-}
\end{align}

does not contribute any terms with $x_{+-}$ components, while canceling all terms arising from the square bracket, except the mixed ones,

\begin{equation}
\label{5.82}
\begin{split}
2r_{+-}(z_+z_--t_+t_-)-z_{+-}(z_+r_-+z_-r_+)+t_{+-}(t_+r_-+t_-r_+)=-2r_{+-}(x_+\cdot x_-)
\end{split}
\end{equation}

(for the last step using (\ref{5.81})). Hence (using again the crucial orthonormality (\ref{5.81}))

\begin{align}
\label{5.83}
x^{\mu}_{+-}=\frac{m_1}{2r}n^{\mu}, \quad n^{\mu}:=\frac{m^{\mu}}{x_+\cdot x_-},
\end{align}

consistent with (\ref{4.46}), as (\ref{5.83}) implies

\begin{align}
\label{5.84}
t_{+-}=r_{+-}\frac{n^0}{n^1}=\frac{m_1}{2r}\left( \frac{r_+z_--r_-z_+}{x_+\cdot x_-}\right) =-\frac{1}{2r}(t_+r_-+t_-r_+).
\end{align}

Similarly,

\begin{align}
\label{5.85}
z_{+-}=-\frac{1}{2r}(r_+z_-+r_-z_+),
\end{align}

while $r_{+-}$ can also be expressed as

\begin{align}
\label{5.86}
r_{+-}=-\frac{1}{2r}\frac{(t_+r_-+t_-r_+)(r_+z_-+r_-z_+)}{(z_+t_-+z_-t_+)},
\end{align}

hence the completely symmetric equalities

\begin{align}
\label{5.87}
t_{+-}(r_+z_-+r_-z_+)=z_{+-}(t_+r_-+t_-r_+)=r_{+-}(z_+t_-+z_-t_+).
\end{align}

Such quasi-linear relations reconcile the result (\ref{5.83}) with \cite{CH}/\cite{BS2} expressing the 'free' Born-Infeld solution in terms of the Lorentz-invariant determinant

\begin{align}
\label{5.88}
D_{+-}=\begin{vmatrix}
t_{+-} & r_{+-} & z_{+-} \\ 
t_{+} & r_{+} & z_{+} \\ 
t_{-} & r_{-} & z_{-}
\end{vmatrix}=x^{\mu}_{+-}\cdot m_{\mu}=\frac{1}{2}\epsilon_{\mu\nu\rho}\cdot x^{\mu}_{+-}\cdot x^{\nu}_{+}x^{\rho}_{-} ,
\end{align}

resp. (32)$_{\cite{CH}}$ giving (for (\ref{5.61})) - note their strange absolute value -

\begin{equation}
\label{5.89}
\begin{split}
-D_{+-}&=(t_+z_--t_-z_+)^2\frac{1}{2r}\sqrt{1-\dot{r}^2+r'^2}\\
&=\frac{1}{2r}\left| (t_+z_--t_-z_+)(x_+\cdot x_-)\right| .
\end{split}
\end{equation}

Note the simplicity of (\ref{5.84})/(\ref{5.85}) compared to (\ref{5.89}).

\medskip

\bigskip\textbf{Lagrangian Null-Vector formulation}
\bigskip

The simplest derivation of (\ref{5.84})/(\ref{5.85}) consists of using the characteristic coordinates directly in the Lagrangian, yielding

\begin{align}
\label{6.90}
\int r(t_+t_--z_+z_--r_+r_-)\ d\theta_{+}d\theta_{-}.
\end{align}

Varying this action gives

\begin{equation}
\begin{split}
\label{6.91}
2rt_{+-}+r_{+}t_{-}+r_{-}t_{+}&=0\left( =\partial_{+}(rt_-)+\partial_{-}(rt_+)\right) \\
2rz_{+-}+r_{+}z_{-}+r_{-}z_{+}&=0\left( =\partial_{+}(rz_-)+\partial_{-}(rz_+)\right) \\
2rr_{+-}+r_{+}r_{-}+t_{-}t_{+}-z_{+}z_{-}&=0\left( =\partial_{+}(rr_-)+\partial_{-}(rr_+)+(x_+\cdot x_-)\right). 
\end{split}
\end{equation}

Note that the first 2 are linear in each of the unknown, and that each of the three equations implies the other two, once imposing/adding the quadratic constraints

\begin{align}
\label{6.92}
t^2_{\pm}-\left( r^2_{\pm}+z^2_{\pm}\right) =0,
\end{align}

as differentiating them gives

\begin{align}
\label{6.93}
t_{+-}\begin{pmatrix}
t_+ \\ 
t_-
\end{pmatrix}=r_{+-}\begin{pmatrix}
r_+ \\ 
r_-
\end{pmatrix}+z_{+-}\begin{pmatrix}
z_+ \\ 
z_-
\end{pmatrix},
\end{align}

which implies

\begin{equation}
\begin{split}
\label{6.94}
m_0r_{+-}+m_1t_{+-}&=0\\
m_0z_{+-}+m_2t_{+-}&=0\\
m_2r_{+-}&=m_1z_{+-},
\end{split}
\end{equation}

resp. (using $m_0S_0=m_1S_1=m_2S_2$, cp. (\ref{5.76})/(\ref{5.87}))

\begin{align}
\label{6.95}
S^0t_{+-}=S^1r_{+-}=S^2z_{+-}.
\end{align}

As $r$ appears, undifferentiated, in all 3 second-order equations it seems better/necessary to keep $r_{\pm}$ when solving the constraints (\ref{6.92}) explicitly. One possibility is to take

\begin{equation}
\begin{split}
\label{6.96}
t_{+}&=r_{+}\cosh p\quad z_{+}=r_{+}\sinh p\\
t_{-}&=r_{-}\cosh q\quad z_{-}=r_{-}\sinh p,
\end{split}
\end{equation}

the compatability conditions $t_{+-}=t_{-+}$, $z_{+-}=z_{-p}$ giving

\begin{equation}
\begin{split}
\label{6.97}
r_{+-}+\frac{1}{e^{p-q}-1}r_+p_-+\frac{1}{e^{q-p}-1}r_-q_+&=0\\
r_{+-}+\frac{1}{1-e^{q-p}}r_+p_-+\frac{1}{1-e^{p-q}}r_-q_+&=0,
\end{split}
\end{equation}

which implies

\begin{align}
\label{6.98}
r_+p_-+r_-q_+=0.
\end{align}

Note that (cp.(\ref{5.67}))

\begin{equation}
\begin{split}
\label{6.99}
p&=\text{arctanh}\frac{z_+}{t_+}=-\text{arctanh}\lambda_-\left(= \frac{1}{2}\ln \frac{t_++z_+}{t_--z_-}\right) \\
q&=\text{arctanh}\frac{z_-}{t_-}=-\text{arctanh}\lambda_+\left(= \frac{1}{2}\ln \frac{t_-+z_-}{t_--z_-}\right).
\end{split}
\end{equation}

Another way (as we will see, allowing one to reduce the equations) is to use light-cone variables (cp. \cite{H82})

\begin{align}
\label{6.100}
\tau:=\frac{t+z}{2}, \quad \xi:=t-z,
\end{align}

\begin{align}
\label{6.101}
t_+t_--z_+z_-=\tau_+\xi_-+\tau_-\xi_+
\end{align}

and explicitly solve the non-linear constraints (\ref{6.92}) as

\begin{align}
\label{6.102}
\xi_+=\frac{r^2_+}{2\tau_+}, \quad \xi_-=\frac{r^2_-}{2\tau_+},
\end{align}

the compatibility condition $\xi_{+-}=\xi_{-+}$ reading

\begin{align}
\label{6.103}
2r_{+-}\tau_+\tau_-=\tau_{+-}(r_+\tau_-+r_-\tau_+),
\end{align}

which together with (\ref{6.91})$_3$,

\begin{equation}
\label{6.104}
\begin{split}
\begin{tabular}{cc}
 & $2rr_{+-}+r_+r_-+r_-^2\frac{\tau_-}{2\tau_+}+r^2_+\frac{\tau_-}{2\tau_+}=0,$ \\ 
or (\ref{6.91})$_{1+2}$ &  \\  
& $2r\tau_{+-}+r_+\tau_-+r_-\tau_+=0,$  \\ 
\end{tabular}
\end{split}
\end{equation}

gives e.g.

\begin{equation}
\label{6.105}
\begin{split}
2r\tau_{+-}+(r_+\tau_-+r_-\tau_+)&=0\\
4\tau_+\tau_-rr_{+-}+(r_+\tau_-+r_-\tau_+)^2&=0,
\end{split}
\end{equation}

which consists of 2 (quasi-linear, second-order) equations for 2 (unconstrained) functions, $r$ and $\tau$, of 2 variables, the first equation linear in each $r$ and $\tau$, the second one homogeneously quadratic in each of the 2 unknown functions. Note that all $\pm $ equations (all the way, and explicit in (\ref{6.105})) are invariant under the change of variables (reparametrization) $\theta_+\rightarrow f_+(\theta_+)$ $\theta_-\rightarrow f_-(\theta_-)$ as well as scaling $r$ and $\tau$ (in (\ref{6.105}): independently; before: simultaneously $x^{\mu}\rightarrow\lambda x^{\mu}$). Also, finally, one could (at this stage) eliminate the appearance of the undifferentiated function $r$ by putting

\begin{align}
\label{6.106}
r=e^{-2\rho},
\end{align}

arriving then at (truly linear) 

\begin{align}
\label{6.107}
\tau_{+-}=\left( \rho_+\tau_-+\rho_-\tau_+\right). 
\end{align}

and

\begin{align*}
2\rho_{+-}=(\frac{\rho_+\tau_-+\rho_-\tau_+}{\tau_+\tau_-})^2+4\rho_+\rho_-.
\end{align*}

At this point it is perhaps worth mentioning that a particularly simple formulation for axially symmetric time-like extremal hypersurfaces in $\mathbb{R}^{1,3}$ was noticed already 40 years ago (\cite{H82}), but until now could not be shown to be integrable, namely (\cite{Hoppe94})

\begin{align}
\label{6.108}
\ddot{R}=R^2R''+RR'^2=R(RR')',
\end{align}

which is Hamilonian (\cite{H82}) with

\begin{align}
\label{6.109}
H=\frac{1}{2}\int(p^2+R^2R'^2)\ d\varphi.
\end{align}

The standard procedure concerning characteristics (cp. (\ref{5.62})-(\ref{5.64})) gives

\begin{align}
\label{6.110}
\lambda_{\pm}=\pm R, \quad \varphi _+-R\tau_+=0, \quad \varphi_{-}+R\tau_+=0,
\end{align}

\begin{align*}
\dot{R}&=\frac{\left\lbrace \varphi,R \right\rbrace }{\left\lbrace \varphi,\tau \right\rbrace }=\frac{1}{2\tau_+\tau_-}(\tau_+R_-+\tau_-R_+)\\
R'&=\frac{\left\lbrace R,\tau\right\rbrace }{\left\lbrace \varphi,\tau\right\rbrace }=\frac{R_+\tau_--R_-\tau_+}{2R\tau_+\tau_-}
\end{align*}

\begin{equation}
\label{6.111}
\begin{split}
\frac{1}{2}\int (\dot{R}^2-R^2R'^2)\ d\varphi d\tau=\frac{1}{2}\int\frac{\left\lbrace \varphi,R\right\rbrace ^2-R^2\left\lbrace R,\tau\right\rbrace ^2}{\left\lbrace\varphi,\tau \right\rbrace }\ d\theta_+d\theta_-
\end{split}
\end{equation}

($=\int RR_+R_- \ d\theta_+d\theta_-$ on solutions of (\ref{6.110})). Explicit use of (cp. (\ref{6.110})) \\ $\varphi_{\pm}=_{\pm}R\tau_{\pm}$, implying

\begin{equation}
\label{6.112}
\begin{split}
2\varphi_{+-}&=\tau_+R_--R_+\tau_-\\
-2R\tau_{+-}&=\tau_+R_-+R_+\tau_-
\end{split}
\end{equation}

in the equations of motion (\ref{6.108}), resp.

\begin{equation*}
\label{6.113}
\begin{split}
0&\stackrel{!}{=}\left\lbrace \varphi,\frac{\left\lbrace \varphi,R\right\rbrace }{\left\lbrace \varphi,\tau\right\rbrace }\right\rbrace -R\left\lbrace R\frac{\left\lbrace R,\tau\right\rbrace }{\left\lbrace \varphi,\tau\right\rbrace },\tau\right\rbrace \\
&=\frac{1}{\left\lbrace \varphi,\tau\right\rbrace ^2}\begin{pmatrix}
\left[\left\lbrace \varphi,\left\lbrace \varphi;R\right\rbrace \right\rbrace-R\left\lbrace R\left\lbrace R,\tau\right\rbrace ,\tau\right\rbrace  \right]\left\lbrace \varphi,\tau\right\rbrace \\ 
-\left\lbrace \varphi,R\right\rbrace \left\lbrace \varphi,\left\lbrace \varphi,\tau\right\rbrace \right\rbrace +R^2\left\lbrace R,\tau\right\rbrace \left\lbrace \left\lbrace \varphi,\tau\right\rbrace, \tau\right\rbrace 
\end{pmatrix} 
\end{split}
\end{equation*}

tediously gives

\begin{align}
\label{6.114}
4RR_{+-}\tau_+\tau_-+(R_+\tau_-+R_-\tau_+)^2=0,
\end{align}

hence fully coinciding with (\ref{6.105}); note that p.165 of \cite{CH} would give

\begin{align*}
-D_{+-}=-\left( R_{+-}\left\lbrace \varphi,\tau\right\rbrace +\varphi_{+-}\left\lbrace \tau,R\right\rbrace +\tau_{+-}\left\lbrace R,\varphi\right\rbrace \right) \stackrel{!}{=}-\frac{1}{2}\left\lbrace R,\tau\right\rbrace ^2,
\end{align*}

in this case without sign-ambiguity due to the radius R naturally being positive.

As noted already in \cite{H82}, the embedding function $\zeta=t-z$, needed to reconstruct the axially symmetric 3-manifold once $R$ is known, satisfies

\begin{align*}
\label{new114}
\zeta'=\dot{R}R',\quad \dot{\zeta}=\frac{1}{2}\left( \dot{R}^2+R^2R'^2\right) ,
\end{align*}

the compatibility being precisely (\ref{6.108}). In the characteristic coordinates, the above becomes

\begin{align*}
\left\lbrace \zeta,\tau\right\rbrace \left\lbrace \tau,\varphi\right\rbrace +\left\lbrace \varphi,R\right\rbrace \left\lbrace R,\tau\right\rbrace =0,\\
2\tau_+\tau_-(\tau_+\zeta_--\tau_-\zeta_+)=\tau_+^2R_-^2-\tau_-^2R_+^2
\end{align*}

and

\begin{align*}
2\left\lbrace \varphi,\zeta\right\rbrace \left\lbrace \varphi,\tau\right\rbrace &=\left\lbrace \varphi,R\right\rbrace ^2+R^2\left\lbrace R,\tau\right\rbrace ^2,\\
2\tau_+\tau_-(\tau_+\zeta_-+\tau_-\zeta_+)&=\tau^2_+R_-^2+\tau_-^2R_+^2,
\end{align*}

i.e.

\begin{align*}
4\tau_+^2\tau_-^2\zeta_-&=2\tau_+^2R_-^2\\
4\tau_-^2\tau_+^2\zeta_+&=2\tau_-^2R_+^2,
\end{align*}

resp.

\begin{align*}
\zeta_{\pm}=\frac{R_{\pm}}{2\tau_{\pm}},
\end{align*}

- the compatibility amounting to (cp. (\ref{6.103}))

\begin{align*}
2R_{+-}\tau_+\tau_-=\tau_{+-}(R_+\tau_-+R_-\tau_+).
\end{align*}

\noindent\hspace{10mm} Yet another way to derive (\ref{5.84})-(\ref{5.86}) resp. (\ref{6.105}) in characteristic coordinates (\ref{5.69}) is to start with (cp. \cite{EHHS} \footnote{Thanks to J.Eggers and M.Hynek for renewed discussions} , with $\varphi$ chosen as to make $\rho(\varphi)=\rho_{0}=r^2_0=\const$)

\begin{align}
\tag{$\ast$}
\dot{r}=-z'h, \quad \dot{z}=r'h, \quad h:=\sqrt{\frac{1}{z'^2+r'^2}-\frac{r^2}{r^2_0}}
\end{align}

(note that when solving the quadratic constraints, cp. (\ref{4.16}), this way, the opposite overall sign would do as well in ($\ast$); below this freedom, corresponding to $t\rightarrow -t$, or $\varphi \rightarrow -\varphi$, is used for consistency with previous formulae). First interchanging the role of dependent and independent coordinates (\ref{4.18}), obtaining

\begin{align*}
\varphi_{z}=-ht_r,\quad \varphi_{r}=ht_z,
\end{align*}

noting

\begin{align*}
h=\sqrt{\frac{(t_r\varphi_z-t_z\varphi_r)^2}{t^2_z+t_r^2}-\frac{r^2}{r_0^2}}=\sqrt{h(t^2_z+t_r^2)-\frac{r^2}{r_0^2}}
\end{align*}

implying

\begin{align*}
h=\frac{r}{r_0\sqrt{t_r^2+t_z^2-1}}
\end{align*}

and then changing to characteristic coordinates,

\begin{align}
\tag{$\sim$}
(r,z)\rightarrow(\theta_+,\theta_-), \quad t_{\pm}^2=r_{\pm}^2+z_{\pm}^2,
\end{align}

one obtains

\begin{equation*}
\begin{split}
r\left\lbrace t,z\right\rbrace \left\lbrace z,r\right\rbrace =r_0\left\lbrace \varphi,r\right\rbrace \sqrt{\left\lbrace t,r\right\rbrace ^2+\left\lbrace z,t\right\rbrace ^2-\left\lbrace r,z\right\rbrace ^2},\\
r\left\lbrace t,r\right\rbrace \left\lbrace r,z\right\rbrace =r_0\left\lbrace \varphi,z\right\rbrace \sqrt{\left\lbrace t,r\right\rbrace ^2+\left\lbrace z,t\right\rbrace ^2-\left\lbrace r,z\right\rbrace ^2},
\end{split}
\end{equation*}

which together with ($\sim$) provides 4 first order equations for the 4 unknown functions (of $\theta_{\pm}$) $\varphi, t, r, z$. In the previously introduced notation for the 3 Poisson-brackets of $t, r$ and $z$ the above reads

\begin{equation*}
\begin{split}
rm_1m_0=r_0\left\lbrace \varphi,r\right\rbrace \mid x_+x_-\mid\\
rm_2m_0=r_0\left\lbrace \varphi,z\right\rbrace \mid x_+x_-\mid
\end{split}
\end{equation*}

hence (using again (\ref{5.76}))

\begin{equation*}
\begin{split}
\left\lbrace \varphi,r\right\rbrace =\varphi_+r_--\varphi_-r_+=\pm \frac{r}{r_0}(t_+r_-+t_-r_+)\\
\left\lbrace \varphi,z\right\rbrace =\varphi_+z_--\varphi_-z_+=\pm \frac{r}{r_0}(z_+t_-+z_-t_+)
\end{split}
\end{equation*}

(the sign depending on whether $x_+x_-\lessgtr 0$), which (solving for $\varphi_{\pm}$) simply becomes

\begin{align*}
r_0\varphi_+=rt_+,\quad r_0\varphi_-=-rt_-
\end{align*}

(here, for consistency with (\ref{4.45}), having chosen the upper sign). Integrability ($\varphi_{+-}=\varphi_{-+}$) gives

\begin{align*}
2rt_{+-}+(t_+r_-+t_-r_+)=0,
\end{align*}

which together with the 2 non-linear constraints ($\sim$) seems to be the 'universal'/'simplest' formulation of the problem (apart from the purely first-order formulation ($\ast$) that was used in \cite{EHHS} to obtain numerical results). 

To verify in yet another way that the previous equation is the minimal surface equation for $U(1)$-symmetric 3 manifolds in $x_{\pm}^2=0$ coordinates, consider

\begin{align*}
\hat{x}^{\mu}(\theta_+,\theta_-)=\begin{pmatrix}
t \\ 
r\cos\theta \\ 
r\sin\theta \\ 
z
\end{pmatrix} ,
\end{align*}

implying

\begin{equation*}
\begin{split}
\hat{x}^{\mu}_{\pm}&=\begin{pmatrix}
t_{\pm} \\ 
r_{\pm}c \\ 
r_{\pm}s \\ 
z_{\pm}
\end{pmatrix}, \quad \partial_{\theta}\hat{x}^{\mu}=\begin{pmatrix}
0 \\ 
-sr \\ 
cr \\ 
0
\end{pmatrix}, \quad \hat{n}_{\mu}=\frac{1}{\hat{x}_+\hat{x}_-}\begin{pmatrix}
m_0 \\ 
m_1\stackrel{c}{s} \\ 
m_2
\end{pmatrix} =\frac{\hat{m}_{\mu}}{\hat{x}_+\hat{x}_-}\\
\hat{x}^{\mu}_{\pm\pm}&=\begin{pmatrix}
t_{\pm\pm} \\ 
r_{\pm\pm}\stackrel{c}{s} \\ 
z_{\pm\pm}
\end{pmatrix}, \quad \hat{x}^{\mu}_{\theta\theta}=-\begin{pmatrix}
0 \\ 
cr \\ 
sr \\ 
0
\end{pmatrix}, \quad \left( \hat{G}_{\hat{\alpha}\hat{\beta}}\right)=\begin{pmatrix}
0 & x_+x_- & 0 \\ 
x_+x_- & 0 & 0 \\ 
0 & 0 & -r^2
\end{pmatrix} \\
\left( \hat{G}^{\hat{\alpha}\hat{\beta}}\right)&=\frac{1}{\hat{x}_+\hat{x}_-}\begin{pmatrix}
0 & 1 & 0 \\ 
1 & 0 & 0 \\ 
0 & 0 & \frac{-\hat{x}_+\hat{x}_-}{r^2}
\end{pmatrix}, \quad\hat{H}_{\hat{\alpha}\hat{\beta}}=\begin{pmatrix}
\hat{x}_{++}\hat{n} & \hat{x}_{+-}\hat{n} & 0 \\ 
\hat{x}_{-+}\hat{n} & \hat{x}_{--}\hat{n} & 0 \\ 
0 & 0 & -\frac{m_1r}{\hat{x}_+\hat{x}_-}
\end{pmatrix} ,
\end{split}
\end{equation*}

hence

\begin{align*}
3H_3=G^{\hat{\alpha}\hat{\beta}}H_{\hat{\alpha}\hat{\beta}}=2\frac{\hat{x}_{+-}\hat{n}}{\hat{x}_{+}\hat{x}_{-}}+\frac{z_+t_--z_-t_+}{r(\hat{x}_{+}\hat{x}_{-})}=2H_2+H_1.
\end{align*}

For the mean curvature of the 2-manifold in $\mathbb{R}^{1,2}(t,r,z)$ one finds

\begin{align*}
H_2=\frac{\hat{x}_{+-}\hat{n}}{\hat{x}_{+}\hat{x}_{-}}=\frac{\hat{x}_{+-}\hat{m}}{(\hat{x}_{+}\hat{x}_{-})^2}=\frac{D_{+-}}{(x_+x_-)^2}=-\frac{t_{+-}}{m_0},
\end{align*}

as

\begin{align*}
D_{+-}=t_{+-}m_0+r_{+-}m_1+z_{+-}m_2=\frac{1}{m_0}t_{+-}\left( m_0^2-m_1^2-m_2^2\right) ;
\end{align*}

and then (cp. (\ref{5.76}))

\begin{align*}
-m_0r(3H_3=2H_2+H_1)=2rt_{+-}-\frac{m_0m_1}{x_+x_-}=2rt_{+-}+(r_+t_{-}+r_-t_+).
\end{align*}

\bigskip\textbf{Minimal Surfaces in Certain Non-Constant-Curvature Spaces}
\bigskip

As the action integral (\ref{6.90}) is the volume of a 2-dimensional parametrized surface $\sum=\left\lbrace x^{\mu}(\theta_{+},\theta_{-})\right\rbrace $ in a 3-dimensional Lorentzian space $\mathcal{L}$ with metric\footnote{Many thanks to M. Bordemann and M. Kontsevich for related discussions.} 

\begin{align}
\label{4.116}
\overline{\eta}_{\mu\nu}(x):=r\begin{pmatrix}
1 &  & 0 \\ 
& -1 &  \\ 
0 &  & -1
\end{pmatrix},
\end{align}

which on $\sum$ induces 

\begin{equation}
\label{4.117}
\begin{split}
\left( g_{\alpha\beta}\right)_{\alpha\beta=+-}=\left( x_{\alpha}\circ x_{\beta}:=\overline{\eta}_{\mu\nu}\partial_{\alpha}x^{\mu}\partial_{\beta}x^{\nu}\right) &=rx_+x_-\begin{pmatrix}
0 & 1 \\ 
1 & 0
\end{pmatrix} \\
&=\sqrt{-g}\begin{pmatrix}
0 & 1 \\ 
1 & 0
\end{pmatrix}=:e^u\begin{pmatrix}
0 & 1 \\ 
1 & 0
\end{pmatrix}
\end{split}
\end{equation}

when choosing coordinates $\theta_{\pm}$ on $\sum$ for which $x_{+}^{2}=0=x_{-}^{2}$, (\ref{6.91}) - which can also be written as

\begin{align}
\label{4.118}
2rx_{+-}^{\mu}+r_{+}x_{-}^{\mu}+r_{-}x_{+}^{\mu}=-x_{+}x_{-}\begin{pmatrix}
0 \\ 
1 \\ 
0
\end{pmatrix} ,
\end{align}

implying (via multiplication with $x_{+}$ resp. $x_{-}$)

\begin{align}
\label{4.119}
\partial_{+}\left( rx^2_-=x_-\circ x_-\right) =0=\partial_{-}\left( rx^{2}_{+}=x_+\circ x_+\right) 
\end{align}

(hence upon reparameterizing $\theta_{\pm}\rightarrow\tilde{\theta}_{\pm}(\theta_{\pm})$ $x_-\circ x_-=\const, x_+\circ x_+=\const$), must simply be the minimal surface equations in $\mathcal{L}$,

\begin{align}
\label{4.120}
\frac{1}{\sqrt{g}}\partial_{\alpha}\left( \sqrt{g}g^{\alpha\beta}\partial_{\beta}x^{\mu}\right)+g^{\alpha\beta}\partial_{\alpha}x^r\partial_{\beta}x^{\rho}\Gamma^{\mu}_{r\rho}=0,
\end{align}

where the Christoffel-Symbols (following from (\ref{4.116})) are

\begin{align}
\label{4.121}
\Gamma^{0}_{10}=\Gamma^{0}_{01}=\Gamma^{2}_{12}=\Gamma^{2}_{21}=\Gamma^{1}_{11}=\Gamma^{1}_{00}=-\Gamma^{1}_{22}=\frac{1}{2r}
\end{align}

(all others $=0$), while those following from (\ref{4.117}) are

\begin{align}
\label{4.122}
\gamma_{++}^{+}=u_{+}, \ \gamma_{--}^{-}=u_{-}
\end{align}

(all others $=0$). Indeed, multiplying (\ref{4.118}) by

\begin{align}
\label{4.123}
\overline{n}_{\mu}:=\frac{\overline{\epsilon}_{\mu\nu\rho}x_+^\nu x_-^{\rho}}{x_+\circ x_-}:=\frac{r^{\frac{3}{2}}\epsilon_{\mu\nu\rho}x_+^{\nu}x_-^{\rho}}{rx_+x_-}=\sqrt{r}\frac{m_{\mu}}{x_+x_-}=\sqrt{r}n_{\mu}
\end{align}

(hence $\overline{n}^2=\overline{n}_{\mu}\overline{n}_{\nu}\overline{\eta}^{\mu\nu}=rn_{\mu}n_{\nu}\frac{1}{r}\eta^{\mu\nu}=n^2=-1$) results in 

\begin{align}
\label{4.124}
g^{ab}\left(h_{ab}:=\overline{n}_{\mu}x_{ab}^{\mu}+\overline{n}_{\mu}\Gamma_{r\rho}^{\mu}\partial_{\alpha}x^{r}\partial_{\beta}x^{\rho}\right) \left( =\frac{2h_{+-}}{x_{+}\circ x_{-}}\right), 
\end{align}

which \underline{is} zero, when using (\ref{6.91}):

\begin{equation}
\begin{split}
\label{4.125}
\frac{1}{\sqrt{r}}\overline{n}_{\mu}x^{\mu}_{+-}&=n_0t_{+-}+n_1r_{+-}+n_2z_{+-}\\
\frac{1}{\sqrt{r}}\overline{n}_{\mu}\Gamma_{\nu\rho}^{\mu}\partial_{a}x^{\nu}\partial_{b}x^{\rho}&=\frac{1}{2r}\left\lbrace n_{0}(t_{+}r_{-}+t_{-}r_{+})\right. \\
&+n_{1}(t_{+}t_{-}+r_{+}r_{-}-z_{+}z_{-})+n_{2}(r_{+}z_{-}+r_{-}z_{+})\left. \right\rbrace .
\end{split}
\end{equation}

The GCMP (Gauss-Peterson-Codazzi-Mainardi) equations follow from (now switching to Latin indices abc for $\alpha\beta\gamma$, to avoid confusion with $\gamma^{\cdot}_{\cdot\cdot}$)

\begin{align}
\label{4.126}
\partial_{ab}^{2}x^{\mu}=\gamma_{ab}^{c}\partial_{c}x^{\mu}-h_{ab}\overline{n}^{\mu}-\partial_{a}x^\nu\partial_{b}x^{\rho}\Gamma_{\nu\rho}^{\mu}
\end{align}

by further differentiation ($\partial_{d},$ and then demanding/using that the $3^{rd}$ derivatives are completely symmetric, e.g. under $b \leftrightarrow d$); splitting the last term in (\ref{4.126}) without which e.g. the $\gamma^{c}_{ab}$ no longer be given by $\gamma_{ab}^{c}=\frac{1}{2}g^{cd}(\partial_{a}g_{db}+\partial_{b}g_{ad}-\partial_{d}g_{ab})$) into tangential and normal part one gets

\begin{equation}
\label{4.127}
\begin{split}
\partial^{2}_{ab}x^{\mu}&=\left( \gamma_{ab}^{c}-Y_{ab}^{c}\right)\partial_{c}x^{\mu}-(h_{ab}-L_{ab})\overline{n}^{\mu}\\
&=:\tilde{\gamma}_{ab}^{c}\partial_{c}x^{\mu}-\tilde{h}_{ab}\overline{n}^{\mu}
\end{split}
\end{equation}

with

\begin{equation}
\label{4.128}
\begin{split}
Y_{ab}^{c}&:=g^{ce}\partial_{a}x^{\nu}\partial_{b}x^{\rho}\partial_{e}x^{\lambda}\overline{\eta}_{\mu\lambda}\Gamma_{\nu\rho}^{\mu}\\
L_{ab}&:=\overline{n}_{\mu}\Gamma_{\nu\rho}^{\mu}\partial_{a}x^{\nu}\partial_{b}x^{\rho}.
\end{split}
\end{equation}

As $\partial_{a}\left( \overline{n}^{\mu}\overline{n}^{\nu}\overline{\eta}_{\mu\nu}\right) =0$ implies

\begin{align}
\label{4.129}
\overline{n}\circ\partial_{a}\overline{n}\left( =\overline{n}^{\mu}\,\overline{\eta}_{\mu\nu}\partial_{a}\overline{n}^{\nu}\right)=-\frac{1}{2}\overline{n}^{\mu}\overline{n}^{\nu}\partial_{a}\overline{\eta}_{\mu\nu} 
\end{align}

and $\partial_{a}(\overline{n}^{\mu}\partial_{c}x^{\lambda}\eta_{\mu\lambda})=0$ 

\begin{align}
\label{4.130}
\partial_{a}\overline{n}^{\mu}\overline{\eta}_{\mu\lambda}\partial_{c}x^{\lambda}=-\overline{n}^{\mu}\left( \partial^{2}_{ac}x^{\lambda}\right) \overline{\eta}_{\mu\lambda}-\overline{n}^{\mu}\partial_{c}x^{\lambda}\left(\partial_{\rho}\overline{\eta}_{\mu\lambda} \right)\partial_{a}x^{\rho} 
\end{align}

(the last term being zero for the metric (\ref{4.116}) as the derivative of the metric is proportional to the metric itself) one has

\begin{align}
\label{4.131}
\partial_{a}\overline{n}^{\mu}=-\tilde{h}_{ac}g^{cb}\partial_{b}x^{\mu}+\frac{1}{2}\left( \overline{n}^{\nu}\overline{n}^{\rho}\partial_{a}\overline{\eta}_{\nu\rho}\right)\overline{n}^{\mu}. 
\end{align}

Note that (\ref{4.126}/\ref{4.127}) and (\ref{4.131}) do\footnote{I thank J. Arnlind and T. Turgut for helpful discussions on this} coincide with the standard decomposition (cp. e.g. (2.1)-(2.5)of \cite{AHH}) 

\begin{equation}
\label{4.132}
\begin{split}
\overline{\nabla}_{a}\partial_{b}x^{\mu}&\left( =\partial_{ab}^{2}x^{\mu}+\Gamma^{\mu}_{\nu\rho}\partial_{a}x^{\nu}\partial_{b}x^{\rho}\right)\\
&=\nabla_{a}\partial_{b}x^{\mu}\left( =\gamma_{ab}^{c}\partial_{c}x^{\mu}\right) -h_{ab}n^{\mu}
\end{split}
\end{equation}

\begin{equation}
\label{4.133}
\begin{split}
\overline{\nabla}_{a}\overline{n}^{\mu}&\left( =\partial_{a}\overline{n}^{\mu}+\Gamma^{\mu}_{\nu\rho}\partial_{a}x^{\nu}n^{\rho}\right)\\
&=-h_{ab}g^{bc}\partial_{c}x^{\mu}+D_an^{\mu},
\end{split}
\end{equation}

due to (\ref{4.127}/\ref{4.128}) and

\begin{align}
\label{4.134}
\Gamma_{\nu\rho}^{\mu}\partial_{a}x^{\nu}n^{\rho}=Y_{a}^{c}\partial_cx^{\mu}+D_{c}n^{\mu},
\end{align}

with indeed (using that (\ref{4.116}) implies $\partial\overline{\eta}_{\mu\nu}=\frac{\partial r}{r}\overline{\eta}_{\mu\nu}$)

\begin{align}
\label{4.135}
Y_{a}^{c}\left( =\overline{\eta}_{\mu\lambda}\Gamma_{\nu\rho}^{\mu}\partial_{a}x^{\nu}n^{\rho}\partial_{b}x^{\lambda}g^{bc}\right) =-L_{ab}g^{bc}.
\end{align}

Differentiating (\ref{4.127}) with respect to $\theta_{d}$ and then using again (\ref{4.127}), as well as (\ref{4.131}), one finds that

\begin{equation}
\label{4.136}
\begin{split}
\partial_{d}\tilde{h}_{ab}+\tilde{\gamma}_{ab}^{c}\tilde{h}_{cd}-\partial_{b}\tilde{h}_{ad}-\tilde{\gamma}_{ad}^{c}\tilde{h}_{cb}=\frac{1}{2r}\left( r_{d}\tilde{h}_{ab}-r_{b}\tilde{h}_{ad}\right) 
\end{split}
\end{equation}

(the lhs of these CMP equations, $\nabla_{d}\tilde{h}_{ab}-\nabla_{b}\tilde{h}_{ad},$ being general, i.e. of canonical form, while the rhs is particularly simple, due to $\overline{\eta}_{\mu\nu}=r\eta_{\mu\nu}$ being conformally flat, i.e. differing from the flat Minkowski-metric only by a factor, $r$) and 

\begin{align}
\label{4.137}
-\left(\tilde{h}_{ab}\tilde{h}_{de}-\tilde{h}_{ad}\tilde{h}_{be}\right) =g_{ec}\left\lbrace \partial_{d}\tilde{\gamma}_{ab}^{c}-\partial_{b}\tilde{\gamma}_{ad}^{c}+\tilde{\gamma}_{ab}^{e}\tilde{\gamma}_{ed}^{c}-\tilde{\gamma}_{ad}^{e}\tilde{\gamma}_{eb}^{c}\right\rbrace 
\end{align}

(note the - sign on the lhs, which is due to $\overline{n}^2=-1$). Due to being antisymmetric under $b\leftrightarrow d$, as well as w.r.t. ($e\leftrightarrow a$), the Gauss-equations (\ref{4.137}) consists only of one single equation (while (\ref{4.136}), antisymmetric only under $b \leftrightarrow d$, leaves \underline{2} choices for $a$). Calculating the $Y_{ab}^{c}$ (crucially using $x_+^2=0=x_-^2$) one finds

\begin{align}
\label{4.138}
Y^{+}_{++}=\frac{r_+}{r}, \ Y^-_{--}=\frac{r_-}{r},
\end{align}

with all others (just as for the $\gamma_{ab}^{c}$, cp. (\ref{4.122})) vanishing; hence

\begin{align}
\label{4.139}
\tilde{\gamma}_{\pm\pm}^{\pm}\left( =u_{\pm}-\frac{r_{\pm}}{r}\right) =\left( \ln (x_+x_-)\right) _{\pm}.
\end{align}

(\ref{4.137}) then gives

\begin{align}
\label{4.140}
\tilde{h}\left( :=\tilde{h}_{11}\tilde{h}_{22}-\tilde{h}_{12}^2\right) =(x_+x_-)(\ln x_+x_-)_{+-},
\end{align}

while (\ref{4.136}), due to most of the $\tilde{\gamma}_{ab}^{c}$ vanishing, becomes

\begin{align*}
\partial_{-}\tilde{h}_{++}-\partial_{+}\tilde{h}_{+-}+\tilde{\gamma}_{++}+\tilde{h}_{+-}=\frac{r_-}{2r}\tilde{h}_{++}-\frac{r_+}{2r}\tilde{h}_{+-}
\end{align*}

i.e.

\begin{equation}
\label{4.141}
\begin{split}
\partial_{-}\left(\frac{\tilde{h}_{++}}{\sqrt{r}}\right)-\partial_{+}\left( \frac{\tilde{h}_{+-}}{\sqrt{r}}\right) +\tilde{\gamma}^{+}_{++}\frac{\tilde{h}_{+-}}{\sqrt{r}}&=0 \\
\partial_{-}\left(\frac{\tilde{h}_{+-}}{\sqrt{r}}\right)-\partial_{+}\left( \frac{\tilde{h}_{--}}{\sqrt{r}}\right) +\tilde{\gamma}^{-}_{--}\frac{\tilde{h}_{+-}}{\sqrt{r}}&=0,
\end{split}
\end{equation}

resp.

\begin{align*}
\partial_{-}\tilde{h}_{+-}-\partial_{+}\tilde{h}_{--}-\tilde{\gamma}^{-}_{--}\tilde{h}_{-+}=\frac{1}{2r}\left( r_{-}\tilde{h}_{+-}-r_{+}\tilde{h}_{--}\right). 
\end{align*}

While this derivation was following a standard route for the Gauss-Codazzi-Mainardi-Peterson consistency equations, a much simpler way is to \underline{not} use the curved metric $\overline{\eta}_{\mu\nu}(x)$, but simply writing, apart from noting $x_+^2=0=x_-^2$,

\begin{equation}
\label{4.142}
\begin{split}
x_{++}&=ax_++bn\\
x_{--}&=cx_-+dn\\
x_{+-}&=en
\end{split}
\end{equation}

where, by definition,

\begin{align}
\label{4.143}
\left( H_{\alpha\beta}\right) :=n\cdot\begin{pmatrix}
x_{++} & x_{+-} \\ 
x_{-+} & x_{--}
\end{pmatrix} =\left( n\cdot x_{\alpha\beta} \right) =\begin{pmatrix}
-b & -e \\ 
-e & -d
\end{pmatrix} 
\end{align}

\begin{align*}
a=\frac{(x_+x_-)_+}{x_+x_-}=(\ln w)_+, \quad c=\frac{(x_+x_-)_-}{x_+x_-}=\frac{w_-}{w}.
\end{align*}

Using the 'flat-space' Weingarten equations (which simply follow from $nx_{\pm}=0$)

\begin{equation}
\label{4.144}
\begin{split}
n_+&=\frac{e}{x_+x_-}x_++\frac{b}{x_+x_-}x_-\\
n_-&=\frac{d}{x_+x_-}x_++\frac{e}{x_+x_-}x_-,
\end{split}
\end{equation}

$(x_{++})_-=(x_{+-})_+$ and $(x_{--})_+=(x_{+-})_-$ trivially give

\begin{equation}
\label{4.145}
\begin{split}
e_+-ae&=b_-,\quad e_--ce=d_+\\
c_+&=\frac{e^2-bd}{x_+x_-}=a_-
\end{split}
\end{equation}

(which are easily seen to coincide with (\ref{4.140}/\ref{4.141}), as $\tilde{h}_{\alpha\beta}=\sqrt{r}H_{\alpha\beta}$.) The minimality condition $h_{12}=\tilde{h}_{12}+L_{12}\stackrel{!}{=}0$ gives (cp. (\ref{5.83}) (\ref{5.74}), (\ref{5.76}), (\ref{4.125}))

\begin{equation}
\label{4.146}
\begin{split}
e&=-\frac{\tilde{h}_{12}}{\sqrt{r}}=+\frac{L_{12}}{\sqrt{r}}\\
&=n_{\mu}\Gamma^{\mu}_{\nu\rho}\partial_{+}x^{\nu}\partial_{-}x^{\rho}\\
&=\frac{1}{2r}\left\lbrace n_0(t_+r_-+t_-r_+)+n_2(r_+z_-+r_-z_+)\right. \\
&\left.+n_1(t_+t_-+r_+r_--z_+z_-) \right\rbrace \\
&=\frac{m_1}{2r}=\frac{z_+t_--z_-t_+}{2r}.
\end{split}
\end{equation}

\bigskip\textbf{Zero Curvature Condition}
\bigskip

Using the null-vectors $x_{\pm}$ and the unit normal $n(n^2=-1)$ one may form a $\theta_{\pm}$-dependent $3\times 3$ matrix

\begin{align}
\label{4.147}
M:=\begin{pmatrix}
\frac{\lambda x_++\frac{1}{\lambda}x_-}{\sqrt{2w}}, & \frac{\lambda x_+-\frac{1}{\lambda}x_-}{\sqrt{2w}}, & n
\end{pmatrix} ^{\text{T}},
\end{align}

which is an element of 

\begin{align}
\label{4.148}
O(1,2):=\left\lbrace S\mid\left\langle Sx,Sy\right\rangle =\left\langle x,y\right\rangle (:=\eta_{\mu\nu}x^{\mu}y^{\nu}) \forall x,y\in\mathbb{R}^{1,2}\right\rbrace
\end{align}

(resp. $G=SO(1,2)$, when choosing an appropriate orientation for $n$); it satisfies the linear system of equations (cf. (\ref{4.142}/\ref{4.144}))

\begin{align}
\label{4.149}
\partial_{+}M=AM, \quad \partial_{-}M=BM,
\end{align}

with 

\begin{equation}
\label{4.150}
\begin{split}
A&=\frac{\lambda b}{\sqrt{2w}}(T_{1}-T_{0})+\frac{a}{2}T_{2}+\frac{e}{\lambda\sqrt{2w}}(T_{1}+T_{0})\\
B&=\frac{\lambda e}{\sqrt{2w}}(T_{1}-T_{0})-\frac{c}{2}T_{2}+\frac{d}{\lambda\sqrt{2w}}(T_{1}+T_{0}),
\end{split}
\end{equation}

$\lambda$ a nonzero constant, and

\begin{equation}
\label{4.151}
\begin{split}
T_{1}:=\begin{pmatrix}
0 & 0 & 1 \\ 
0 & 0 & 0 \\ 
1 & 0 & 0
\end{pmatrix} , \quad T_{2}:=\begin{pmatrix}
0 & 1 & 0 \\ 
1 & 0 & 0 \\ 
0 & 0 & 0
\end{pmatrix} , \quad T_{0}:=\begin{pmatrix}
0 & 0 & 0 \\ 
0 & 0 & -1 \\ 
0 & 1 & 0
\end{pmatrix} 
\end{split}
\end{equation}

satisfying $so(1,2)$ commutation-relations

\begin{align*}
\left[ T_1,T_2\right] =T_0, \quad \left[ T_2,T_0\right] =-T_{1}, \quad \left[ T_0,T_1\right] =-T_{2},
\end{align*}

resp.

\begin{align}
\label{4.152}
\left[ T_2,T_1\pm T_0\right]&=\mp (T_1\pm T_0)\\
\left[ T_1-T_0,T_1 + T_0\right]&=2T_2.\nonumber
\end{align}

One can easily check that the zero-curvature conditions implied by (\ref{4.149}),

\begin{equation}
\label{4.153}
\begin{split}
0&=\partial_{-+}^2M-\partial_{+-}^2M=\partial_{-}(AM)-\partial_{+}(BM)\\
&=A_-M+A(\partial_{-}M=BM)-B_+M-B(\partial_{+}M=AM)\\
&=(A_--B_++\left[ A,B\right] )M
\end{split}
\end{equation}

is equivalent to the GCMP equations (\ref{4.145}). One would like to take $\lambda$ as a spectral parameter, but (obvious from (\ref{4.152})) any constant $\lambda$ can be absorbed into a redefinition of $T_{\pm}:=T_1\mp T_0$, as $\lambda T_+,\frac{1}{\lambda}T_-,T_2$ satisfy (\ref{4.152}) just as well. One then may try to use that the equations to be solved (\ref{6.91}/\ref{6.92}) are invariant under reparametrizations,

\begin{align}
\label{4.154}
\theta_{\pm}\rightarrow f_{\pm}(\theta_{\pm})=\tilde{\theta}_{\pm},
\end{align}

\begin{equation}
\label{4.155}
\begin{split}
v_{\pm}&=\frac{\lambda f'_+\tilde{x}_{+}\pm \frac{f'_-}{\lambda}\tilde{x}_-}{\sqrt{f'_+f'_-}\sqrt{2}\tilde{w}}\\
&=\frac{1}{\sqrt{2\tilde{w}}}\left( \lambda\sqrt{\frac{f'_{+}}{f'_-}}\tilde{x}_+\pm \frac{1}{\lambda}\sqrt{\frac{f'_-}{f'_+}}\tilde{x}_-\right), 
\end{split}
\end{equation}

when $x(\theta_+,\theta_-)=\tilde{x}(f_+(\theta_+),f_-(\theta_-)),$ one is led to replace the constant $\lambda$ by a $\theta-$dependent function $\tilde{\lambda}(\theta_+,\theta_-)$, which gives (leaving out all $\sim$'s)

\begin{equation}
\label{4.156}
\begin{split}
A&=\frac{\lambda b}{\sqrt{2w}}T_++\left( \frac{a}{2}+\frac{\lambda_+}{\lambda}\right) T_2+\frac{e}{\lambda\sqrt{2w}}T_-\\
B&=\frac{\lambda e}{\sqrt{2w}}T_+-\left( \frac{c}{2}-\frac{\lambda_-}{\lambda}\right) T_2+\frac{d}{\lambda\sqrt{2w}}T_-.
\end{split}
\end{equation}

While (\ref{4.156}) now seems to genuinely depend on a ($\theta$-dependent) spectral parameter, the equations resulting from $M_+=AM,$ $M_-=BM$, resp (\ref{4.153}) do \underline{not}, as in

\begin{equation}
\label{4.157}
\begin{split}
0&=A_--B_++\left[ A,B\right] \\
&=\left( \frac{\lambda b}{\sqrt{2w}}\right) _-T_++\left( \frac{a_-}{2}+(\ln\lambda)_{+-}\right) T_2+\left( \frac{e}{\lambda\sqrt{2w}}\right) _-T_-\\
&-\left( \frac{\lambda e}{\sqrt{2w}}\right) _+T_++\left( \frac{c_+}{2}-(\ln\lambda)_{+-}\right) T_2-\left( \frac{d}{\lambda\sqrt{2w}}\right) _+T_-\\
&+\frac{\lambda b}{\sqrt{2w}}\left( \frac{c}{2}-\frac{\lambda_-}{\lambda}\right) T_++\frac{db}{w}T_2-\frac{d}{\sqrt{2w}\lambda}\left( \frac{a}{2}+\frac{\lambda_+}{\lambda}\right) T_-\\
&+\frac{\lambda e}{\sqrt{2w}}\left( \frac{a}{2}+\frac{\lambda_+}{\lambda}\right) T_+-\frac{e^2}{w}T_2-\frac{e}{\lambda\sqrt{2w}}\left( c-\frac{\lambda_-}{\lambda}\right) T_-
\end{split}
\end{equation}

all terms containing derivatives of $\lambda(\theta_+\theta_-)$ cancel! 

Using the fundamental representation of $sl(2,\mathbb{R})$, i.e.

\begin{equation}
\label{new156}
\begin{split}
T_1-T_0&\rightarrow E_+=\begin{pmatrix}
0 & 1 \\ 
0 & 0
\end{pmatrix} ,\quad T_1+T_0\rightarrow E_-=\begin{pmatrix}
0 & 0 \\ 
1 & 0
\end{pmatrix}\\
2T_2&\rightarrow H=\begin{pmatrix}
1 & 0 \\ 
0 & -1
\end{pmatrix} ,
\end{split}
\end{equation}

(\ref{4.156}) reads

\begin{align}
\label{new157}
A=\begin{pmatrix}
\frac{a}{4}+\frac{\lambda_+}{2\lambda} & \frac{\lambda b}{\sqrt{2w}} \\ 
\frac{e}{\lambda\sqrt{2w}} & -\frac{a}{4}-\frac{\lambda_+}{2\lambda}
\end{pmatrix}, \quad B=\begin{pmatrix}
-\frac{c}{4}+\frac{\lambda_-}{2\lambda} & \frac{\lambda e}{\sqrt{2w}} \\ 
\frac{d}{\lambda\sqrt{2w}} & \frac{c}{4}-\frac{\lambda_-}{2\lambda}
\end{pmatrix}. 
\end{align}

Unfortunately and somewhat unexpectedly, e.g. in contrast with the Lund-Regge model  \cite{LR76}  (where a constant spectral parameter related to residual reparametrization invariance is nontrivial, as in SO(3) it can not be absorbed by rescaling the generators ) one again finds that the $\lambda-$dependence is a pure gauge, i.e. can be gauged away, as if $M_+=AM,$ $M_-=BM,$ $ \tilde{M}:=e^{\frac{\sigma_3}{2}\ln \lambda}M$ will satisfy

\begin{equation}
\label{new158}
\begin{split}
\tilde{M}_+&=e^{-\frac{\sigma_3}{2}\ln\lambda}\left( -\frac{\lambda_+}{2\lambda}\sigma_3+A\right) M\\
&=\left( e^{-\frac{\sigma_3}{2}\ln\lambda}\left( A-\frac{\lambda_+}{2\lambda}\sigma_3\right) e^{+\frac{\sigma_3}{2}\ln\lambda}\right) \tilde{M}=\tilde{A}\tilde{M},\\
\left( e^{-\frac{\sigma_3}{2}\ln\lambda}M\right)_- &=\cdots=\tilde{B}\left( e^{-\frac{\sigma_3}{2}\ln\lambda}M\right) ,
\end{split}
\end{equation}

with

\begin{align}
\label{new159}
\tilde{A}=\begin{pmatrix}
\frac{1}{\sqrt{\lambda}} & 0 \\ 
0 & \sqrt{\lambda}
\end{pmatrix} 
\begin{pmatrix}
\frac{a}{4} & \frac{\lambda b}{\sqrt{2w}} \\ 
\frac{e}{\lambda\sqrt{2w}} & -\frac{a}{4}
\end{pmatrix} 
\begin{pmatrix}
 \sqrt{\lambda}& 0 \\ 
0 & \frac{1}{\sqrt{\lambda}}
\end{pmatrix} \text{ and } \tilde{B}=\cdots
\end{align}



\newpage
\textbf{Acknowledgement:}\\
Parts of these notes have been presented at AEI Potsdam, Bo\u{g}azi\c{c}i University, Braunschweig University, Copenhagen University, ETH, IHES, the Istanbul Center for Mathematical Sciences, the Korea Institute for Advanced Study, KTH, and Sogang University.
I am grateful for kind hospitality and useful comments that I received, and to S.Karimi-Dehbokri, M.Sandhaus, and C.Gourgues for helping with the manuscript.

\end{document}